# Orthonormal Basis in Minkowski Space

Aleks Kleyn

Alexandre Laugier

*E-mail address*: Aleks_Kleyn@MailAPS.org

*E-mail address*: Laugier.Alexandre@orange.fr

Abstract. Finsler space is differentiable manifold for which Minkowski space is the fiber of the tangent bundle. To understand structure of the reference frame in Finsler space, we need to understand the structure of orthonormal basis in Minkowski space.

In this paper, we considered the definition of orthonormal basis in Minkowski space, the structure of metric tensor relative to orthonormal basis, procedure of orthogonalization. Linear transformation of Minkowski space mapping at least one orthonormal basis into orthonormal basis is called motion. The set of motions of Minkowski space $\overline{V}$ generates not complete group $SO(\overline{V})$ which acts single transitive on the basis manifold.

Passive transformation of Minkowski space mapping at least one orthonormal basis into orthonormal basis is called quasimotion of Minkowski space. The set of passive transformations of Minkowski space generates passive representation of not complete group $SO(\overline{V})$ on basis manifold. Since twin representations (active and passive) of not complete group $SO(\overline{V})$ on basis manifold are single transitive, then we may consider definition of geometric object.

# Contents

CHAPTER 1

# Preface

## 1.1. Structure of the Set of Motions of Minkowski Space

This book appeared as a result of cooperation of authors on some of the issues that emerged as a result of the analysis of paper [8]. Euclidean space is the limiting case of Minkowski space, when the metric tensor does not depend on direction. However, although the matrix of the metric tensor with respect to orthogonal basis of Euclidean space is diagonal, the analog (2.3.9) of this matrix in Minkowski space is triangular matrix. This statement is extremely unsatisfactory.

Since the theorem 2.3.8 does not restrict value of entries of bottom triangle of the matrix (2.3.9), then we must assume that this value is arbitrary. This leads us to the statement that basis manifold of Minkowski space is larger than basis manifold of Euclidean space. When we finished calculations for the section 3.1, we realized that although this value is unknown, but it is not arbitrary.

The opportunity to simplify the equation [8]-(4.3) turned out to be more significant. The equation (2.4.8) resembles similar equation for the Euclidean space. From the equation (2.4.2), it follows that a matrix of a motion of Minkowski space is close to orthogonal matrix.

Since the set $SO(\overline{V})$ of motions of Minkowski space of dimension $n$ is subset of the group $GL(n)$ of linear transformation, then we may consider the product of motions. However the set $SO(\overline{V})$ is not close relative to product.[1.1] In other words, the set $SO(\overline{V})$ is not complete group.[1.2] At present it is difficult to predict the physical consequences of this statement. However, since not complete group $SO(\overline{V})$ has single transitive representation on basis manifold, and twin representation is also single transitive, then we may consider definition of geometric object.

## 1.2. Measurement of angle

In Euclidean space there are two ways of measuring the angle.

Consider angle generated by rays $a$ and $b$ that have a common endpoint called vertex of angle.

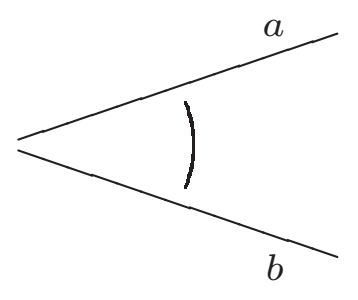

[1.1] We expect the paper where it is shown that in Minkowski space with norm (3.4.1) there exist motions whose product is not motion.

[1.2] Not complete group is particular case of not complete $\Omega$-algebra.

DEFINITION 1.1.1. Since operation $\omega \in \Omega$ is not defined for any tuple of elements of the set $A$, then the structure of **not complete $\Omega$-algebra** is defined on the set $A$. □

For instance, let associative product with inverse element and unit is defined on the set $A$. Thus, for any $a \in A$ there exists $b \in A$, such that product $ab$ is defined. However there exist $a$, $b \in A$, such that product is not defined. Then the set $A$ is called **not complete group**.

Consider a unit circle centered at the vertex of the angle. We can define the value of the angle as length of the arc contained between rays $a$ and $b$. This method of measuring of the angle is used in different areas of human activities including Astronomy (the section [1]-1.9).

When we perform analytical calculations it is convenient to define cosine of the angle using the scalar product of vectors that are parallel rays $a$ and $b$. We assume that measurements were performed in pseudo-Euclidean space, where scalar product is symmetric. Therefore, the results of these methods of measurement are the same.

Observer may choose any method of measurement of intervals of length and time. However, observers must agree on the standards of length and time in order to exchange measurement results.

Therefore, each observer chooses its basis in such a way that the length of vectors of basis is equal to the length of the standard. We will call corresponding basis normal. Each basis is associated by the coordinate system arranged so that the coordinate $x^i$ of points on the axis of coordinates $X^i$ is the distance from the origin. In Euclidean space, we identify every point and radius vector of this point. The length of this vector is the distance from origin to the considered point. Square of the length of this vector is a quadratic form of coordinates. In case of Euclidean space, this quadratic form is positively defined.

Observer can choose any normal basis. Accordingly, quadratic form defining the metric can be arbitrary form. However, we can reduce any quadratic form to canonical form ([3], p. 169 - 172). Corresponding basis is orthogonal basis. Orthogonal normal basis is called orthonormal basis.

Thus, observer uses orthonormal basis to measure spatial and temporal intervals.[1.3]

In Minkowski space, metric tensor is not quadratic form. The procedure of experimental determination of the metric for the given basis is more difficult task. By measuring the distance in different directions, we can construct an analog of the unit sphere in Euclidean space, namely surface

$$|F(\overline{v})| = 1$$

which is called indicatrix. In the paper [10], Asanov considered how, based on certain physical hypotheses, we can determine the form of Finsler metric.

According to the remark 2.2.4, scalar product in Minkowski space is not commutative. Hence the cosine of the angle depends on the direction of measurement. Therefore, in Minkowski space, the results of discussed above methods of measurement do not match. Correspondingly, we may consider two different definitions of orthonormal basis.

According to [15], vector $\overline{w}$ is orthogonal to vector $\overline{v}$, if vector $\overline{w}$ is tangent to the hypersurface

$$F^2(\overline{x}) = F^2(\overline{v}) \tag{1.2.1}$$

homothetic to indicatrix. From the equation (2.2.9), it follows that vector $\overline{w}$ which is tangent to the hypersurface (1.2.1) holds the equation

$$g_{ij}(\overline{v})v^i w^j = 0 \tag{1.2.2}$$

From the remark 2.2.4 it follows that vector $\overline{w}$ is orthogonal to vector $\overline{v}$, if scalar product

$$\overline{v} \cdot \overline{w} = 0$$

In this paper, we consider the definition 2.3.2 of orthogonal basis based on the definition (1.2.2) of orthogonality.

However, there are other definitions of orthogonality, for instance [11, 13].

[1.3] This step in reasoning seems like logical leap. Observing the free movement of solid and distracting from friction forces, Galileo discovers a law of inertia. When we move from a flat to a curved space, it is naturally to suggests motion along a geodesic instead of moving along a straight line. Request of orthonormality like Lagrange principle does not follow directly from the experiment. It seems surprising that the experiment confirms this view.

## 1.3. Measuring the Speed of Light

Observer in general relativity uses local basis in tangent space which is event space of special relativity. So the transition of observer from one local basis to another can be described by a Lorentz transformation.

When geometry of event space changes so dramatically that transformation of basis changes its structure, some authors accept this change as Lorentz invariance violation ([14]). If we see in the experiment, that is a deviation in the structure of the Lorentz transformation, then it is an argument in favor that we should consider the new geometry. However, if the new geometry is correct, we should not restrict ourselves by Lorentz invariance. Our task is to understand what kind of algebra is generated by transformations of new geometry and how it affects the invariance principle.

However there exists another concept of Lorentz invariance violation. This concept refers to the dependence of the speed of light on the direction ([12], p. 7, [13], p. 12). Standards of length and time and the speed of light are inextricably linked in the theory of relativity. If we know two of these parameters, we can determine the third. The choice of geometry is also important to determine the relationship between considered parameters.

Light propagates in the Finsler space. However, since the operation of speed detection is the local operation, then we may consider this operation in the Minkowski space. We assume that vector $\overline{e}_0$ is timelike vector, and the remaining vectors of the basis are spacelike vectors.

To measure speed of light we need to assume that we specified standards of length and time. Let Finsler metric of event space have form ([12]-(3))

$$F^2(\overline{v})=-c^2(v^0)^2+F_3^2(v^1,v^2,v^3) \tag{1.3.1}$$

where $c$ is factor binding standards of length and time and $F_3^2$ is positively defined Finsler metric of 3-dimensional space. To find speed of light in direction in space $(v^1,v^2,v^3)$, observer puts in a given direction segment whose length is equal to the standard of length

$$F_3(v^1,v^2,v^3)=1$$

and determines the time of propagation of light signal. Since the vector of propagation of light is isotropic, then from the equation (1.3.1) it follows that

$$-c^2(v^0)^2+F_3^2(v^1,v^2,v^3)=0 \tag{1.3.2}$$

For any selected direction in space $(v^1,v^2,v^3)$

$$v^0=\frac{F_3(v^1,v^2,v^3)}{c}=\frac{1}{c} \tag{1.3.3}$$

Therefore, for a given distance, time of propagation of light $v^0$ does not depend on direction in space.

The speed of light in considered model does not depend on direction and equal $c$. In Finsler space, Finsler metric depends on point of manifold. So the procedure of measurement of speed of light should be done in infinitesimal area.

An attempt to understand how the interaction of particles and quantum fields can change the geometry of event space is one of the reasons to study Finsler geometry ([10, 12]). Without doubt, only experiment can confirm whether the geometry of the event space is Finsler geometry and whether the speed of light depends on the direction.[1.4] Our task is to find an effective tool for experimental verification of the nature of geometry.

[1.4] Physicists conduct experiments for testing the isotropy of the speed of light ([4, 5, 6]). The goal of experiments is to obtain limits on Lorentz invariance violation.

Even if it turns out that the speed of light depends on the direction, this does not contradict to the statement that speed of light in given direction is maximal speed of signal transmission in given direction.

CHAPTER 2

# Minkowski Space

## 2.1. Homogeneous Function

In this paper, we consider vector space over real field $R$.

DEFINITION 2.1.1. Let $\overline{V}$ be vector space. Function $f(\overline{x})$, $\overline{x} \in \overline{V}$, is called **homogeneous of degree** $k$, if

$$f(a\overline{x}) = a^k f(\overline{x})$$

□

THEOREM 2.1.2 (Euler theorem). *Function $f(\overline{x})$, homogeneous of degree $k$, satisfies the differential equation*

$$\frac{\partial f(\overline{x})}{\partial x^i} x^i = k f(\overline{x}) \tag{2.1.1}$$

PROOF. We differentiate the equation (2.1.1) with respect to $a$

$$\frac{df(a\overline{x})}{da} = \frac{da^k}{da} f(\overline{x}) \tag{2.1.2}$$

According to chain rule, we get

$$\frac{df(a\overline{x})}{da} = \frac{\partial f(a\overline{x})}{\partial a x^i} \frac{d a x^i}{da} = \frac{\partial f(a\overline{x})}{\partial a x^i} x^i \tag{2.1.3}$$

From equations (2.1.2), (2.1.3), it follows that

$$\frac{\partial f(a\overline{x})}{\partial a x^i} x^i = k a^{k-1} f(\overline{x}) \tag{2.1.4}$$

Equation (2.1.1) follows from equations (2.1.4) if we assume $a = 1$. □

THEOREM 2.1.3. *If $f(\overline{x})$ is homogeneous function of degree $k$, then partial derivatives $\dfrac{\partial f(\overline{x})}{\partial x^i}$ are homogeneous functions of degree $k-1$.*

PROOF. Consider equation[2.1]

$$f(a\overline{x}) = a^k f(\overline{x}) \tag{2.1.5}$$

We differentiate the equation (2.1.5) with respect to $x^i$

$$a \frac{\partial f(a\overline{x})}{\partial x^i} = a^k \frac{\partial f(\overline{x})}{\partial x^i} \tag{2.1.6}$$

From equation (2.1.6), it follows that

$$\frac{\partial f(a\overline{x})}{\partial x^i} = a^{k-1} \frac{\partial f(\overline{x})}{\partial x^i} \tag{2.1.7}$$

[2.1] There is similar proof in [17], p. 265.

Therefore, derivatives $\frac{\partial f(\overline{x})}{\partial x^i}$ are homogeneous functions of degree $k-1$. □

If mapping

$$f : R \to R$$

is homogeneous of degree 0, then $f(x) = \text{const}$. In general, we can only state that

$$f(t\overline{x}) = f(\overline{x}) \quad t \in R$$

We will also say that the mapping $f$ is constant in the direction of $\overline{x}$.

## 2.2. Finsler Space

In this section, I made definitions similar to definitions in [12].

Definition 2.2.1. A vector space $V$ is called **Minkowski space**[2.2] when, for the vector space $V$, we define a **Finsler structure**

$$F : \overline{V} \to R$$

such that

2.2.1.1: The **Finsler metric** $F^2$ is not necessarily positive definite[2.3]

2.2.1.2: Function $F(\overline{x})$ is homogeneous of degree 1

$$F(a\overline{x}) = aF(\overline{x}) \quad a > 0 \tag{2.2.1}$$

2.2.1.3: Let $\overline{\overline{e}}$ be the basis of vector space $A$. Coordinates of **metric tensor**

$$g_{ij}(\overline{v}) = \frac{1}{2}\frac{\partial^2 F^2(\overline{v})}{\partial v^i \partial v^j} \tag{2.2.2}$$

form a nonsingular symmetric matrix. □

Let $\overline{V}$, $\overline{W}$ be Minkowski space of the same dimension. Let

$$f : \overline{V} \to \overline{W}$$

be nonsingular linear mapping. We accept as general rule ([2], p. 176), that if either

$$F_V^2(\overline{v}) = F_W^2(f(\overline{v}))$$

or

$$F_V^2(\overline{v}) = -F_W^2(f(\overline{v}))$$

where $F_V$ is Finsler structure in Minkowski space $\overline{V}$, $F_W$ is Finsler structure in Minkowski space $\overline{W}$, then we will not consider appropriate geometries significantly different and will examine only one of them. In particular, if $F(\overline{v}) < 0$ for any vector $\overline{v}$, then, without loss of generality, we may consider Finsler structure

$$F'(\overline{v}) = -F(\overline{v})$$

instead of Finsler structure $F$.

Definition 2.2.2. If $F(\overline{v}) = 0$, then vector $\overline{v}$ is called **isotropic vector**. If $F^2(\overline{v}) < 0$, then vector $\overline{v}$ is called **timelike vector**. In this case, $F(\overline{v})$ is imaginary positive number. If $F^2(\overline{v}) > 0$, then vector $\overline{v}$ is called **spacelike vector**. In this case, $F(\overline{v})$ is real positive number. □

[2.2] I considered the definition of Minkowski space according to the definition in [15], [16], p. 44. Although this term calls some association with special relativity, usually it is clear from the context which geometry is referred.

[2.3] This requirement is due to the fact that we consider applications in general relativity.

Definition 2.2.3. For any vector $\overline{v} \in \overline{V}$, $F(\overline{v}) \neq 0$, vector

$$\overline{u} = \frac{1}{|F(\overline{v})|}\overline{v}$$

is called **unit vector** □

Remark 2.2.4. According to the definition (2.2.2), metric tensor depends on vector $\overline{v}$. Therefore, the scalar product defined by the equation

$$\overline{v} \cdot \overline{w} = g_{ij}(\overline{v}) v^i w^j$$

in general, is not commutative. □

Let $\overline{\overline{e}}$ be the basis of Minkowski space. Since a vector $\overline{v}$ has expansion

$$\overline{v} = v^i \overline{e}_i$$

relative to the basis $\overline{\overline{e}}$, then we can present the mapping $F$ as

$$F(v^1, ..., v^n) = F(\overline{v})$$

Theorem 2.2.5. *The Finsler structure of Minkowski space satisfies the differential equations*

$$\frac{\partial F(\overline{a})}{\partial a^i} a^i = F(\overline{a}) \tag{2.2.3}$$

$$\frac{\partial^2 F(\overline{a})}{\partial a^i \partial a^j} a^i = 0 \tag{2.2.4}$$

$$\frac{1}{2}\frac{\partial^2 F^2(\overline{a})}{\partial a^i \partial a^j} a^i a^j = F^2(\overline{a}) \tag{2.2.5}$$

Proof. The equation (2.2.3) follows from the statement (2.2.1.2) of the definition 2.2.1 and the theorem 2.1.2. According to the theorem 2.1.2, derivative $\dfrac{\partial F(\overline{x})}{\partial x^i}$ is homogeneous function of degree 0, whence equation (2.2.4) follows.

Successively differentiating function $F^2$, we get[2.4]

$$\frac{\partial F^2(\overline{x})}{\partial x^i} = 2F(\overline{x})\frac{\partial F(\overline{x})}{\partial x^i}$$

$$\frac{1}{2}\frac{\partial F^2(\overline{x})}{\partial x^j \partial x^i} = \frac{\partial F(\overline{x})}{\partial x^j}\frac{\partial F(\overline{x})}{\partial x^i} + F(\overline{x})\frac{\partial^2 F(\overline{x})}{\partial x^j \partial x^i} \tag{2.2.6}$$

From equations (2.2.3), (2.2.4), (2.2.6) it follows that

$$\frac{1}{2}\frac{\partial F^2(\overline{x})}{\partial x^j \partial x^i} x^i = \frac{\partial F(\overline{x})}{\partial x^j}\frac{\partial F(\overline{x})}{\partial x^i} x^i + F(\overline{x})\frac{\partial^2 F(\overline{x})}{\partial x^j \partial x^i} x^i = \frac{\partial F(\overline{x})}{\partial x^j} F(\overline{x}) \tag{2.2.7}$$

From the equation (2.2.7) it follows that

$$\frac{1}{2}\frac{\partial F^2(\overline{x})}{\partial x^j \partial x^i} x^i x^j = \frac{\partial F(\overline{x})}{\partial x^j} x^j F(\overline{x}) \tag{2.2.8}$$

The equation (2.2.5) follows from equations (2.2.3), (2.2.8). □

Theorem 2.2.6.

$$g_{ij}(\overline{v}) v^i v^j = F^2(\overline{v}) \tag{2.2.9}$$

Proof. The equation (2.2.9) follows from equations (2.2.2), (2.2.5). □

[2.4] see also [15]

Theorem 2.2.7. *Metric tensor $g_{ij}(\overline{a})$ is homogenious function of degree 0 and satisfies the equation*

(2.2.10)
$$\frac{\partial g_{ij}(\overline{a})}{\partial a^k}a^k=0$$

Proof. From the statement (2.2.1.2) of the definition 2.2.1, it follows that the mapping $F^2(\overline{v})$ is homogeneous of degree 2. From the theorem 2.1.3, it follows that function

$$\frac{\partial F^2(\overline{x})}{\partial x^i}$$

is homogeneous of degree 1. From the theorem 2.1.3 and the definition (2.2.2), it follows that the function $g_{ij}(\overline{x})$ is homogeneous of degree 0. The equation (2.2.10) follows from the theorem 2.1.2. □

Remark 2.2.8. As noted by Rund in [15], tensor

$$C_{ijk}(\overline{a})=\frac{\partial g_{ij}(\overline{a})}{\partial a^k}=\frac{1}{2}\frac{\partial^3 F^2(\overline{a})}{\partial a^k\partial a^i\partial a^j}$$

is symmetric with respect to all indices. Its components are homogeneous function of degree $-1$ and satisfy to following equations

$$\begin{array}{ccccccc}C_{ijk}(\overline{a})a^k&=&C_{kij}(\overline{a})a^k&=&C_{ikj}(\overline{a})a^k&=0\\ \frac{\partial C_{kij}(\overline{a})}{\partial a^h}a^k&=&\frac{\partial C_{ikj}(\overline{a})}{\partial a^h}a^k&=&\frac{\partial C_{ijk}(\overline{a})}{\partial a^h}a^k&=0\end{array}$$

□

Theorem 2.2.9. *Let coordinates of tensor $g$ be defined relative to basis $\overline{\overline{e}}$. Then*

(2.2.11)
$$\left.\frac{\partial g_{ij}(\overline{a})}{\partial a^l}\right|_{\overline{a}=\overline{e}_l}=0$$

Proof. The equation (2.2.11) follows from the equation (2.2.10), since $e^m_i=\delta^m_i$. □

Consider infinitesimal transformation

$$\overline{a}'=\overline{a}+d\overline{a}$$

Then there is infinitesimal transformation of coordinates of metric tensor

(2.2.12)
$$g_{ij}(\overline{a}')=g_{ij}(\overline{a})+\frac{\partial g_{ij}(\overline{a})}{\partial a^k}da^k$$

Definition 2.2.10. The manifold $M$ is called **Finsler space**, if its tangent space is Minkowski space and the Finsler structure $F(x,\overline{v})$ depends continuously on point of tangency $x\in M$. □

Remark 2.2.11. Due to the fact that the Finsler structure in the tangent space depends continuously on point of tangency, there is the ability to determine the differential of length of the curve on the manifold

$$dl=|F(x,d\overline{x})|$$

We usually define Finsler space, and after this we consider tangent to it Minkowski space. In fact, the order of definitions is insignificant. In this paper, Minkowski space is main object of research. □

## 2.3. Orthogonality

As Rund noted in [15], there are various definitions of trigonometric functions in Minkowski space. We are primarily interested in the concept of orthogonality.

Definition 2.3.1. Vector $v_2$ **is orthogonal** to vector $v_1$, if

$$g_{ij}(v_1)v_1^i v_2^j = 0$$

□

As we can see from the definition 2.3.1, relation of orthogonality is noncommutative. This statement has following consequence. Let $n$ be dimension of Minkowski space $\overline{V}$. The set of vectors $\overline{v}$ orthogonal to vector $\overline{a}$ satisfies the linear equation

$$g_{ij}(\overline{a})a^i v^j = 0$$

and, therefore, is vector space of dimension $n-1$. The set of vectors $\overline{v}$ whom vector $\overline{a}$ is orthogonal to satisfies the equation

$$g_{ij}(\overline{v})v^i a^j = 0$$

This set, in general, is not vector space.

Definition 2.3.2. The set of vectors $\overline{e}_1$, ..., $\overline{e}_p$ is called orthogonal if

$$\begin{aligned} &g_{ij}(\overline{e}_k)e_k^i e_k^j \neq 0 \\ &g_{ij}(\overline{e}_k)e_k^i e_l^j = 0 \quad k < l \end{aligned} \tag{2.3.1}$$

The basis $\overline{\overline{e}}$ is called **orthogonal**, if its vectors form orthogonal set. □

Definition 2.3.3. The basis $\overline{\overline{e}}$ is called **orthonormal**, if this is orthogonal basis and its vectors have unit length. □

Since relation of orthogonality is noncommutative, the order of vectors is important when determining an orthogonal basis. There exist different procedures of orthogonalization in Minkowski space, for instance, [18], p. 39. Below we consider the orthogonalization procedure proposed in [3], p. 213 - 214.

Theorem 2.3.4. *Let $\overline{e}_1$, ..., $\overline{e}_p$ be orthogonal set of not isotropic vectors. Then vectors $\overline{e}_1$, ..., $\overline{e}_p$ are linear independent.*

Proof. Consider equation

$$a_1\overline{e}_1 + ... + a_p\overline{e}_p = 0 \tag{2.3.2}$$

From the equation (2.3.2), it follows that

$$a_1 g_{ij}(\overline{e}_1)e_1^i e_1^j + ... + a_p g_{ij}(\overline{e}_1)e_1^i e_p^j = 0 \tag{2.3.3}$$

Since vector $\overline{e}_1$ is not isotropic, then

$$g_{ij}(\overline{e}_1)e_1^i e_1^j \neq 0 \tag{2.3.4}$$

From the conditions (2.3.1), (2.3.4) and the equation (2.3.3), it follows that $a_1 = 0$.

Since we have proved that $a^1 = ... = a^{m-1} = 0$, then the equation (2.3.2) gets form

$$a_m\overline{e}_m + ... + a_p\overline{e}_p = 0 \tag{2.3.5}$$

From the equation (2.3.5), it follows that

$$a_m g_{ij}(\overline{e}_m)e_m^i e_m^j + ... + a_p g_{ij}(\overline{e}_m)e_m^i e_p^j = 0 \tag{2.3.6}$$

Since vector $\overline{e}_m$ is not isotropic, then

$$g_{ij}(\overline{e}_m)e_m^i e_m^j \neq 0 \tag{2.3.7}$$

From the conditions (2.3.1), (2.3.7) and the equation (2.3.6), it follows that $a_m=0$. □

THEOREM 2.3.5. *There exists orthonormal basis in Minkowski space with positively defined Finsler metric.*

PROOF. Let $n$ be dimension of Minkowski space. Let $\overline{\overline{e}}'$ be a basis in Minkowski space.
We set

$$\overline{e}_1=\overline{e}'_1$$

Suppose we have defined the set of vectors $\overline{e}_1$, ..., $\overline{e}_m$. In addition we assume that for every $i$, $1\le i\le m$, the vector $\overline{e}_i$ is linear combination of vectors $\overline{e}'_1$, ..., $\overline{e}'_m$. This assumption also holds for the vector $\overline{e}_{m+1}$, if we represent this vector as

$$\overline{e}_{m+1}=a_1\overline{e}_1+...+a_m\overline{e}_m+\overline{e}'_{m+1}$$

$\overline{e}_{m+1}\ne\overline{0}$, because $\overline{\overline{e}}'$ is basis and vector $\overline{e}'_{m+1}$ is not included in the expansion of vectors $\overline{e}_1$, ..., $\overline{e}_m$. For the choice of the vector $\overline{e}_{m+1}$ we require that the vector $\overline{e}_{m+1}$ is orthogonal to the vectors $\overline{e}_1$, ..., $\overline{e}_m$.

$$\begin{aligned}
&g_{ij}(\overline{e}_1)e^i_1e^j_{m+1}=0\\
&...\\
&g_{ij}(\overline{e}_m)e^i_me^j_{m+1}=0
\end{aligned}\tag{2.3.8}$$

The system of linear equations (2.3.8) has form

$$\begin{aligned}
&a_1g_{ij}(\overline{e}_1)e^i_1e^j_1 &&=-g_{ij}(\overline{e}_1)e^i_1e'^j_{m+1}\\
&a_1g_{ij}(\overline{e}_2)e^i_2e^j_1+a_2g_{ij}(\overline{e}_2)e^i_2e^j_2 &&=-g_{ij}(\overline{e}_2)e^i_2e'^j_{m+1}\\
&...\\
&a_1g_{ij}(\overline{e}_m)e^i_me^j_1+a_2g_{ij}(\overline{e}_m)e^i_me^j_2+...+a_mg_{ij}(\overline{e}_m)e^i_me^j_m &&=-g_{ij}(\overline{e}_m)e^i_me'^j_{m+1}
\end{aligned}$$

Therefore, the solution of the system of linear equations (2.3.8) has form

$$\begin{aligned}
a_1&=-\frac{g_{ij}(\overline{e}_1)e^i_1e'^j_{m+1}}{g_{ij}(\overline{e}_1)e^i_1e^j_1}\\
a_2&=-\frac{g_{ij}(\overline{e}_2)e^i_2e'^j_{m+1}+a_1g_{ij}(\overline{e}_2)e^i_2e^j_1}{g_{ij}(\overline{e}_2)e^i_2e^j_2}\\
...\\
a_m&=-\frac{g_{ij}(\overline{e}_m)e^i_me'^j_{m+1}+\sum_{k=1}^{m-1}a_kg_{ij}(\overline{e}_m)e^i_me^j_k}{g_{ij}(\overline{e}_m)e^i_me^j_m}
\end{aligned}$$

By continuing this process, we obtain an orthogonal basis $\overline{\overline{e}}$. We can normalize vectors of the basis $\overline{\overline{e}}$ according to rule

$$E_k=g_{ij}(\overline{e}_k)e^i_ke^j_k$$

$$\overline{e}_k\rightarrow\frac{1}{\sqrt{|E_k|}}\overline{e}_k$$

$$k=1,...,n$$

□

Consider Minkowski space $\overline{V}$ of dimension $n$. Let Finsler metric $F^2$ of Minkowski space $\overline{V}$ is not positively defined. The set of spacelike vectors is not vector space, however this set contains a maximal set of linear independent vectors

$$\overline{\overline{e}}_+=(\overline{e}_{+\cdot 1},...,\overline{e}_{+\cdot p})$$

such that any linear combination of these vectors is spacelike vector. The set $\overline{\overline{e}}'_+$ is basis of vector space $\overline{V}_+$. The set of timelike vectors is not vector space, however this set contains a maximal set of linear independent vectors

$$\overline{\overline{e}}_- = (\overline{e}_{-\cdot 1}, ..., \overline{e}_{-\cdot q})$$

such that any linear combination of these vectors is timelike vector. The set $\overline{\overline{e}}'_-$ is basis of vector space $\overline{V}_-$. The choice of sets $\overline{\overline{e}}'_+$, $\overline{\overline{e}}'_-$, as well the choice of vector spaces $\overline{V}_+$, $\overline{V}_-$ is not unique.

THEOREM 2.3.6. *There exists orthonormal basis $\overline{\overline{e}}_+$ in Minkowski space $\overline{V}_+$ with Finsler structure $F$. The set of vectors which are orthogonal to vectors of the basis $\overline{\overline{e}}_+$ is vector space $\overline{V}_+^{\perp}$ of dimension $n-p$.*

PROOF. The existence of an orthonormal basis $\overline{\overline{e}}_+$ in Minkowski space $\overline{V}_+$ follows from the theorem 2.3.5. Vector $\overline{v}$ which is orthogonal to vectors of the basis $\overline{\overline{e}}_+$ satisfies the system of linear equations

$$g_{ij}(\overline{e}_{+\cdot 1})e^i_{+\cdot 1}v^j = 0$$
$$...$$
$$g_{ij}(\overline{e}_{+\cdot p})e^i_{+\cdot p}v^j = 0$$

Therefore, the set of vectors $\overline{v}$ generates vector space. □

THEOREM 2.3.7. *There exists orthonormal basis $\overline{\overline{e}}$ in Minkowski space $\overline{V}$ with not positively defined Finsler metric $F^2$ if one of the following conditions is met*

2.3.7.1: *Finsler metric $F^2$ of Minkowski space $\overline{V}$ has form*

$$F^2(\overline{v}) = F^2_+(v^1_+, ..., v^p_+) - F^2_-(v^1_-, ..., v^q_-)$$

*where $F^2_+$, $F^2_-$ are positively defined Finsler metrics. The set of variables $v^1_+$, ..., $v^p_+$, $v^1_-$, ..., $v^q_-$ divides the set of coordinates $v^1$, ..., $v^n$ of vector $\overline{v}$ into two disjoint sets.*

2.3.7.2: *There exists vector space $\overline{V}_+$ such that vector space $\overline{V}_+^{\perp}$ does not contain isotropic vectors*[2.5]

2.3.7.3: *There exist vector spaces $\overline{V}_+$, $\overline{V}_-$ such that following equation is true*

$$\overline{V}_- = \overline{V}_+^{\perp}$$

PROOF. Since the condition 2.3.7.1 is satisfied, then it is natural to assume that $p+q=n$; otherwise, there exists an isotropic plane.[2.6] Since we assume

$$v^1_- = ... = v^q_- = 0$$

we see that vector space $\overline{V}_+$ has dimension $p$. Since we assume

$$v^1_+ = ... = v^p_+ = 0$$

---

[2.5] Although the mapping $F$ is not linear one, we can formally define the set $\ker F$ as the set of isotropic vectors. Then we can write condition 2.3.7.2 as following condition. There exists vector space $\overline{V}_+$ such that

$$\overline{V}_+^{\perp} \cap \ker F = \emptyset$$

[2.6] For instance, if variable $v^1$ does not belong to the set $v^1_+$, ..., $v^p_+$, $v^1_-$, ..., $v^q_-$ then any vector $\overline{v} = \begin{pmatrix} v^1 \\ 0 \\ ... \\ 0 \end{pmatrix}$ is isotropic vector.

This contradicts the requirement 2.2.1.3 of nonsingularity of Finsler metric.

However, the singularity of the Finsler metric is not necessarily the existence of an isotropic plane. In pseudo-Euclidean space with the metric

$$F^2(\overline{v}) = (v^1)^2 - (v^2)^2 + (v^3)^2 - (v^4)^2$$

we see that vector space $\overline{V}_-$ has dimension $q$. According to the theorem 2.3.6, there exists orthonormal basis $\overline{\overline{e}}_+$ in Minkowski space $\overline{V}_+$ with Finsler structure $F$ and dimension of vector space $\overline{V}_+^\perp$ equal $n-p=q$. Therefore, the condition 2.3.7.1 is a particular case of the condition 2.3.7.3.

Spacelike vector cannot be vector of vector space $\overline{V}_+^\perp$. Otherwise the set $\overline{\overline{e}}_+$ is not maximal set of linear independent spacelike vectors. If condition 2.3.7.2 is true, then vector space $\overline{V}_+^\perp$ contains only timelike vectors. Therefore, the condition 2.3.7.2 is equivalent to condition 2.3.7.3.

Consider the condition 2.3.7.3. Let $\overline{\overline{e}}_+$ be orthonormal basis of Minkowski space $\overline{V}_+$. Let $\overline{\overline{e}}_-$ be orthonormal basis of Minkowski space $\overline{V}_-$. According to the condition 2.3.7.3, every vector of the basis $\overline{\overline{e}}_-$ is orthogonal to all vectors of the basis $\overline{\overline{e}}_+$. According to the definition 2.3.2, the basis

$$\overline{\overline{e}} = (\overline{e}_{+\cdot 1}, ..., \overline{e}_{+\cdot p}, \overline{e}_{-\cdot 1}, ..., \overline{e}_{-\cdot q})$$

is orthogonal basis. Since every vector of the basis $\overline{\overline{e}}$ is unit vector, then according to the definition 2.3.3, the basis $\overline{\overline{e}}$ is orthonormal basis. □

THEOREM 2.3.8. *Let $\overline{\overline{e}}$ be orthonormal basis of Minkowski space. If we write the coordinates of the metric tensor $g_{ij}(\overline{e}_k)$ relative to basis $\overline{\overline{e}}$ as matrix*

$$\begin{pmatrix} g_{11}(\overline{e}_1) & ... & g_{1n}(\overline{e}_1) \\ ... & ... & ... \\ g_{n1}(\overline{e}_n) & ... & g_{nn}(\overline{e}_n) \end{pmatrix} \tag{2.3.9}$$

*then the matrix (2.3.9) is a triangular matrix whose diagonal elements are either $1$ or $-1$.*

---

vectors

$$z_1 = \begin{pmatrix} 1 \\ 1 \\ 0 \\ 0 \end{pmatrix} \quad z_2 = \begin{pmatrix} 0 \\ 0 \\ 1 \\ 1 \end{pmatrix}$$

are isotropic vectors. Vector

$$z = a^1 z_1 + a^2 z_2$$

also is isotropic vector. However vectors

$$e_{+\cdot 1} = \begin{pmatrix} 1 \\ 0 \\ 0 \\ 0 \end{pmatrix} \quad e_{+\cdot 2} = \begin{pmatrix} 0 \\ 0 \\ 1 \\ 0 \end{pmatrix}$$

generate vector space $\overline{V}_+$, and vectors

$$e_{-\cdot 1} = \begin{pmatrix} 0 \\ 1 \\ 0 \\ 0 \end{pmatrix} \quad e_{-\cdot 2} = \begin{pmatrix} 0 \\ 0 \\ 0 \\ 1 \end{pmatrix}$$

generate vector space $\overline{V}_-$.

Proof. If we consider coordinates of basis $\overline{\overline{e}}$ relative to the basis $\overline{\overline{e}}$, then $e^i_j = \delta^i_j$. According to the definition 2.3.2, 2.3.3,

(2.3.10)
$$\begin{aligned} g_{ij}(\overline{e}_k)\delta^i_k\delta^j_k &= -1 \quad \overline{e}_k \text{ is timelike vector}\\ g_{ij}(\overline{e}_k)\delta^i_k\delta^j_k &= \ \ 1 \quad \overline{e}_k \text{ is spacelike vector}\\ g_{ij}(\overline{e}_k)\delta^i_k\delta^j_l &= \ \ 0 \quad k<l \end{aligned}$$

From equations (2.3.10), it follows that[2.7]

$$\begin{aligned} g_{kk}(\overline{e}_k) &= -1 \quad \overline{e}_k \text{ is timelike vector}\\ g_{kk}(\overline{e}_k) &= \ \ 1 \quad \overline{e}_k \text{ is spacelike vector}\\ g_{kl}(\overline{e}_k) &= \ \ 0 \quad k<l \end{aligned}$$

Therefore, $g_{kl}(\overline{e}_k)$ is arbitrary, when $k > l$. □

If we relax requirements in the theorem 2.3.5 and assume arbitrary Finsler metric, then we will meet major challenge. The set of vectors $\overline{e}_{m+1}$ satisfying the system of linear equations (2.3.8) generates vector space $\overline{V}_{n-m}$ of dimension $n-m$. If relationship of orthogonality is symmetric, then the system of linear equations (2.3.8) is equivalent to the system of linear equations

(2.3.11)
$$\begin{aligned} g_{ij}(\overline{e}_{m+1})e^i_{m+1}e^j_1 &= 0\\ &\dots\\ g_{ij}(\overline{e}_{m+1})e^i_{m+1}e^j_m &= 0 \end{aligned}$$

From the requirement 2.2.1.3 of nonsingularity of Finsler metric and equations (2.3.11) it follows that there exists vector

(2.3.12)
$$\overline{v} \in \overline{V}_{n-m} \quad F(\overline{v}) \neq 0$$

The statement (2.3.12) guarantees the next step of procedure of orthogonalization. For arbitrary Finsler metric, we cannot guarantee that the procedure of orthogonalization can be finished for initial choice of basis $\overline{\overline{e}}'$, as the next vector may turn out isotropic.

Task of finding of an orthogonal basis and reducing a quadratic form to the sum of the squares are closely related in the geometry of Euclidean space. In Minkowski space, the theorem 2.3.8 makes it impossible for the applicability of the algorithm of reducing a quadratic form to canonical form.

## 2.4. Motion of Minkowski Space

The structure of Minkowski space is close to the structure of Euclidean space. Automorphism of Euclidean space is called motion. Motion of Euclidean space is linear transformation which maps orthonormal basis into orthonormal basis.

Automorphism of Minkowski space also is called motion. Motion of Minkowski space is linear transformation. If $\overline{\overline{e}}_1$, $\overline{\overline{e}}_2$ are orthonormal bases, then there exists unique linear transformation which maps the basis $\overline{\overline{e}}_1$ into basis $\overline{\overline{e}}_2$. However from the theorem 3.4.8, it follows that basis manifold of Minkowski space is not closed

[2.7] It is evident that the equation (2.3.10) says nothing about value of $g_{ij}(\overline{e}_k)$ when $i \neq k$, $j \neq k$, since the coefficient of $g_{ij}(\overline{e}_k)$ in the considered sum is equal to 0. Since $g$ is symmetric tensor, then

$$g_{kl}(\overline{e}_k) = g_{lk}(\overline{e}_k) \quad k \leq l$$

Moreover, if $k < l$, then $g_{kl}(\overline{e}_k) = 0$. However, the value of $g_{lk}(\overline{e}_l)$ is not defined.

relative to linear transformation. So we need to relax requirements to the definition of motion of Minkowski space.

Definition 2.4.1. Linear transformation of Minkowski space mapping at least one orthonormal basis into orthonormal basis is called **motion**. □

Let motion $\overline{A}$ of Minkowski space maps the basis $\overline{\overline{e}}_1$ into the basis $\overline{\overline{e}}_2$. Let $e_i$ be matrix of coordinates of the basis $\overline{\overline{e}}_i$. Then following equation is true

$$e_2 = Ae_1$$

Theorem 2.4.2. *Let motion $\overline{A}$ of Minkowski space maps orthonormal basis $\overline{\overline{e}}_1$ into orthonormal basis $\overline{\overline{e}}_2$. Let motion $\overline{B}$ of Minkowski space maps orthonormal basis $\overline{\overline{e}}_2$ into orthonormal basis $\overline{\overline{e}}_3$. Then the product of motions $\overline{AB}$ is also motion.*

Proof. Since the product of linear transformations is linear transformation, then the theorem follows from the definition 2.4.1 and the statement that linear transformation $\overline{AB}$ maps orthonormal basis $\overline{\overline{e}}_1$ into orthonormal basis $\overline{\overline{e}}_3$. □

Theorem 2.4.3. *Let linear transformation $\overline{A}$ is motion. Then linear transformation $\overline{A}^{-1}$ is motion.*

Proof. Let motion $\overline{A}$ of Minkowski space maps orthonormal basis $\overline{\overline{e}}_1$ into orthonormal basis $\overline{\overline{e}}_2$. Since linear transformation $\overline{A}^{-1}$ maps orthonormal basis $\overline{\overline{e}}_2$ into orthonormal basis $\overline{\overline{e}}_1$, then the theorem follows from the definition 2.4.1. □

According to the definition 2.4.1, the set of motions $SO(\overline{V})$ of Minkowski space $\overline{V}$ is subset of the group of linear transformations $GL$. From theorems 2.4.2, reftheorem: inverse motion it follows that the associative product which has unit and inverse element is defined on the set of motions $SO(\overline{V})$.

Let motion $\overline{A}$ of Minkowski space maps orthonormal basis $\overline{\overline{e}}_1$ into orthonormal basis $\overline{\overline{e}}_2$. Let motion $\overline{B}$ of Minkowski space maps orthonormal basis $\overline{\overline{e}}_2$ into basis $\overline{\overline{e}}_3$. If the basis $\overline{\overline{e}}_3$ is not orthonormal, then we do not know whether or not the product $\overline{AB}$ is motion. Therefore, product of motions is not always defined and the set of motions $SO(\overline{V})$ of Minkowski space $\overline{V}$ is not complete group.

Theorem 2.4.4. *Let coordinates of tensor $g$ be defined relative to orthonormal basis $\overline{\overline{e}}$. Let motion of Minkowski space $\overline{V}$*

$$v'^i = A^i_j v^j \tag{2.4.1}$$

*map the basis $\overline{\overline{e}}$ into orthonormal basis $\overline{\overline{e}}'$. Then*[2.8]

$$g_{ij}(\overline{e}'_k)A^i_k A^j_l = g_{kl}(\overline{e}_k) \quad k \le l \tag{2.4.2}$$

Proof. According to definitions 2.3.2, 2.3.3, vectors of the basis $\overline{\overline{e}}$ satisfy to equation

$$\begin{aligned} &g_{ij}(\overline{e}_k)e^i_k e^j_l = 0 \quad k < l \\ &g_{ij}(\overline{e}_k)e^i_k e^j_k = 1 \end{aligned} \tag{2.4.3}$$

and vectors of the basis $\overline{\overline{e}}'$ satisfy to equation

$$\begin{aligned} &g_{ij}(\overline{e}'_k)e'^i_k e'^j_l = 0 \quad k < l \\ &g_{ij}(\overline{e}'_k)e'^i_k e'^j_k = 1 \end{aligned} \tag{2.4.4}$$

From equations (2.4.1), (2.4.3), (2.4.4), it follows that

$$g_{ij}(\overline{e}_k)e^i_k e^j_l = g_{ij}(\overline{e}'_k)e'^i_k e'^j_l = g_{ij}(\overline{e}'_k)A^i_p e^p_k A^j_q e^q_l \quad k \le l \tag{2.4.5}$$

[2.8] Since the motion does not change the basis, then mappings $g_{ij}$ also do not change.

Since $e^i_j=\delta^i_j$, then equation (2.4.2) follows from equation (2.4.5). □

Corollary 2.4.5. *Let coordinates of tensor $g$ be defined relative to orthonormal basis $\overline{\overline{e}}$. Let motion of Minkowski space $\overline{V}$*

$$v'^i=A^i_jv^j$$

*map the basis $\overline{\overline{e}}$ into orthonormal basis $\overline{\overline{e}}'$. Then*

$$\begin{aligned}g_{ij}(\overline{e}'_k)A^i_kA^j_k&=1\\ g_{ij}(\overline{e}'_k)A^i_kA^j_l&=0\quad k<l\end{aligned}\tag{2.4.6}$$

□

The theorem 2.4.4 follows also follows from the statement that coordinates of motion $\overline{A}$ relative to basis $\overline{\overline{e}}$ coincide with coordinates of basis $\overline{\overline{e}}'$

$$e'^i_j=A^i_j$$

Corollary 2.4.6. *Let $\overline{\overline{e}}_1$, $\overline{\overline{e}}_2$ be orthonormal bases. Let $e_{21\cdot}{}^k_l$ be coordinates of the basis $\overline{\overline{e}}_2$ relative to the basis $\overline{\overline{e}}_1$. Then*

$$g_{1\cdot ij}(\overline{e}'_{2\cdot k})e_{21\cdot}{}^i_k e_{21\cdot}{}^j_l=g_{1\cdot kl}(\overline{e}_{1\cdot k})\quad k\le l$$

*where $g_{1\cdot ij}(\overline{v})$ are coordinates of metric tensor relative to the basis $\overline{\overline{e}}_1$.* □

Theorem 2.4.7. *Let infinitesimal motion*

$$a'^i=a^j(\delta^i_j+A^i_jdt)$$

*map orthonormal basis $\overline{\overline{e}}$ to orthonormal basis $\overline{\overline{e}}'$*

$$e'^i_k=e^j_k(\delta^i_j+A^i_jdt)\tag{2.4.7}$$

*Then $(k\le l)$*

$$g_{kp}(\overline{e}_k)A^p_l+g_{pl}(\overline{e}_k)A^p_k=0\tag{2.4.8}$$

Proof. For $k\le l$, from equations (2.4.7), (2.2.12), it follows that

$$\begin{aligned}g_{ij}(\overline{e}'_k)e'^i_ke'^j_l&=\left(g_{ij}(\overline{e}_k)+\left.\frac{\partial g_{ij}(\overline{a})}{\partial a^r}\right|_{\overline{a}=\overline{e}_k}e^m_kA^r_mdt\right)(e^i_k+e^m_kA^i_mdt)(e^j_l+e^p_lA^j_pdt)\\ &=g_{ij}(\overline{e}_k)e^i_ke^j_l+g_{ij}(\overline{e}_k)e^i_ke^j_pA^p_ldt+g_{ij}(\overline{e}_k)e^i_mA^m_kdte^j_l+\left.\frac{\partial g_{ij}(\overline{a})}{\partial a^r}\right|_{\overline{a}=\overline{e}_k}e^m_kA^r_mdte^i_ke^j_l\end{aligned}\tag{2.4.9}$$

From equations (2.4.3), (2.4.4), (2.4.9), it follows that $(k\le l)$

$$0=g_{ij}(\overline{e}_k)e^i_ke^j_pA^p_ldt+g_{ij}(\overline{e}_k)e^i_mA^m_kdte^j_l+\left.\frac{\partial g_{ij}(\overline{a})}{\partial a^r}\right|_{\overline{a}=\overline{e}_k}e^m_kA^r_mdte^i_ke^j_l\tag{2.4.10}$$

Since $e^i_k=\delta^i_k$, then from equation (2.4.10), it follows that $(k\le l)$

$$g_{kp}(\overline{e}_k)A^p_l+g_{pl}(\overline{e}_k)A^p_k+\left.\frac{\partial g_{kl}(\overline{a})}{\partial a^r}\right|_{\overline{a}=\overline{e}_k}A^r_kdt=0\tag{2.4.11}$$

The equation (2.4.8) follows from equations (2.2.11), (2.4.11). □

Theorem 2.4.8. *The product of infinitesimal motions of Mincovsky space is infinitesimal motion of Mincovsky space.*

Proof. Let

$$f(\overline{a})^i = a^k(\delta^i_k + A^i_k dt)$$
$$g(\overline{a})^i = a^k(\delta^i_k + B^i_k dt)$$

be infinitesimal motions of Mincovsky space. Transformation $fg$ has form

$$\begin{aligned} f(g(\overline{a}))^i &= (g(\overline{a}))^l(\delta^i_l + A^i_l dt) \\ &= a^p(\delta^l_p + B^l_p dt)(\delta^i_l + A^i_l dt) \\ &= a^p(\delta^i_p + A^i_p dt + B^i_p dt) \\ &= a^p(\delta^i_p + (A^i_p + B^i_p)dt) \end{aligned}$$

From the theorem 2.4.7, it follows that coordinates of mappings $f$ and $g$ satisfy equation

$$\begin{aligned} g_{kp}(\overline{e}_k)A^p_l + g_{pl}(\overline{e}_k)A^p_k = 0 \\ g_{kp}(\overline{e}_k)B^p_l + g_{pl}(\overline{e}_k)B^p_k = 0 \end{aligned} \qquad k \le l \tag{2.4.12}$$

From the equation (2.4.12), it follows that $(k \le l)$

$$\begin{aligned} &g_{kp}(\overline{e}_k)(A^p_l + B^p_l) + g_{pl}(\overline{e}_k)(A^p_k + B^p_k) \\ =&g_{kp}(\overline{e}_k)A^p_l + g_{pl}(\overline{e}_k)A^p_k + g_{kp}(\overline{e}_k)B^p_l + g_{pl}(\overline{e}_k)B^p_k \\ =&0 \end{aligned}$$

Therefore, the mapping $fg$ is infinitesimal motion of Mincovsky space. □

Theorem 2.4.9. *Let*

$$a'^i = a^j(\delta^i_j + A^i_j dt)$$

*be infinitesimal motion. For any $r \in R$, the transformation*

$$a'^i = a^j(\delta^i_j + rA^i_j dt)$$

*is infinitesimal motion.*

Proof. From the equation (2.4.8), it follows that $(k \le l)$

$$g_{kp}(\overline{e}_k)(rA^p_l) + g_{pl}(\overline{e}_k)(rA^p_k) = 0 \tag{2.4.13}$$

The theorem follows from the equation (2.4.13) and the theorem 2.4.7. □

Theorem 2.4.10. *The set of infinitesimal motions is vector space over real field.*

Proof. According to the theorem 2.4.8, the set of infinitesimal motions generates Abelian group $g$. According to the theorem 2.4.9, there is representation of real field in group $g$. Therefore, the group $g$ is vector space. □

Since the orthogonality relation is not symmetric, then the structure of metric tensor in an orthonormal basis changes. In particular, since scalar product

$$g_{ij}(\overline{e}_k)e^i_k e^j_l \quad k > l$$

is not defined, then we cannot require that automorphism of Minkowski space preserves scalar product.

THEOREM 2.4.11. *Infinitesimal motion of Mincovsky space*

$$a'^i = a^j(\delta^i_j + A^i_j dt) \tag{2.4.14}$$

*does not preserve the scalar product. For given vectors $\overline{a}$, $\overline{b}$, we can evaluate how much scalar product is not preserved by the expression*

$$g_{ij}(\overline{a}')a'^i b'^j - g_{ij}(\overline{a})a^i b^j = \left(g_{pj}(\overline{a})A^p_i + g_{ip}(\overline{a})A^p_j + \frac{\partial g_{ij}(\overline{a})}{\partial a^p}A^p_l a^l\right)a^i b^j dt \tag{2.4.15}$$

PROOF. Let infinitesimal motion of Mincovsky space map vectors $\overline{a}$, $\overline{b}$ into vectors $\overline{a}'$, $\overline{b}'$

$$\begin{aligned} a'^i &= a^j(\delta^i_j + A^i_j dt)\\ b'^i &= b^j(\delta^i_j + A^i_j dt) \end{aligned} \tag{2.4.16}$$

Scalar product of vectors $\overline{a}'$, $\overline{b}'$ has form

$$\begin{aligned} g_{ij}(\overline{a}')a'^i b'^j &= g_{ij}((a^k + A^k_l a^l dt)\overline{e}_k)(a^i + A^i_p a^p dt)(b^j + A^j_q b^q dt)\\ &= (g_{ij}(\overline{a}) + \frac{\partial g_{ij}(\overline{a})}{\partial a^r}A^r_l a^l dt)(a^i + A^i_p a^p dt)(b^j + A^j_q b^q dt)\\ &= g_{ij}(\overline{a})a^i b^j + g_{ij}(\overline{a})A^i_p a^p b^j dt + g_{ij}(\overline{a})a^i A^j_p b^p dt + \frac{\partial g_{ij}(\overline{a})}{\partial a^r}A^r_l a^l dt a^i b^j\\ &= g_{ij}(\overline{a})a^i b^j + \left(g_{pj}(\overline{a})A^p_i + g_{ip}(\overline{a})A^p_j + \frac{\partial g_{ij}(\overline{a})}{\partial a^p}A^p_l a^l\right)a^i b^j dt \end{aligned} \tag{2.4.17}$$

The equation (2.4.15) follows from the equation (2.4.17). □

## 2.5. Quasimotion of Minkowski Space

Even the structure of the set $SO(\overline{V})$ of motions of Minkowski space $\overline{V}$ is not completely defined group, the set $SO(\overline{V})$ is subset of group $GL$ and acts single transitive on the basis manifold of Minkowski space $\overline{V}$. Hence in the set of passive transformations of vector space $\overline{V}$, we can distinguish the transformations that preserve the structure of Minkowski space $\overline{V}$.

DEFINITION 2.5.1. Passive transformation of Minkowski space

$$\overline{e}'_i = A^j_i \overline{e}_j$$

mapping at least one orthonormal basis into orthonormal basis is called **quasimotion of Minkowski space**. □

Let quasimotion $\overline{A}$ of Minkowski space maps the basis $\overline{\overline{e}}_1$ into the basis $\overline{\overline{e}}_2$. Let $e_i$ be matrix of coordinates of the basis $\overline{\overline{e}}_i$. Then following equation is true

$$e_2 = e_1 A$$

Coordinates of vector

$$\overline{a} = a^i \overline{e}_i$$

transform according to the rule

$$a'^j = A^{-1}{}^j_i a^i \tag{2.5.1}$$

From equations (2.2.9), (2.5.1) it follows that

$$g_{ij}(\overline{a})a^i a^j = g'_{kl}(\overline{a})a'^k a'^l = g'_{kl}(\overline{a})A^{-1}{}^k_i a^i A^{-1}{}^l_j a^j \tag{2.5.2}$$

From the equation (2.5.2), it follows that

$$g'_{kl}(\overline{a}) = g_{ij}(\overline{a})A^i_k A^j_l \tag{2.5.3}$$

Therefore, $g_{ij}(\overline{a})$ is tensor.

REMARK 2.5.2. Until now we have considered the metric of Minkowski space $\overline{V}$ as set of functions $g_{ij}(\overline{v})$ of vector $\overline{v}\in\overline{V}$. Since the basis $\overline{\overline{e}}$ is given, then we have unique expansion of vector $\overline{v}$ relative to the basis $\overline{\overline{e}}$

$$\overline{v}=v^i\overline{e}_i$$

Therefore, we can consider the mapping $g_{ij}(\overline{v})$ as function $g_{ij}(v^1,...,v^n)$ of variables $v^1$, ..., $v^n$

$$g_{ij}(v^1,...,v^n)=g_{ij}(\overline{v}) \tag{2.5.4}$$

□

THEOREM 2.5.3. *Let coordinates of tensor $g$ be defined relative to orthonormal basis $\overline{\overline{e}}$. Let quasimotion of Minkowski space $\overline{V}$*

$$\overline{e}'_i=A_i^j\overline{e}_j \tag{2.5.5}$$

*map the basis $\overline{\overline{e}}$ into orthonormal basis $\overline{\overline{e}}'$. Then*

$$g_{ij}(\overline{e}'_k)A_k^iA_l^j=g_{kl}(\overline{e}_k)\quad k\le l \tag{2.5.6}$$

PROOF. According to definitions 2.3.2, 2.3.3, vectors of the basis $\overline{\overline{e}}$ satisfy to equation

$$\begin{aligned}&g_{ij}(\overline{e}_k)e_k^ie_l^j=0\quad k<l\\&g_{ij}(\overline{e}_k)e_k^ie_k^j=1\end{aligned} \tag{2.5.7}$$

and vectors of the basis $\overline{\overline{e}}'$ satisfy to equation

$$\begin{aligned}&g_{ij}(\overline{e}'_k)e'^i_ke'^j_l=0\quad k<l\\&g_{ij}(\overline{e}'_k)e'^i_ke'^j_k=1\end{aligned} \tag{2.5.8}$$

From equations (2.5.5), (2.5.7), (2.5.8), it follows that

$$g_{ij}(\overline{e}_k)e_k^ie_l^j=g_{ij}(\overline{e}'_k)e'^i_ke'^j_l=g_{ij}(\overline{e}'_k)A_k^pe_p^iA_l^qe_q^j\quad k\le l \tag{2.5.9}$$

Since $e_j^i=\delta_j^i$, then equation (2.5.6) follows from equation (2.5.9). □

From theorems 2.4.4, 2.5.3, it follows that the set of quasimotions of Minkowski space and the set of motions of Minkowski space coincide.

REMARK 2.5.4. The set of quasimotions generates passive representations of not complete group $SO(\overline{V})$.

Let $\overline{\overline{e}}_1$, $\overline{\overline{e}}_2$ be orthonormal bases of Minkowski space $\overline{V}$. According to the theorem 2.5.3, there exists quasimotion mapping the basis $\overline{\overline{e}}_1$ into basis $\overline{\overline{e}}_2$. Since not complete group $SO(\overline{V})$ is subgroup of $GL$, this quasimotion is uniquely determined. Therefore, the representation of not complete group $SO(\overline{V})$ on basis manifold of Minkowski space $\overline{V}$ is single transitive.

Therefore, the definition of a geometric object (the section [7]-5.3) remains valid in Minkowski space as well. □

THEOREM 2.5.5. *Let $\overline{A}$ be quasimotion mapping the basis $\overline{\overline{e}}$ into the basis $\overline{\overline{e}}'$*

$$\overline{e}'_k=A_k^i\overline{e}_i \tag{2.5.10}$$

*Quasimotion $\overline{A}$ maps the mapping $g_{ij}(v^1,...,v^n)$ into mapping $g'_{ij}(v'^1,...,v'^n)$ according to rule*

$$g'_{kl}(v'^1,...,v'^n)=g_{ij}(v'^kA_k^1,...,v'^kA_k^n)A_k^iA_l^j \tag{2.5.11}$$

Proof. From equations (2.5.3), (2.5.4), it follows that

$$g'_{kl}(v'^1,...,v'^n)=g'_{kl}(\overline{v})=g_{ij}(\overline{v})A^i_kA^j_l=g_{ij}(v^1,...,v^n)A^i_kA^j_l \tag{2.5.12}$$

Since quasimotion $\overline{A}$ does not change the vector $\overline{v}$, then

$$\overline{v}=v^k\overline{e}_k=v'^i\overline{e}'_i=v'^kA^i_k\overline{e}_i \tag{2.5.13}$$

From the equation (2.5.13) it follows that

$$v^i=v'^kA^i_k \tag{2.5.14}$$

The equation (2.5.11) follows from equations (2.5.12), (2.5.14). □

Theorem 2.5.6. *Let infinitesimal quasimotion*

$$\overline{e}'_k=\overline{e}_l(\delta^l_k+A^l_kdt)$$

*map orthonormal basis $\overline{\overline{e}}$ to orthonormal basis $\overline{\overline{e}}'$*

$$e'^i_k=e^i_l(\delta^l_k+A^l_kdt) \tag{2.5.15}$$

*Let coordinates of tensor $g$ be defined relative to orthonormal basis $\overline{\overline{e}}$. Let coordinates of tensor $g'$ be defined relative to orthonormal basis $\overline{\overline{e}}'$. Then $(k\le l)$*

$$g_{kp}(\overline{e}_k)A^p_l+g_{lp}(\overline{e}_k)A^p_k=0 \tag{2.5.16}$$

Proof. According to (2.5.3), (2.2.12), it follows that there is infinitesimal transformation of metric tensor[2.9]

$$\begin{aligned}g'_{kl}(\overline{e}'_p)&=g_{ij}(\overline{e}'_p)(\delta^i_k\delta^j_l+\delta^i_kA^j_ldt+A^i_k\delta^j_ldt)\\&=\left(g_{ij}(\overline{e}_p)+\left.\frac{\partial g_{ij}(\overline{a})}{\partial a^r}\right|_{\overline{a}=\overline{e}_p}\delta^r_mA^m_pdt\right)(\delta^i_k\delta^j_l+\delta^i_kA^j_ldt+A^i_k\delta^j_ldt)\\&=g_{kl}(\overline{e}_p)+g_{kj}(\overline{e}_p)A^j_ldt+g_{il}(\overline{e}_p)A^i_kdt+\left.\frac{\partial g_{kl}(\overline{a})}{\partial a^m}\right|_{\overline{a}=\overline{e}_p}A^m_pdt\end{aligned} \tag{2.5.17}$$

For $p=k$, from the equation (2.5.17), it follows that

$$g'_{kl}(\overline{e}'_k)=g_{kl}(\overline{e}_k)+g_{kj}(\overline{e}_k)A^j_ldt+g_{il}(\overline{e}_k)A^i_kdt+\left.\frac{\partial g_{kl}(\overline{a})}{\partial a^m}\right|_{\overline{a}=\overline{e}_k}A^m_pdt \tag{2.5.18}$$

For $k\le l$, according to the theorem 2.3.8,

$$g_{kl}(\overline{e}_k)=g'_{kl}(\overline{e}'_k) \tag{2.5.19}$$

The equation (2.5.16) follows from equations (2.2.11), (2.5.18), (2.5.19). □

From theorems 2.4.7, 2.5.6, it follows that the set of infinitesimal quasimotions of Minkowski space and the set of infinitesimal motions of Minkowski space coincide.

[2.9] We use equation $e^m_i=\delta^m_i$.

CHAPTER 3

# Examples of Minkowski Space

## 3.1. Example of Minkowski Space of Dimension 2

We consider example of the Finslerian metric of Minkowski space of dimension 2 (considered metric is a special case of metric [12]-(3))

$$\begin{aligned} F^2(\overline{v}) &= ((v^1)^2+(v^2)^2)\left(1+\frac{(v^1)^3v^2}{((v^1)^2+(v^2)^2)^2}\right) \\ &= (v^1)^2+(v^2)^2+\frac{(v^1)^3v^2}{(v^1)^2+(v^2)^2} \end{aligned} \tag{3.1.1}$$

Theorem 3.1.1. *The Finslerian metric defined by equation (3.1.1) generates metric of Minkowski space*

$$g_{11}(\overline{v}) = 1+\frac{3v^1(v^2)^5-(v^1)^3(v^2)^3}{((v^1)^2+(v^2)^2)^3} \tag{3.1.2}$$

$$g_{12}(\overline{v}) = \frac{1}{2}\frac{(v^1)^6+6(v^1)^4(v^2)^2-3(v^1)^2(v^2)^4}{((v^1)^2+(v^2)^2)^3} \tag{3.1.3}$$

$$g_{22}(\overline{v}) = 1+\frac{(v^1)^3(v^2)^3-3(v^1)^5v^2}{((v^1)^2+(v^2)^2)^3} \tag{3.1.4}$$

Proof. From the equation (3.1.1) it follows that

$$\begin{aligned} \frac{\partial F^2(\overline{v})}{\partial v^1} &= 2v^1+\frac{3(v^1)^2v^2((v^1)^2+(v^2)^2)-(v^1)^3v^2 2v^1}{((v^1)^2+(v^2)^2)^2} \\ &= 2v^1+\frac{3(v^1)^4v^2+3(v^1)^2(v^2)^3-2(v^1)^4v^2}{((v^1)^2+(v^2)^2)^2} \\ &= 2v^1+\frac{(v^1)^4v^2+3(v^1)^2(v^2)^3}{((v^1)^2+(v^2)^2)^2} \end{aligned} \tag{3.1.5}$$

$$\begin{aligned} \frac{\partial F^2(\overline{v})}{\partial v^2} &= 2v^2+\frac{(v^1)^3((v^1)^2+(v^2)^2)-(v^1)^3v^2 2v^2}{((v^1)^2+(v^2)^2)^2} \\ &= 2v^2+\frac{(v^1)^5-(v^1)^3(v^2)^2}{((v^1)^2+(v^2)^2)^2} \end{aligned} \tag{3.1.6}$$

From equations (2.2.2), (3.1.5), (3.1.6), it follows that

$$
\begin{aligned}
g_{11}(\overline{v}) =& 1+\frac{1}{2}\frac{(4(v^1)^3v^2+6v^1(v^2)^3)((v^1)^2+(v^2)^2)^2}{((v^1)^2+(v^2)^2)^4}\\
&-\frac{1}{2}\frac{((v^1)^4v^2+3(v^1)^2(v^2)^3)2((v^1)^2+(v^2)^2)2v^1}{((v^1)^2+(v^2)^2)^4}\\
=& 1+\frac{2(v^1)^5v^2+3(v^1)^3(v^2)^3+2(v^1)^3(v^2)^3+3v^1(v^2)^5}{((v^1)^2+(v^2)^2)^3}-\frac{2(v^1)^5v^2+6(v^1)^3(v^2)^3}{((v^1)^2+(v^2)^2)^3}
\end{aligned}
\tag{3.1.7}
$$

$$
\begin{aligned}
g_{12}(\overline{v}) =& \frac{1}{2}\frac{(5(v^1)^4-3(v^1)^2(v^2)^2)((v^1)^2+(v^2)^2)^2}{((v^1)^2+(v^2)^2)^4}\\
&-\frac{1}{2}\frac{((v^1)^5-(v^1)^3(v^2)^2)2((v^1)^2+(v^2)^2)2v^1}{((v^1)^2+(v^2)^2)^4}\\
=& \frac{1}{2}\frac{5(v^1)^6-3(v^1)^4(v^2)^2+5(v^1)^4(v^2)^2-3(v^1)^2(v^2)^4}{((v^1)^2+(v^2)^2)^3}-\frac{1}{2}\frac{4(v^1)^6-4(v^1)^4(v^2)^2}{((v^1)^2+(v^2)^2)^3}
\end{aligned}
\tag{3.1.8}
$$

$$
\begin{aligned}
g_{22}(\overline{v}) =& 1+\frac{1}{2}\frac{-2(v^1)^3v^2((v^1)^2+(v^2)^2)^2}{((v^1)^2+(v^2)^2)^4}\\
&-\frac{1}{2}\frac{((v^1)^5-(v^1)^3(v^2)^2)2((v^1)^2+(v^2)^2)2v^2}{((v^1)^2+(v^2)^2)^4}\\
=& 1-\frac{(v^1)^5v^2+(v^1)^3(v^2)^3}{((v^1)^2+(v^2)^2)^3}-\frac{2(v^1)^5v^2-2(v^1)^3(v^2)^3}{((v^1)^2+(v^2)^2)^3}
\end{aligned}
\tag{3.1.9}
$$

The equation (3.1.2) follows from the equation (3.1.7). The equation (3.1.3) follows from the equation (3.1.8). The equation (3.1.4) follows from the equation (3.1.9). □

We will start the procedure of orthogonalization[3.1] from the basis

$$\overline{e}_1=\begin{pmatrix}1\\0\end{pmatrix}\quad \overline{e}_2=\begin{pmatrix}0\\1\end{pmatrix}$$

Let

$$\overline{e}'_1=\overline{e}_1\quad \overline{e}'_2=\begin{pmatrix}e'^1_2\\e'^2_2\end{pmatrix}$$

be orthogonal basis. Since vector $\overline{e}'_2$ is orthogonal to vector $\overline{e}'_1$, then according to the definition 2.3.1

$$
\begin{aligned}
&g_{ij}(\overline{e}'_1)e'^i_1e'^j_2=0\\
&g_{1j}(\overline{e}'_1)e'^j_2=0
\end{aligned}
\tag{3.1.10}
$$

According to equations (3.1.2), (3.1.3),

$$g_{11}(\overline{e}'_1)=1\quad g_{12}(\overline{e}'_1)=\frac{1}{2}\quad g_{22}(\overline{e}'_1)=1 \tag{3.1.11}$$

From equations (3.1.10), (3.1.11), it follows that

$$e'^1_2+\frac{1}{2}e'^2_2=0 \tag{3.1.12}$$

[3.1] See section [8]-3.

From the equation (3.1.12), it follows that we can assume

$$\overline{e}'_2=\begin{pmatrix}-1\\2\end{pmatrix} \tag{3.1.13}$$

According to equations (3.1.2), (3.1.3), (3.1.4), (3.1.13),

$$\left\{\begin{aligned}g_{11}(\overline{e}'_2)&=1+\frac{3(-1)(2)^5-(-1)^3(2)^3}{((-1)^2+(2)^2)^3}\\&=1+\frac{-96+8}{(1+4)^3}=\frac{37}{125}\\g_{12}(\overline{e}'_2)&=\frac{1}{2}\frac{(-1)^6+6(-1)^4(2)^2-3(-1)^2(2)^4}{((-1)^2+(2)^2)^3}\\&=\frac{1}{2}\frac{1+24-48}{(1+4)^3}=-\frac{23}{250}\\g_{22}(\overline{e}'_2)&=1+\frac{(-1)^3(2)^3-3(-1)^52}{((-1)^2+(2)^2)^3}\\&=1+\frac{-8+6}{(1+4)^3}=\frac{123}{125}\end{aligned}\right. \tag{3.1.14}$$

From equations (3.1.13), (3.1.14), it follows that

$$\begin{aligned}g_{ij}(\overline{e}'_2)e'^i_2e'^j_1&=g_{i1}(\overline{e}'_2)e'^i_2\\&=\frac{37}{125}(-1)+(-\frac{23}{250})2=-\frac{37+23}{125}\\&=-\frac{12}{25}\end{aligned}$$

Therefore, the vector $\overline{e}'_1$ is not orthogonal to vector $\overline{e}'_2$.

Since

$$F^2(\overline{e}'_1)=1$$

then the length of the vector

$$\overline{e}''_1=\overline{e}'_1=\begin{pmatrix}1\\0\end{pmatrix}$$

equal 1. Since

$$F^2(\overline{e}'_2)=(-1)^2+(2)^2+\frac{(-1)^32}{(-1)^2+(2)^2}=1+4-\frac{2}{1+4}=\frac{23}{5}$$

then the length of the vector

$$\overline{e}''_2=\sqrt{\frac{5}{23}}\overline{e}'_2=\begin{pmatrix}-\sqrt{\frac{5}{23}}\\2\sqrt{\frac{5}{23}}\end{pmatrix}$$

equal 1. Therefore, the basis $\overline{\overline{e}}''$ is orthonormal basis. There is passive transformation between bases $\overline{\overline{e}}''$ and $\overline{\overline{e}}$

$$\begin{pmatrix}\overline{e}''_1 & \overline{e}''_2\end{pmatrix} = \begin{pmatrix}\overline{e}_1 & \overline{e}_2\end{pmatrix}\begin{pmatrix}1 & -\sqrt{\frac{5}{23}}\\ 0 & 2\sqrt{\frac{5}{23}}\end{pmatrix} \tag{3.1.15}$$

THEOREM 3.1.2. *Coordinates of metric tensor relative to the basis $\overline{\overline{e}}''$ have form*

$$g''(\overline{e}''_1) = \begin{pmatrix}1 & 0\\ 0 & \frac{15}{23}\end{pmatrix} \tag{3.1.16}$$

$$g''(\overline{e}''_2) = \begin{pmatrix}\frac{37}{125} & -\frac{12}{25}\sqrt{\frac{5}{23}}\\ -\frac{12}{25}\sqrt{\frac{5}{23}} & 1\end{pmatrix} \tag{3.1.17}$$

PROOF. From equations (2.5.3), (3.1.11), (3.1.15), it follows that

$$\left\{\begin{aligned}
&g''_{11}(\overline{e}''_1) = g_{11}(\overline{e}''_1) = 1\\
&g''_{12}(\overline{e}''_1) = g_{11}(\overline{e}''_1)A^1_2 + g_{12}(\overline{e}''_1)A^2_2 = -\sqrt{\frac{5}{23}} + \frac{1}{2}2\sqrt{\frac{5}{23}} = 0\\
&g''_{22}(\overline{e}''_1) = g_{11}(\overline{e}''_1)A^1_2A^1_2 + 2g_{12}(\overline{e}''_1)A^1_2A^2_2 + g_{22}(\overline{e}''_1)A^2_2A^2_2\\
&\qquad = \left(-\sqrt{\frac{5}{23}}\right)\left(-\sqrt{\frac{5}{23}}\right) + 2\frac{1}{2}\left(-\sqrt{\frac{5}{23}}\right)2\sqrt{\frac{5}{23}} + 2\sqrt{\frac{5}{23}}2\sqrt{\frac{5}{23}} = \frac{15}{23}
\end{aligned}\right. \tag{3.1.18}$$

The equation (3.1.17) follows from (3.1.18). From equations (2.5.3), (3.1.14), (3.1.15), it follows that

$$\left\{\begin{aligned}
&g''_{11}(\overline{e}''_2) = g_{11}(\overline{e}''_2) = \frac{37}{125}\\
&g''_{12}(\overline{e}''_2) = g_{11}(\overline{e}''_2)A^1_2 + g_{12}(\overline{e}''_2)A^2_2\\
&\qquad = -\frac{37}{125}\sqrt{\frac{5}{23}} + \left(-\frac{23}{250}\right)2\sqrt{\frac{5}{23}} = -\frac{12}{25}\sqrt{\frac{5}{23}}\\
&g''_{22}(\overline{e}''_2) = g_{11}(\overline{e}''_2)A^1_2A^1_2 + 2g_{12}(\overline{e}''_2)A^1_2A^2_2 + g_{22}(\overline{e}''_2)A^2_2A^2_2\\
&\qquad = \frac{37}{125}\left(-\sqrt{\frac{5}{23}}\right)\left(-\sqrt{\frac{5}{23}}\right) + 2\left(-\frac{23}{250}\right)\left(-\sqrt{\frac{5}{23}}\right)2\sqrt{\frac{5}{23}}\\
&\qquad + \frac{123}{125}2\sqrt{\frac{5}{23}}2\sqrt{\frac{5}{23}}\\
&\qquad = \left(\frac{37}{125} + \frac{46}{125}\right)\frac{5}{23} + \frac{492}{125}\frac{5}{23} = \frac{575}{125}\frac{5}{23} = 1
\end{aligned}\right. \tag{3.1.19}$$

The equation (3.1.2) follows from (3.1.19). □

COROLLARY 3.1.3. *Relative to the basis (3.1.15), the matrix (2.3.9) corresponding to the Finslerian metric defined by equation (3.1.1) has form*

$$\begin{pmatrix}1 & 0\\ -\frac{12}{25}\sqrt{\frac{5}{23}} & 1\end{pmatrix} \tag{3.1.20}$$

□

## 3.2. Structure of Metric of Minkowski Space

The theorem 2.3.8 defines value of $g_{kl}(\overline{e}_k)$ when $k \ge l$. However the theorem 2.3.8 tells nothing about value $g_{kl}(\overline{e}_k)$ when $k < l$. The impression that this value is arbitrary creates a problem. For instance, if metric does not depend on direction, then Minkowski space becomes Euclid space and matrix (2.3.9) becomes diagonal matrix; thus it is not clear how consequence of arbitrary values converges to 0.

From the equation (3.1.20), it follows that the value $g_{kl}(\overline{e}_k)$, $k < l$, is unknown, but it is not arbitrary. In the section 3.1, we have made calculations for the Finslerian metric defined by the equation (3.1.1). Our goal is to reproduce these calculations for arbitrary norm $F$ of Minkowski space of dimension 2.

Let $F(\overline{v})$ be the Finsler structure of Minkowski space of dimension 2. Let $g(\overline{v})$ be metric defined by the equation (2.2.2).

We will start the procedure of orthogonalization[3.2] from the basis

$$\overline{e}_1 = \begin{pmatrix}1\\0\end{pmatrix} \quad \overline{e}_2 = \begin{pmatrix}0\\1\end{pmatrix}$$

Let

$$\overline{e}'_1 = \overline{e}_1 \quad \overline{e}'_2 = \begin{pmatrix}e'^1_2\\e'^2_2\end{pmatrix}$$

be orthogonal basis. Since vector $\overline{e}'_2$ is orthogonal to vector $\overline{e}'_1$, then according to the definition 2.3.1

$$g_{ij}(\overline{e}'_1)e'^i_1 e'^j_2 = 0$$

$$g_{1j}(\overline{e}'_1)e'^j_2 = g_{11}(\overline{e}'_1)e'^1_2 + g_{12}(\overline{e}'_1)e'^2_2 = 0 \tag{3.2.1}$$

From the equation (3.2.1), it follows that we can assume

$$\overline{e}'_2 = \begin{pmatrix}-g_{12}(\overline{e}'_1)\\g_{11}(\overline{e}'_1)\end{pmatrix} \tag{3.2.2}$$

From the equation (3.2.2), it follows that

$$g_{ij}(\overline{e}'_2)e'^i_2 e'^j_1 = g_{i1}(\overline{e}'_2)e'^i_2 = -g_{11}(\overline{e}'_2)g_{12}(\overline{e}'_1) + g_{21}(\overline{e}'_2)g_{11}(\overline{e}'_1)$$
$$= \begin{vmatrix}g_{11}(\overline{e}'_1) & g_{11}(\overline{e}'_2)\\g_{12}(\overline{e}'_1) & g_{12}(\overline{e}'_2)\end{vmatrix}$$

Let

$$F(v^1, v^2) = F(\overline{v})$$

Then the orthonormal basis $\overline{\overline{e}}''_2$ has form

$$\begin{pmatrix}\overline{e}''_1 & \overline{e}''_2\end{pmatrix} = \begin{pmatrix}\overline{e}_1 & \overline{e}_2\end{pmatrix}\begin{pmatrix}\frac{1}{F(1,0)} & -\frac{g_{12}(1,0)}{F(-g_{12}(1,0), g_{11}(1,0))}\\0 & \frac{g_{11}(1,0)}{F(-g_{12}(1,0), g_{11}(1,0))}\end{pmatrix} \tag{3.2.3}$$

[3.2] See section [8]-3.

### 3.3. Example 1 of Minkowski Space of Dimension 3

We consider example of the Finslerian metric of Minkowski space of dimension 3 (considered metric is a special case of metric [12]-(3))

$$
\begin{aligned}
F^2(\overline{v}) &= ((v^1)^2+(v^2)^2+(v^3)^2)\left(1+\frac{(v^1)^2v^2v^3}{((v^1)^2+(v^2)^2+(v^3)^2)^2}\right)\\
&= (v^1)^2+(v^2)^2+(v^3)^2+\frac{(v^1)^2v^2v^3}{(v^1)^2+(v^2)^2+(v^3)^2}
\end{aligned}
\tag{3.3.1}
$$

Theorem 3.3.1. *The Finslerian metric defined by equation (3.3.1) generates metric of Minkowski space*

$$
g_{11}(\overline{v}) = 1+\frac{-3(v^1)^2(v^2)^3v^3-3(v^1)^2v^2(v^3)^3+(v^2)^5v^3+2(v^2)^3(v^3)^3+v^2(v^3)^5}{((v^1)^2+(v^2)^2+(v^3)^2)^3}
\tag{3.3.2}
$$

$$
g_{12}(\overline{v}) = \frac{3(v^1)^3(v^2)^2v^3-v^1(v^2)^4v^3+(v^1)^3(v^3)^3+v^1(v^3)^5}{((v^1)^2+(v^2)^2+(v^3)^2)^3}
\tag{3.3.3}
$$

$$
g_{13}(\overline{v}) = \frac{(v^1)^3(v^2)^3+v^1(v^2)^5+3(v^1)^3v^2(v^3)^2-v^1v^2(v^3)^4}{((v^1)^2+(v^2)^2+(v^3)^2)^3}
\tag{3.3.4}
$$

$$
g_{22}(\overline{v}) = 1+\frac{-3(v^1)^4v^2v^3+(v^1)^2(v^2)^3v^3-3(v^1)^2v^2(v^3)^3}{((v^1)^2+(v^2)^2+(v^3)^2)^3}
\tag{3.3.5}
$$

$$
g_{23}(\overline{v}) = \frac{1}{2}\frac{(v^1)^6-(v^1)^2(v^3)^4-(v^1)^2(v^2)^4+6(v^1)^2(v^2)^2(v^3)^2}{((v^1)^2+(v^2)^2+(v^3)^2)^3}
\tag{3.3.6}
$$

$$
g_{33}(\overline{v}) = 1+\frac{(v^1)^2v^2(v^3)^3-3(v^1)^4v^2v^3-3(v^1)^2(v^2)^3v^3}{((v^1)^2+(v^2)^2+(v^3)^2)^3}
\tag{3.3.7}
$$

Proof. From the equation (3.3.1) it follows that

$$
\begin{aligned}
\frac{\partial F^2(\overline{v})}{\partial v^1} &= 2v^1+\frac{2v^1v^2v^3((v^1)^2+(v^2)^2+(v^3)^2)-(v^1)^2v^2v^3 2v^1}{((v^1)^2+(v^2)^2+(v^3)^2)^2}\\
&= 2v^1+\frac{2v^1(v^2)^3v^3+2v^1v^2(v^3)^3}{((v^1)^2+(v^2)^2+(v^3)^2)^2}
\end{aligned}
\tag{3.3.8}
$$

$$
\begin{aligned}
\frac{\partial F^2(\overline{v})}{\partial v^2} &= 2v^2+\frac{(v^1)^2v^3((v^1)^2+(v^2)^2+(v^3)^2)-(v^1)^2v^2v^3 2v^2}{((v^1)^2+(v^2)^2+(v^3)^2)^2}\\
&= 2v^2+\frac{(v^1)^4v^3-(v^1)^2(v^2)^2v^3+(v^1)^2(v^3)^3}{((v^1)^2+(v^2)^2+(v^3)^2)^2}
\end{aligned}
\tag{3.3.9}
$$

$$
\begin{aligned}
\frac{\partial F^2(\overline{v})}{\partial v^3} &= 2v^3+\frac{(v^1)^2v^2((v^1)^2+(v^2)^2+(v^3)^2)-(v^1)^2v^2v^3 2v^3}{((v^1)^2+(v^2)^2+(v^3)^2)^2}\\
&= 2v^3+\frac{(v^1)^4v^2+(v^1)^2(v^2)^3-(v^1)^2v^2(v^3)^2}{((v^1)^2+(v^2)^2+(v^3)^2)^2}
\end{aligned}
\tag{3.3.10}
$$

From equations (2.2.2), (3.3.8), (3.3.9), (3.3.10), it follows that

$$
(3.3.11)\quad
\left\{
\begin{aligned}
g_{11}(\overline{v}) =& 1+\frac{1}{2}\frac{(2(v^2)^3v^3+2v^2(v^3)^3)((v^1)^2+(v^2)^2+(v^3)^2)^2}{((v^1)^2+(v^2)^2+(v^3)^2)^4}\\
&-\frac{1}{2}\frac{(2v^1(v^2)^3v^3+2v^1v^2(v^3)^3)2((v^1)^2+(v^2)^2+(v^3)^2)2v^1}{((v^1)^2+(v^2)^2+(v^3)^2)^4}\\
&=1+\frac{((v^2)^3v^3+v^2(v^3)^3)((v^1)^2+(v^2)^2+(v^3)^2)}{((v^1)^2+(v^2)^2+(v^3)^2)^3}\\
&-\frac{(v^1(v^2)^3v^3+v^1v^2(v^3)^3)4v^1}{((v^1)^2+(v^2)^2+(v^3)^2)^3}\\
&=1+\frac{(v^1)^2(v^2)^3v^3+(v^1)^2v^2(v^3)^3+(v^2)^5v^3+(v^2)^3(v^3)^3}{((v^1)^2+(v^2)^2+(v^3)^2)^3}\\
&+\frac{(v^2)^3(v^3)^3+v^2(v^3)^3(v^3)^2}{((v^1)^2+(v^2)^2+(v^3)^2)^3}-\frac{4(v^1)^2(v^2)^3v^3+4(v^1)^2v^2(v^3)^3}{((v^1)^2+(v^2)^2+(v^3)^2)^3}
\end{aligned}
\right.
$$

$$
(3.3.12)\quad
\left\{
\begin{aligned}
&g_{12}(\overline{v})\\
&=\frac{1}{2}\frac{(4(v^1)^3v^3-2v^1(v^2)^2v^3+2v^1(v^3)^3)((v^1)^2+(v^2)^2+(v^3)^2)^2}{((v^1)^2+(v^2)^2+(v^3)^2)^4}\\
&-\frac{1}{2}\frac{((v^1)^4v^3-(v^1)^2(v^2)^2v^3+(v^1)^2(v^3)^3)2((v^1)^2+(v^2)^2+(v^3)^2)2v^1}{((v^1)^2+(v^2)^2+(v^3)^2)^4}\\
&=\frac{(2(v^1)^3v^3-v^1(v^2)^2v^3+v^1(v^3)^3)((v^1)^2+(v^2)^2+(v^3)^2)}{((v^1)^2+(v^2)^2+(v^3)^2)^3}\\
&-\frac{((v^1)^4v^3-(v^1)^2(v^2)^2v^3+(v^1)^2(v^3)^3)2v^1}{((v^1)^2+(v^2)^2+(v^3)^2)^3}\\
&=\frac{2(v^1)^5v^3-(v^1)^3(v^2)^2v^3+(v^1)^3(v^3)^3}{((v^1)^2+(v^2)^2+(v^3)^2)^3}\\
&+\frac{2(v^1)^3(v^2)^2v^3-v^1(v^2)^4v^3+v^1(v^2)^2(v^3)^3}{((v^1)^2+(v^2)^2+(v^3)^2)^3}\\
&+\frac{2(v^1)^3(v^3)^3-v^1(v^2)^2(v^3)^3+v^1(v^3)^5}{((v^1)^2+(v^2)^2+(v^3)^2)^3}\\
&-\frac{2(v^1)^5v^3-2(v^1)^3(v^2)^2v^3+2(v^1)^3(v^3)^3}{((v^1)^2+(v^2)^2+(v^3)^2)^3}
\end{aligned}
\right.
$$

$$
(3.3.13)\quad \begin{cases}
g_{13}(\overline{v}) = \frac{1}{2}\frac{(4(v^1)^3v^2+2v^1(v^2)^3-2v^1v^2(v^3)^2)((v^1)^2+(v^2)^2+(v^3)^2)^2}{((v^1)^2+(v^2)^2+(v^3)^2)^4}\\
\quad -\frac{1}{2}\frac{((v^1)^4v^2+(v^1)^2(v^2)^3-(v^1)^2v^2(v^3)^2)2((v^1)^2+(v^2)^2+(v^3)^2)2v^1}{((v^1)^2+(v^2)^2+(v^3)^2)^4}\\
\quad =\frac{(2(v^1)^3v^2+v^1(v^2)^3-v^1v^2(v^3)^2)((v^1)^2+(v^2)^2+(v^3)^2)}{((v^1)^2+(v^2)^2+(v^3)^2)^3}\\
\quad -\frac{((v^1)^4v^2+(v^1)^2(v^2)^3-(v^1)^2v^2(v^3)^2)2v^1}{((v^1)^2+(v^2)^2+(v^3)^2)^3}\\
\quad =\frac{2(v^1)^5v^2+(v^1)^3(v^2)^3-(v^1)^3v^2(v^3)^2}{((v^1)^2+(v^2)^2+(v^3)^2)^3}+\frac{2(v^1)^3(v^2)^3+v^1(v^2)^5-v^1(v^2)^3(v^3)^2}{((v^1)^2+(v^2)^2+(v^3)^2)^3}\\
\quad +\frac{2(v^1)^3v^2(v^3)^2+v^1(v^2)^3(v^3)^2-v^1v^2(v^3)^4}{((v^1)^2+(v^2)^2+(v^3)^2)^3}-\frac{2(v^1)^5v^2+2(v^1)^3(v^2)^3-2(v^1)^3v^2(v^3)^2}{((v^1)^2+(v^2)^2+(v^3)^2)^3}
\end{cases}
$$

$$
(3.3.14)\quad \begin{cases}
g_{22}(\overline{v}) = 1+\frac{1}{2}\frac{(-2(v^1)^2v^2v^3)((v^1)^2+(v^2)^2+(v^3)^2)^2}{((v^1)^2+(v^2)^2+(v^3)^2)^4}\\
\quad -\frac{1}{2}\frac{((v^1)^4v^3-(v^1)^2(v^2)^2v^3+(v^1)^2(v^3)^3)2((v^1)^2+(v^2)^2+(v^3)^2)2v^2}{((v^1)^2+(v^2)^2+(v^3)^2)^4}\\
\quad =1+\frac{-(v^1)^2v^2v^3((v^1)^2+(v^2)^2+(v^3)^2)}{((v^1)^2+(v^2)^2+(v^3)^2)^3}\\
\quad +\frac{(-(v^1)^4v^3+(v^1)^2(v^2)^2v^3-(v^1)^2(v^3)^3)2v^2}{((v^1)^2+(v^2)^2+(v^3)^2)^3}\\
\quad =1+\frac{-(v^1)^4v^2v^3-(v^1)^2(v^2)^3v^3-(v^1)^2v^2(v^3)^3}{((v^1)^2+(v^2)^2+(v^3)^2)^3}\\
\quad +\frac{-2(v^1)^4v^2v^3+2(v^1)^2(v^2)^3v^3-2(v^1)^2v^2(v^3)^3}{((v^1)^2+(v^2)^2+(v^3)^2)^3}
\end{cases}
$$

$$
(3.3.15)\quad
\left\{
\begin{aligned}
g_{23}(\overline{v}) &= \frac{1}{2}\frac{((v^1)^4+3(v^1)^2(v^2)^2-(v^1)^2(v^3)^2)((v^1)^2+(v^2)^2+(v^3)^2)^2}{((v^1)^2+(v^2)^2+(v^3)^2)^4}\\
&-\frac{1}{2}\frac{((v^1)^4v^2+(v^1)^2(v^2)^3-(v^1)^2v^2(v^3)^2)2((v^1)^2+(v^2)^2+(v^3)^2)2v^2}{((v^1)^2+(v^2)^2+(v^3)^2)^4}\\
&=\frac{1}{2}\frac{((v^1)^4+3(v^1)^2(v^2)^2-(v^1)^2(v^3)^2)((v^1)^2+(v^2)^2+(v^3)^2)}{((v^1)^2+(v^2)^2+(v^3)^2)^3}\\
&+\frac{1}{2}\frac{(-(v^1)^4v^2-(v^1)^2(v^2)^3+(v^1)^2v^2(v^3)^2)4v^2}{((v^1)^2+(v^2)^2+(v^3)^2)^3}\\
&=\frac{1}{2}\frac{(v^1)^6+3(v^1)^4(v^2)^2-(v^1)^4(v^3)^2}{((v^1)^2+(v^2)^2+(v^3)^2)^3}\\
&+\frac{1}{2}\frac{(v^1)^4(v^2)^2+3(v^1)^2(v^2)^4-(v^1)^2(v^2)^2(v^3)^2}{((v^1)^2+(v^2)^2+(v^3)^2)^3}\\
&+\frac{1}{2}\frac{(v^1)^4(v^3)^2+3(v^1)^2(v^2)^2(v^3)^2-(v^1)^2(v^3)^4}{((v^1)^2+(v^2)^2+(v^3)^2)^3}\\
&+\frac{1}{2}\frac{-4(v^1)^4(v^2)^2-4(v^1)^2(v^2)^4+4(v^1)^2(v^2)^2(v^3)^2}{((v^1)^2+(v^2)^2+(v^3)^2)^3}
\end{aligned}
\right.
$$

$$
(3.3.16)\quad
\left\{
\begin{aligned}
&g_{33}(\overline{v})\\
&=1+\frac{1}{2}\frac{(-2(v^1)^2v^2v^3)((v^1)^2+(v^2)^2+(v^3)^2)^2}{((v^1)^2+(v^2)^2+(v^3)^2)^4}\\
&-\frac{1}{2}\frac{((v^1)^4v^2+(v^1)^2(v^2)^3-(v^1)^2v^2(v^3)^2)2((v^1)^2+(v^2)^2+(v^3)^2)2v^3}{((v^1)^2+(v^2)^2+(v^3)^2)^4}\\
&=1+\frac{-(v^1)^2v^2v^3((v^1)^2+(v^2)^2+(v^3)^2)}{((v^1)^2+(v^2)^2+(v^3)^2)^4}\\
&+\frac{(-(v^1)^4v^2-(v^1)^2(v^2)^3+(v^1)^2v^2(v^3)^2)2v^3}{((v^1)^2+(v^2)^2+(v^3)^2)^3}\\
&=1+\frac{-(v^1)^4v^2v^3-(v^1)^2(v^2)^3v^3-(v^1)^2v^2(v^3)^3}{((v^1)^2+(v^2)^2+(v^3)^2)^4}\\
&+\frac{-2(v^1)^4v^2v^3-2(v^1)^2(v^2)^3v^3+2(v^1)^2v^2(v^3)^3}{((v^1)^2+(v^2)^2+(v^3)^2)^3}
\end{aligned}
\right.
$$

The equation (3.3.2) follows from the equation (3.3.11). The equation (3.3.3) follows from the equation (3.3.12). The equation (3.3.4) follows from the equation (3.3.13). The equation (3.3.5) follows from the equation (3.3.14). The equation (3.3.6) follows from the equation (3.3.15). The equation (3.3.7) follows from the equation (3.3.16).

□

THEOREM 3.3.2. *The Finslerian metric defined by equation* (3.3.1) *generates symmetric relation of orthogonality.*

Proof. We will start the procedure of orthogonalization[3.3] from the basis

$$\overline{e}_1=\begin{pmatrix}1\\0\\0\end{pmatrix}\quad \overline{e}_2=\begin{pmatrix}0\\1\\0\end{pmatrix}\quad \overline{e}_3=\begin{pmatrix}0\\0\\1\end{pmatrix}$$

Let

$$\overline{e}'_1=\overline{e}_1\quad \overline{e}'_2=\begin{pmatrix}e'^1_2\\e'^2_2\\e'^3_2\end{pmatrix}\quad \overline{e}'_3=\begin{pmatrix}e'^1_3\\e'^2_3\\e'^3_3\end{pmatrix}$$

be orthogonal basis. Since vector $\overline{e}'_2$ is orthogonal to vector $\overline{e}'_1$, then according to the definition 2.3.1

$$g_{ij}(\overline{e}'_1)e'^i_1e'^j_2=0$$

$$g_{1j}(\overline{e}'_1)e'^j_2=0 \tag{3.3.17}$$

According to equations (3.3.2), (3.3.3), (3.3.4),

$$\begin{matrix}g_{11}(\overline{e}'_1)=1 & g_{12}(\overline{e}'_1)=0 & g_{13}(\overline{e}'_1)=0\\ & g_{22}(\overline{e}'_1)=1 & g_{23}(\overline{e}'_1)=\dfrac{1}{2}\\ & & g_{33}(\overline{e}'_1)=1\end{matrix} \tag{3.3.18}$$

From equations (3.3.17), (3.3.18), it follows that

$$e'^1_2=0 \tag{3.3.19}$$

Since values $e'^2_2$, $e'^3_2$ are arbitrary, we assume[3.4]

$$e'^2_2=a\quad e'^3_2=b \tag{3.3.20}$$

According to the theorem 3.3.1 and the equations (3.3.19), (3.3.20),

$$\begin{matrix}g_{11}(\overline{e}'_2)=1+\dfrac{a^5b+2a^3b^3+ab^5}{(a^2+b^2)^3}=1+\dfrac{ab}{a^2+b^2} & g_{12}(\overline{e}'_2)=0 & g_{13}(\overline{e}'_2)=0\\ & g_{22}(\overline{e}'_2)=1 & g_{23}(\overline{e}'_2)=0\\ & & g_{33}(\overline{e}'_2)=1\end{matrix} \tag{3.3.21}$$

From equations (3.3.20), (3.3.21), it follows that

$$g_{ij}(\overline{e}'_2)e'^i_2e'^j_1=g_{i1}(\overline{e}'_2)e'^i_2=0$$

Therefore, the vector $\overline{e}'_1$ is orthogonal to vector $\overline{e}'_2$.

Since vector $\overline{e}'_3$ is orthogonal to vectors $\overline{e}'_1$, $\overline{e}'_2$, then according to the definition 2.3.1

$$\begin{matrix}g_{ij}(\overline{e}'_1)e'^i_1e'^j_3=0 & g_{ij}(\overline{e}'_2)e'^i_2e'^j_3=0\\ g_{1j}(\overline{e}'_1)e'^j_3=0 & g_{2j}(\overline{e}'_2)ae'^j_3+g_{3j}(\overline{e}'_2)be'^j_3=0\end{matrix} \tag{3.3.22}$$

---

[3.3] See section [8]-3.

[3.4] If we assume $a=1$, $b=0$, then it is easy to see that vector $\overline{e}_1$ is orthogonal to vector $\overline{e}_2$ and vector $\overline{e}_2$ is orthogonal to vector $\overline{e}_1$.

From equations (3.3.18), (3.3.21), (3.3.22), it follows that

$$e'^1_3 = 0 \quad ae'^2_3 + be'^3_3 = 0 \tag{3.3.23}$$

From the equation (3.3.23), it follows that we can assume

$$e'^1_3 = 0 \quad e'^2_3 = -b \quad e'^3_3 = a \tag{3.3.24}$$

According to the theorem 3.3.1 and the equation (3.3.24),

$$\begin{aligned} & g_{11}(\overline{e}'_3) \\ &= 1 + \frac{(-b)^5 a + 2(-b)^3(a)^3 + (-b)(a)^5}{(a^2+b^2)^3} \qquad g_{12}(\overline{e}'_3) = 0 \quad g_{13}(\overline{e}'_3) = 0 \\ &= 1 - \frac{ab}{a^2+b^2} \\ & \qquad g_{22}(\overline{e}'_3) = 1 \quad g_{23}(\overline{e}'_3) = 0 \\ & \qquad g_{33}(\overline{e}'_3) = 1 \end{aligned} \tag{3.3.25}$$

From equations (3.3.20), (3.3.24), (3.3.25), it follows that

$$\begin{aligned} g_{ij}(\overline{e}'_3)e'^i_3 e'^j_1 &= g_{i1}(\overline{e}'_3)e'^i_3 = 0 \\ g_{ij}(\overline{e}'_3)e'^i_3 e'^j_2 &= g_{i2}(\overline{e}'_3)e'^i_3 a + g_{i3}(\overline{e}'_3)e'^i_3 b \\ &= g_{22}(\overline{e}'_3)e'^2_3 a + g_{33}(\overline{e}'_3)e'^3_3 b \\ &= (-b)a + ab = 0 \end{aligned}$$

Therefore, the vectors $\overline{e}'_1$, $\overline{e}'_2$ are orthogonal to vector $\overline{e}'_3$. □

Since

$$F^2(\overline{e}'_1) = 1$$

then the length of the vector

$$\overline{e}''_1 = \overline{e}'_1 = \begin{pmatrix} 1 \\ 0 \\ 0 \end{pmatrix}$$

equal 1. Since

$$F^2(\overline{e}'_2) = a^2 + b^2$$

then the length of the vector

$$\overline{e}''_2 = \frac{1}{\sqrt{a^2+b^2}}\overline{e}'_2 = \begin{pmatrix} 0 \\ \frac{a}{\sqrt{a^2+b^2}} \\ \frac{b}{\sqrt{a^2+b^2}} \end{pmatrix}$$

equal 1. Since

$$F^2(\overline{e}'_3) = (-b)^2 + (a)^2 = b^2 + a^2$$

then the length of the vector

$$\overline{e}''_3 = \frac{1}{\sqrt{a^2+b^2}}\overline{e}'_3 = \begin{pmatrix} 0 \\ -\dfrac{b}{\sqrt{a^2+b^2}} \\ \dfrac{a}{\sqrt{a^2+b^2}} \end{pmatrix}$$

equal 1. Therefore, the basis $\overline{\overline{e}}''$ is orthonormal basis. There is passive transformation between bases $\overline{\overline{e}}''$ and $\overline{\overline{e}}$

$$\begin{pmatrix}\overline{e}''_1 & \overline{e}''_2 & \overline{e}''_3\end{pmatrix} = \begin{pmatrix}\overline{e}_1 & \overline{e}_2 & \overline{e}_3\end{pmatrix} \begin{pmatrix} 1 & 0 & 0 \\ 0 & \dfrac{a}{\sqrt{a^2+b^2}} & -\dfrac{b}{\sqrt{a^2+b^2}} \\ 0 & \dfrac{b}{\sqrt{a^2+b^2}} & \dfrac{a}{\sqrt{a^2+b^2}} \end{pmatrix} \tag{3.3.26}$$

Theorem 3.3.3. *Coordinates of metric tensor relative to the basis $\overline{\overline{e}}''$ have form*

$$g''(\overline{e}''_1) = \begin{pmatrix} 1 & 0 & 0 \\ 0 & 1+\dfrac{ab}{a^2+b^2} & \dfrac{1}{2}\dfrac{a^2-b^2}{a^2+b^2} \\ 0 & \dfrac{1}{2}\dfrac{a^2-b^2}{a^2+b^2} & 1-\dfrac{ab}{a^2+b^2} \end{pmatrix} \tag{3.3.27}$$

$$g''(\overline{e}''_2) = \begin{pmatrix} 1+\dfrac{ab}{a^2+b^2} & 0 & 0 \\ 0 & 1 & 0 \\ 0 & 0 & 1 \end{pmatrix} \tag{3.3.28}$$

$$g''(\overline{e}''_3) = \begin{pmatrix} 1-\dfrac{ab}{a^2+b^2} & 0 & 0 \\ 0 & 1 & 0 \\ 0 & 0 & 1 \end{pmatrix} \tag{3.3.29}$$

Proof. From equations (2.5.3), (3.3.18), (3.3.26), it follows that

$$
\left\{
\begin{aligned}
&g''_{11}(\overline{e}''_1)=g_{11}(\overline{e}''_1)=1\\
&g''_{12}(\overline{e}''_1)=g_{12}(\overline{e}''_1)A^2_2+g_{13}(\overline{e}''_1)A^3_2=0\\
&g''_{13}(\overline{e}''_1)=g_{12}(\overline{e}''_1)A^2_3+g_{13}(\overline{e}''_1)A^3_3=0\\
&g''_{22}(\overline{e}''_1)=g_{22}(\overline{e}''_1)A^2_2A^2_2+2g_{23}(\overline{e}''_1)A^2_2A^3_2+g_{33}(\overline{e}''_1)A^3_2A^3_2\\
&\qquad=\frac{a^2}{a^2+b^2}+\frac{ab}{a^2+b^2}+\frac{b^2}{a^2+b^2}=1+\frac{ab}{a^2+b^2}\\
&g''_{23}(\overline{e}''_1)=g_{22}(\overline{e}''_1)A^2_2A^2_3+g_{23}(\overline{e}''_1)A^2_2A^3_3+g_{23}(\overline{e}''_1)A^2_3A^3_2+g_{33}(\overline{e}''_1)A^3_2A^3_3\\
&\qquad=\frac{-ab}{a^2+b^2}+\frac{1}{2}\frac{a^2}{a^2+b^2}-\frac{1}{2}\frac{1}{2}\frac{b^2}{a^2+b^2}+\frac{ab}{a^2+b^2}=\frac{1}{2}\frac{a^2-b^2}{a^2+b^2}\\
&g''_{33}(\overline{e}''_1)=g_{22}(\overline{e}''_1)A^2_3A^2_3+2g_{23}(\overline{e}''_1)A^2_3A^3_3+g_{33}(\overline{e}''_1)A^3_3A^3_3\\
&\qquad=\frac{b^2}{a^2+b^2}-\frac{ab}{a^2+b^2}+\frac{a^2}{a^2+b^2}=1-\frac{ab}{a^2+b^2}
\end{aligned}
\right.
\tag{3.3.30}
$$

The equation (3.3.27) follows from (3.3.30). From equations (2.5.3), (3.3.21), (3.3.26), it follows that

$$
\left\{
\begin{aligned}
&g''_{11}(\overline{e}''_2)=g_{11}(\overline{e}''_2)=1+\frac{ab}{a^2+b^2}\\
&g''_{12}(\overline{e}''_2)=g_{11}(\overline{e}''_2)A^1_2+g_{12}(\overline{e}''_2)A^2_2+g_{13}(\overline{e}''_2)A^3_2=0\\
&g''_{13}(\overline{e}''_2)=g_{11}(\overline{e}''_2)A^1_3+g_{12}(\overline{e}''_2)A^2_3+g_{13}(\overline{e}''_2)A^3_3=0\\
&g''_{22}(\overline{e}''_2)=g_{22}(\overline{e}''_2)A^2_2A^2_2+2g_{23}(\overline{e}''_2)A^2_2A^3_2+g_{33}(\overline{e}''_2)A^3_2A^3_2\\
&\qquad=\frac{a^2}{a^2+b^2}+\frac{b^2}{a^2+b^2}=1\\
&g''_{23}(\overline{e}''_2)=g_{22}(\overline{e}''_2)A^2_2A^2_3+g_{33}(\overline{e}''_2)A^3_2A^3_3=0\\
&g''_{33}(\overline{e}''_2)=g_{22}(\overline{e}''_2)A^2_3A^2_3+2g_{23}(\overline{e}''_2)A^2_3A^3_3+g_{33}(\overline{e}''_2)A^3_3A^3_3\\
&\qquad=\frac{b^2}{a^2+b^2}+\frac{a^2}{a^2+b^2}=1
\end{aligned}
\right.
\tag{3.3.31}
$$

The equation (3.3.28) follows from (3.3.31). From equations (2.5.3), (3.3.25), (3.3.26), it follows that

$$
\left\{
\begin{aligned}
&g''_{11}(\overline{e}''_3)=g_{11}(\overline{e}''_3)=1-\frac{ab}{a^2+b^2}\\
&g''_{12}(\overline{e}''_3)=g_{11}(\overline{e}''_3)A^1_2+g_{12}(\overline{e}''_3)A^2_2+g_{13}(\overline{e}''_3)A^3_2=0\\
&g''_{13}(\overline{e}''_3)=g_{11}(\overline{e}''_3)A^1_3+g_{12}(\overline{e}''_3)A^2_3+g_{13}(\overline{e}''_3)A^3_3=0\\
&g''_{22}(\overline{e}''_3)=g_{22}(\overline{e}''_3)A^2_2A^2_2+2g_{23}(\overline{e}''_3)A^2_2A^3_2+g_{33}(\overline{e}''_3)A^3_2A^3_2\\
&\qquad=\frac{a^2}{a^2+b^2}+\frac{b^2}{a^2+b^2}=1\\
&g''_{23}(\overline{e}''_3)=g_{22}(\overline{e}''_3)A^2_2A^2_3+g_{33}(\overline{e}''_3)A^3_2A^3_3=0\\
&g''_{33}(\overline{e}''_3)=g_{22}(\overline{e}''_3)A^2_3A^2_3+2g_{23}(\overline{e}''_3)A^2_3A^3_3+g_{33}(\overline{e}''_3)A^3_3A^3_3\\
&\qquad=\frac{b^2}{a^2+b^2}+\frac{a^2}{a^2+b^2}=1
\end{aligned}
\right.
\tag{3.3.32}
$$

The equation (3.3.29) follows from (3.3.32). □

Corollary 3.3.4. *The matrix* (2.3.9) *corresponding to the Finslerian metric defined by equation* (3.3.1) *is diagonal matrix.* □

## 3.4. Example 2 of Minkowski Space of Dimension 3

We consider example of the Finslerian metric of Minkowski space of dimension 3 (considered metric is a special case of metric [12]-(3))

$$
\begin{aligned}
F^2(\overline{v}) &= ((v^1)^2+(v^2)^2+(v^3)^2)\left(1+\frac{(v^1)^3v^2}{((v^1)^2+(v^2)^2+(v^3)^2)^2}\right)\\
&= (v^1)^2+(v^2)^2+(v^3)^2+\frac{(v^1)^3v^2}{(v^1)^2+(v^2)^2+(v^3)^2}
\end{aligned}
\tag{3.4.1}
$$

Theorem 3.4.1. *The Finslerian metric $F^2$ defined by equation* (3.4.1) *generates metric of Minkowski space*

$$
\begin{aligned}
g_{11}(\overline{v}) = 1 &- \frac{(v^1)^3(v^2)^3+(v^1)^3v^2(v^3)^2}{((v^1)^2+(v^2)^2+(v^3)^2)^3}\\
&+\frac{6v^1(v^2)^3(v^3)^2+3v^1(v^2)^5+3v^1v^2(v^3)^4}{((v^1)^2+(v^2)^2+(v^3)^2)^3}
\end{aligned}
\tag{3.4.2}
$$

$$
g_{12}(\overline{v}) = \frac{1}{2}\frac{(v^1)^6+6(v^1)^4(v^2)^2-3(v^1)^2(v^2)^4+4(v^1)^4(v^3)^2+3(v^1)^2(v^3)^4}{((v^1)^2+(v^2)^2+(v^3)^2)^3}
\tag{3.4.3}
$$

$$
g_{13}(\overline{v}) = \frac{(v^1)^4v^2v^3-3(v^1)^2(v^2)^3v^3-3(v^1)^2v^2(v^3)^3}{((v^1)^2+(v^2)^2+(v^3)^2)^3}
\tag{3.4.4}
$$

$$
g_{22}(\overline{v}) = 1+\frac{(v^1)^3(v^2)^3-3(v^1)^5v^2-3(v^1)^3v^2(v^3)^2}{((v^1)^2+(v^2)^2+(v^3)^2)^3}
\tag{3.4.5}
$$

$$
g_{23}(\overline{v}) = \frac{3(v^1)^3(v^2)^2v^3-(v^1)^5v^3-(v^1)^3(v^3)^3}{((v^1)^2+(v^2)^2+(v^3)^2)^3}
\tag{3.4.6}
$$

$$
g_{33}(\overline{v}) = 1+\frac{3(v^1)^3v^2(v^3)^2-(v^1)^5v^2-(v^1)^3(v^2)^3}{((v^1)^2+(v^2)^2+(v^3)^2)^3}
\tag{3.4.7}
$$

Proof. From the equation (3.4.1) it follows that

$$
\begin{aligned}
\frac{\partial F^2(\overline{v})}{\partial v^1} &= 2v^1+\frac{3(v^1)^2v^2((v^1)^2+(v^2)^2+(v^3)^2)-(v^1)^3v^22v^1}{((v^1)^2+(v^2)^2+(v^3)^2)^2}\\
&= 2v^1+\frac{3(v^1)^4v^2+3(v^1)^2(v^2)^3+3(v^1)^2v^2(v^3)^2-2(v^1)^4v^2}{((v^1)^2+(v^2)^2+(v^3)^2)^2}\\
&= 2v^1+\frac{(v^1)^4v^2+3(v^1)^2(v^2)^3+3(v^1)^2v^2(v^3)^2}{((v^1)^2+(v^2)^2+(v^3)^2)^2}
\end{aligned}
\tag{3.4.8}
$$

$$
\begin{aligned}
\frac{\partial F^2(\overline{v})}{\partial v^2} &= 2v^2+\frac{(v^1)^3((v^1)^2+(v^2)^2+(v^3)^2)-(v^1)^3v^22v^2}{((v^1)^2+(v^2)^2+(v^3)^2)^2}\\
&= 2v^2+\frac{(v^1)^5+(v^1)^3(v^2)^2+(v^1)^3(v^3)^2-2(v^1)^3(v^2)^2}{((v^1)^2+(v^2)^2+(v^3)^2)^2}\\
&= 2v^2+\frac{(v^1)^5+(v^1)^3(v^3)^2-(v^1)^3(v^2)^2}{((v^1)^2+(v^2)^2+(v^3)^2)^2}
\end{aligned}
\tag{3.4.9}
$$

(3.4.10)
$$\frac{\partial F^2(\overline{v})}{\partial v^3}=2v^3+\frac{-2(v^1)^3v^2v^3}{((v^1)^2+(v^2)^2+(v^3)^2)^2}$$

From equations (2.2.2), (3.4.8), (3.4.9), (3.4.10), it follows that

(3.4.11)
$$\begin{cases} g_{11}(\overline{v}) \\ =1+\frac{1}{2}\frac{(4(v^1)^3v^2+6v^1(v^2)^3+6v^1v^2(v^3)^2)((v^1)^2+(v^2)^2+(v^3)^2)^2}{((v^1)^2+(v^2)^2+(v^3)^2)^4} \\ -\frac{1}{2}\frac{(v^1)^4v^2+3(v^1)^2(v^2)^3+3(v^1)^2v^2(v^3)^2}{((v^1)^2+(v^2)^2+(v^3)^2)^4} \\ *2((v^1)^2+(v^2)^2+(v^3)^2)2v^1 \\ =1+\frac{(2(v^1)^3v^2+3v^1(v^2)^3+3v^1v^2(v^3)^2)((v^1)^2+(v^2)^2+(v^3)^2)}{((v^1)^2+(v^2)^2+(v^3)^2)^3} \\ -\frac{((v^1)^4v^2+3(v^1)^2(v^2)^3+3(v^1)^2v^2(v^3)^2)2v^1}{((v^1)^2+(v^2)^2+(v^3)^2)^3} \\ =1+\frac{2(v^1)^5v^2+3(v^1)^3(v^2)^3+3(v^1)^3v^2(v^3)^2}{((v^1)^2+(v^2)^2+(v^3)^2)^3} \\ +\frac{2(v^1)^3(v^2)^3+3v^1(v^2)^5+3v^1(v^2)^3(v^3)^2}{((v^1)^2+(v^2)^2+(v^3)^2)^3} \\ +\frac{2(v^1)^3v^2(v^3)^2+3v^1(v^2)^3(v^3)^2+3v^1v^2(v^3)^4}{((v^1)^2+(v^2)^2+(v^3)^2)^3} \\ -\frac{2(v^1)^5v^2+6(v^1)^3(v^2)^3+6(v^1)^3v^2(v^3)^2}{((v^1)^2+(v^2)^2+(v^3)^2)^3} \end{cases}$$

(3.4.12)
$$\begin{cases} g_{12}(\overline{v}) \\ =\frac{1}{2}\frac{(5(v^1)^4+3(v^1)^2(v^3)^2-3(v^1)^2(v^2)^2)((v^1)^2+(v^2)^2+(v^3)^2)^2}{((v^1)^2+(v^2)^2+(v^3)^2)^4} \\ -\frac{1}{2}\frac{((v^1)^5+(v^1)^3(v^3)^2-(v^1)^3(v^2)^2)2((v^1)^2+(v^2)^2+(v^3)^2)2v^1}{((v^1)^2+(v^2)^2+(v^3)^2)^4} \\ =\frac{1}{2}\frac{(5(v^1)^4+3(v^1)^2(v^3)^2-3(v^1)^2(v^2)^2)((v^1)^2+(v^2)^2+(v^3)^2)}{((v^1)^2+(v^2)^2+(v^3)^2)^3} \\ -\frac{1}{2}\frac{((v^1)^5+(v^1)^3(v^3)^2-(v^1)^3(v^2)^2)4v^1}{((v^1)^2+(v^2)^2+(v^3)^2)^3} \\ =\frac{1}{2}\frac{5(v^1)^6+3(v^1)^4(v^3)^2-3(v^1)^4(v^2)^2}{((v^1)^2+(v^2)^2+(v^3)^2)^3} \\ +\frac{1}{2}\frac{5(v^1)^4(v^2)^2+3(v^1)^2(v^2)^2(v^3)^2-3(v^1)^2(v^2)^4}{((v^1)^2+(v^2)^2+(v^3)^2)^3} \\ +\frac{1}{2}\frac{5(v^1)^4(v^3)^2+3(v^1)^2(v^3)^4-3(v^1)^2(v^2)^2(v^3)^2}{((v^1)^2+(v^2)^2+(v^3)^2)^3} \\ +\frac{1}{2}\frac{-4(v^1)^6-4(v^1)^4(v^3)^2+4(v^1)^4(v^2)^2}{((v^1)^2+(v^2)^2+(v^3)^2)^3} \end{cases}$$

$$\begin{cases} g_{13}(\overline{v}) = \frac{1}{2}\frac{(-6(v^1)^2v^2v^3)((v^1)^2+(v^2)^2+(v^3)^2)^2}{((v^1)^2+(v^2)^2+(v^3)^2)^4} \\ \quad -\frac{1}{2}\frac{(-2(v^1)^3v^2v^3)2((v^1)^2+(v^2)^2+(v^3)^2)2v^1}{((v^1)^2+(v^2)^2+(v^3)^2)^4} \\ \quad = \frac{(-3(v^1)^2v^2v^3)((v^1)^2+(v^2)^2+(v^3)^2)+((v^1)^3v^2v^3)4v^1}{((v^1)^2+(v^2)^2+(v^3)^2)^3} \\ \quad = \frac{-3(v^1)^4v^2v^3-3(v^1)^2(v^2)^3v^3-3(v^1)^2v^2(v^3)^3)+4(v^1)^4v^2v^3}{((v^1)^2+(v^2)^2+(v^3)^2)^3} \end{cases} \tag{3.4.13}$$

$$\begin{cases} g_{22}(\overline{v}) = 1+\frac{1}{2}\frac{(-2(v^1)^3v^2)((v^1)^2+(v^2)^2+(v^3)^2)^2}{((v^1)^2+(v^2)^2+(v^3)^2)^4} \\ \quad -\frac{1}{2}\frac{((v^1)^5+(v^1)^3(v^3)^2-(v^1)^3(v^2)^2)2((v^1)^2+(v^2)^2+(v^3)^2)2v^2}{((v^1)^2+(v^2)^2+(v^3)^2)^4} \\ \quad = 1+\frac{(-(v^1)^3v^2)((v^1)^2+(v^2)^2+(v^3)^2)}{((v^1)^2+(v^2)^2+(v^3)^2)^3} \\ \quad -\frac{((v^1)^5+(v^1)^3(v^3)^2-(v^1)^3(v^2)^2)2v^2}{((v^1)^2+(v^2)^2+(v^3)^2)^3} \\ \quad = 1+\frac{-(v^1)^5v^2-(v^1)^3(v^2)^3-(v^1)^3v^2(v^3)^2}{((v^1)^2+(v^2)^2+(v^3)^2)^3} \\ \quad +\frac{(-2(v^1)^5v^2-2(v^1)^3v^2(v^3)^2+2(v^1)^3(v^2)^3}{((v^1)^2+(v^2)^2+(v^3)^2)^3} \end{cases} \tag{3.4.14}$$

$$\begin{cases} g_{23}(\overline{v}) \\ = \frac{1}{2}\frac{(-2(v^1)^3v^3)((v^1)^2+(v^2)^2+(v^3)^2)^2}{((v^1)^2+(v^2)^2+(v^3)^2)^4} \\ -\frac{1}{2}\frac{(-2(v^1)^3v^2v^3)2((v^1)^2+(v^2)^2+(v^3)^2)2v^2}{((v^1)^2+(v^2)^2+(v^3)^2)^4} \\ = \frac{-(v^1)^3v^3((v^1)^2+(v^2)^2+(v^3)^2)+2(v^1)^3v^2v^32v^2}{((v^1)^2+(v^2)^2+(v^3)^2)^3} \\ = \frac{-(v^1)^5v^3-(v^1)^3(v^2)^2v^3-(v^1)^3(v^3)^3+4(v^1)^3(v^2)^2v^3}{((v^1)^2+(v^2)^2+(v^3)^2)^3} \end{cases} \tag{3.4.15}$$

$$\begin{cases} g_{33}(\overline{v}) \\ = 1+\frac{1}{2}\frac{(-2(v^1)^3v^2)((v^1)^2+(v^2)^2+(v^3)^2)^2}{((v^1)^2+(v^2)^2+(v^3)^2)^4} \\ -\frac{1}{2}\frac{(-2(v^1)^3v^2v^3)2((v^1)^2+(v^2)^2+(v^3)^2)2v^3}{((v^1)^2+(v^2)^2+(v^3)^2)^4} \\ = 1+\frac{(-(v^1)^3v^2)((v^1)^2+(v^2)^2+(v^3)^2)+(2(v^1)^3v^2v^3)2v^3}{((v^1)^2+(v^2)^2+(v^3)^2)^3} \\ = 1+\frac{-(v^1)^5v^2-(v^1)^3(v^2)^3-(v^1)^3v^2(v^3)^2+4(v^1)^3v^2(v^3)^2}{((v^1)^2+(v^2)^2+(v^3)^2)^3} \end{cases} \tag{3.4.16}$$

The equation (3.4.2) follows from the equation (3.4.11). The equation (3.4.3) follows from the equation (3.4.12). The equation (3.4.4) follows from the equation (3.4.13). The equation (3.4.5) follows from the equation (3.4.14). The equation (3.4.6) follows from the equation (3.4.15). The equation (3.4.7) follows from the equation (3.4.16). □

Theorem 3.4.2. *The Finslerian metric $F^2$ defined by equation (3.4.1) generates asymmetrical relation of orthogonality. However, there exist vectors $\overline{v}$, $\overline{w}$ for which the relation of orthogonality is symmetric.*

Proof. We will start the procedure of orthogonalization[3.5] from the basis

$$\overline{e}_1=\begin{pmatrix}1\\0\\0\end{pmatrix}\quad\overline{e}_2=\begin{pmatrix}0\\1\\0\end{pmatrix}\quad\overline{e}_3=\begin{pmatrix}0\\0\\1\end{pmatrix}$$

Let

$$\overline{e}'_1=\overline{e}_1\quad\overline{e}'_2=\begin{pmatrix}e'^1_2\\e'^2_2\\e'^3_2\end{pmatrix}\quad\overline{e}'_3=\begin{pmatrix}e'^1_3\\e'^2_3\\e'^3_3\end{pmatrix}$$

be orthogonal basis. Since vector $\overline{e}'_2$ is orthogonal to vector $\overline{e}'_1$, then according to the definition 2.3.1

$$g_{ij}(\overline{e}'_1)e'^i_1e'^j_2=0$$

$$g_{1j}(\overline{e}'_1)e'^j_2=0 \tag{3.4.17}$$

According to equations (3.4.2), (3.4.3), (3.4.4),

$$\begin{matrix}g_{11}(\overline{e}'_1)=1 & g_{12}(\overline{e}'_1)=\frac{1}{2} & g_{13}(\overline{e}'_1)=0\\ & g_{22}(\overline{e}'_1)=1 & g_{23}(\overline{e}'_1)=0\\ & & g_{33}(\overline{e}'_1)=1\end{matrix} \tag{3.4.18}$$

From equations (3.4.17), (3.4.18), it follows that

$$e'^1_2+\frac{1}{2}e'^2_2=0 \tag{3.4.19}$$

Since value $e'^3_2$ is arbitrary, we assume[3.6]

$$e'^1_2=-a\quad e'^2_2=2a\quad e'^3_2=b \tag{3.4.20}$$

We define coordinates of vector $\overline{e}'_2$ up to arbitrary factor. Therefore, without loss of generality, we can divide coordinates of vector $\overline{e}'_2$ on $a$. We must consider two cases. We consider procedures of orthogonalization when $a=0$ in subsection 3.4.1. We will see that in this case relation of orthogonality is commutative. We consider procedures of orthogonalization when $a\neq 0$ in subsection 3.4.2. We will see that in this case relation of orthogonality is not commutative. □

[3.5] See section [8]-3.

[3.6] If we assume $a=1$, $b=0$, then it is easy to see that vector $\overline{e}_1$ is orthogonal to vector $\overline{e}_2$ and vector $\overline{e}_2$ is orthogonal to vector $\overline{e}_1$.

**3.4.1. Procedures of Orthogonalization, $(a = 0)$.** Since $b$ is arbitrary, we can assume $b = 1$. So, the vector $\overline{e}'_2$ has form

$$\overline{e}'_2 = \begin{pmatrix} 0 \\ 0 \\ 1 \end{pmatrix} \tag{3.4.21}$$

According to the theorem 3.4.1 and the equation (3.4.21),

$$\begin{array}{ccc} g_{11}(\overline{e}'_2) = 1 & g_{12}(\overline{e}'_2) = 0 & g_{13}(\overline{e}'_2) = 0 \\ & g_{22}(\overline{e}'_2) = 1 & g_{23}(\overline{e}'_2) = 0 \\ & & g_{33}(\overline{e}'_2) = 1 \end{array} \tag{3.4.22}$$

From equations (3.4.20), (3.4.22), it follows that

$$g_{ij}(\overline{e}'_2)e'^i_2 e'^j_1 = g_{22}(\overline{e}'_2)e'^2_1 = 0$$

Therefore, the vector $\overline{e}'_1$ is orthogonal to vector $\overline{e}'_2$.

Since vector $\overline{e}'_3$ is orthogonal to vectors $\overline{e}'_1$, $\overline{e}'_2$, then according to the definition 2.3.1

$$\begin{array}{cc} g_{ij}(\overline{e}'_1)e'^i_1 e'^j_3 = 0 & g_{ij}(\overline{e}'_2)e'^i_2 e'^j_3 = 0 \\ g_{1j}(\overline{e}'_1)e'^j_3 = 0 & g_{3j}(\overline{e}'_2)e'^j_3 = 0 \end{array} \tag{3.4.23}$$

From equations (3.4.18) (3.4.22), (3.4.23), it follows that

$$e'^1_3 + \frac{1}{2}e'^2_3 = 0 \quad e'^3_3 = 0 \tag{3.4.24}$$

From the equation (3.4.24), it follows that we can assume

$$e'^1_3 = -1 \quad e'^2_3 = 2 \quad e'^3_3 = 0 \tag{3.4.25}$$

According to the theorem 3.4.1 and the equation (3.4.25),

$$\begin{aligned} g_{11}(\overline{e}'_3) &= 1 + \frac{-(-1)^3 2^3 + 3(-1)2^5}{((-1)^2 + 2^2)^3} \\ g_{12}(\overline{e}'_3) &= \frac{1}{2}\frac{(-1)^6 + 6(-1)^4 2^2 - 3(-1)^2 2^4}{((-1)^2 + 2^2)^3} \\ g_{13}(\overline{e}'_3) &= 0 \\ g_{22}(\overline{e}'_3) &= 1 + \frac{(-1)^3(2)^3 - 3(-1)^5 2}{((-1)^2 + 2^2)^3} \\ g_{23}(\overline{e}'_3) &= 0 \\ g_{33}(\overline{e}'_3) &= 1 + \frac{-(-1)^5 2 - (-1)^3 2^3}{((-1)^2 + 2^2)^3} \end{aligned} \tag{3.4.26}$$

From equation (3.4.26), it follows that

$$g(\overline{e}'_3) = \begin{pmatrix} \frac{37}{125} & -\frac{23}{250} & 0 \\ & \frac{123}{125} & 0 \\ & & \frac{27}{25} \end{pmatrix} \tag{3.4.27}$$

From equations (3.4.20), (3.4.25), (3.4.27), it follows that

$$g_{ij}(\overline{e}'_3)e'^i_3e'^j_1 = g_{11}(\overline{e}'_3)e'^1_3 + g_{21}(\overline{e}'_3)e'^2_3$$
$$= \frac{37}{125}(-1) - \frac{23}{250}2 = -\frac{12}{25}$$
$$g_{ij}(\overline{e}'_3)e'^i_3e'^j_2 = g_{33}(\overline{e}'_3)e'^3_3 = 0$$

Therefore, the vector $\overline{e}'_2$ is orthogonal to vector $\overline{e}'_3$ and the vector $\overline{e}'_1$ is not orthogonal to vector $\overline{e}'_3$.

Since

$$F^2(\overline{e}'_1) = 1$$

then the length of the vector

$$\overline{e}''_1 = \overline{e}'_1 = \begin{pmatrix}1\\0\\0\end{pmatrix}$$

equal 1. Since

$$F^2(\overline{e}'_2) = 1$$

then the length of the vector

$$\overline{e}''_2 = \overline{e}'_2 = \begin{pmatrix}0\\0\\1\end{pmatrix}$$

equal 1. Since

$$F^2(\overline{e}'_3) = (-1)^2 + 2^2 + \frac{(-1)^3 2}{(-1)^2 + 2^2} = \frac{23}{5}$$

then the length of the vector

$$\overline{e}''_3 = \sqrt{\frac{5}{23}}\overline{e}'_3 = \begin{pmatrix}-\sqrt{\frac{5}{23}}\\2\sqrt{\frac{5}{23}}\\0\end{pmatrix}$$

equal 1. Therefore, the basis $\overline{\overline{e}}''$ is orthonormal basis. There is passive transformation between bases $\overline{\overline{e}}''$ and $\overline{\overline{e}}$

$$\begin{pmatrix}\overline{e}''_1 & \overline{e}''_2 & \overline{e}''_3\end{pmatrix} = \begin{pmatrix}\overline{e}_1 & \overline{e}_2 & \overline{e}_3\end{pmatrix}\begin{pmatrix}1 & 0 & -\sqrt{\frac{5}{23}}\\0 & 0 & 2\sqrt{\frac{5}{23}}\\0 & 1 & 0\end{pmatrix} \tag{3.4.28}$$

Theorem 3.4.3. *Coordinates of metric tensor relative to the basis $\overline{\overline{e}}''$ have form*

$$g''(\overline{e}''_1)=\begin{pmatrix}1&0&0\\0&1&0\\0&0&\frac{20}{23}\end{pmatrix}\tag{3.4.29}$$

$$g''(\overline{e}''_2)=\begin{pmatrix}1&0&-\sqrt{\frac{5}{23}}\\0&1&0\\-\sqrt{\frac{5}{23}}&0&\frac{25}{23}\end{pmatrix}\tag{3.4.30}$$

$$g''(\overline{e}''_3)=\begin{pmatrix}\frac{37}{125}&\frac{123}{125}&-\frac{12}{25}\sqrt{\frac{5}{23}}\\\frac{123}{125}&\frac{27}{25}&0\\-\frac{12}{25}\sqrt{\frac{5}{23}}&0&1\end{pmatrix}\tag{3.4.31}$$

Proof. From equations (2.5.3), (3.4.18), (3.4.28), it follows that

$$\begin{cases}g''_{11}(\overline{e}''_1)=g_{11}(\overline{e}''_1)=1\\g''_{12}(\overline{e}''_1)=g_{13}(\overline{e}''_1)A^1_3=0\\g''_{13}(\overline{e}''_1)=g_{11}(\overline{e}''_1)A^1_3+g_{12}(\overline{e}''_1)A^2_3=-\sqrt{\frac{5}{23}}+\frac{1}{2}2\sqrt{\frac{5}{23}}=0\\g''_{22}(\overline{e}''_1)=g_{33}(\overline{e}''_1)A^3_2A^3_2=1\\g''_{23}(\overline{e}''_1)=g_{13}(\overline{e}''_1)A^3_2A^1_3+g_{23}(\overline{e}''_1)A^3_2A^2_3=0\\g''_{33}(\overline{e}''_1)=g_{11}(\overline{e}''_1)A^1_3A^1_3+2g_{12}(\overline{e}''_1)A^1_3A^2_3+g_{22}(\overline{e}''_1)A^2_3A^2_3\\\qquad=\frac{5}{23}+2*\frac{1}{2}\left(-\sqrt{\frac{5}{23}}\right)2\sqrt{\frac{5}{23}}+4*\frac{5}{23}=\frac{20}{23}\end{cases}\tag{3.4.32}$$

The equation (3.4.29) follows from (3.4.32). From equations (2.5.3), (3.4.22), (3.4.28), it follows that

$$\begin{cases}g''_{11}(\overline{e}''_2)=1\\g''_{12}(\overline{e}''_2)=g_{13}(\overline{e}''_2)A^3_2=0\\g''_{13}(\overline{e}''_2)=g_{11}(\overline{e}''_2)A^1_3+g_{12}(\overline{e}''_2)A^2_3=-\sqrt{\frac{5}{23}}\\g''_{22}(\overline{e}''_2)=g_{33}(\overline{e}''_2)A^3_2A^3_2=1\\g''_{23}(\overline{e}''_2)=g_{31}(\overline{e}''_2)A^3_2A^1_3+g_{32}(\overline{e}''_2)A^3_2A^2_3=0\\g''_{33}(\overline{e}''_2)=g_{11}(\overline{e}''_2)A^1_3A^1_3+2g_{12}(\overline{e}''_2)A^1_3A^2_3+g_{22}(\overline{e}''_2)A^2_3A^2_3\\\qquad=\frac{5}{23}+4*\frac{5}{23}=\frac{25}{23}\end{cases}\tag{3.4.33}$$

The equation (3.4.30) follows from (3.4.33). From equations (2.5.3), (3.4.27), (3.4.28), it follows that

$$
\left\{\begin{array}{l}
g''_{11}(\overline{e}''_3)=g_{11}(\overline{e}''_3)=\frac{37}{125}\\
g''_{12}(\overline{e}''_3)=g_{13}(\overline{e}''_3)A^3_2=\frac{123}{125}\\
g''_{13}(\overline{e}''_3)=g_{11}(\overline{e}''_3)A^1_3+g_{12}(\overline{e}''_3)A^2_3=-\frac{37}{125}\sqrt{\frac{5}{23}}-\frac{23}{250}2\sqrt{\frac{5}{23}}\\
g''_{22}(\overline{e}''_3)=g_{33}(\overline{e}''_3)A^3_2A^3_2\\
\qquad=\frac{27}{25}\\
g''_{23}(\overline{e}''_3)=g_{31}(\overline{e}''_3)A^3_2A^1_3+g_{32}(\overline{e}''_3)A^3_2A^2_3=0\\
g''_{33}(\overline{e}''_3)=g_{11}(\overline{e}''_3)A^1_3A^1_3+2g_{12}(\overline{e}''_3)A^1_3A^2_3+g_{22}(\overline{e}''_3)A^2_3A^2_3\\
\qquad=\frac{37}{125}\frac{5}{23}-2(-\frac{23}{250})2*\frac{5}{23}+\frac{123}{125}4*\frac{5}{23}=1
\end{array}\right.
\tag{3.4.34}
$$

The equation (3.4.31) follows from (3.4.34). □

COROLLARY 3.4.4. *Relative to the basis $\overline{\overline{e}}''$ defined by the equation (3.4.28), the matrix (2.3.9) corresponding to the Finslerian metric $F^2$ defined by equation (3.4.1) has form*

$$
\begin{pmatrix}
1 & 0 & 0\\
0 & 1 & 0\\
-\frac{12}{25}\sqrt{\frac{5}{23}} & 0 & 1
\end{pmatrix}
\tag{3.4.35}
$$

□

**3.4.2. Procedures of Orthogonalization, ($a\neq 0$).** Since $b$ is arbitrary, we can assume $a=1$. So, the vector $\overline{e}'_2$ has form

$$
\overline{e}'_2=\begin{pmatrix}-1\\2\\b\end{pmatrix}
\tag{3.4.36}
$$

According to the theorem 3.4.1 and the equation (3.4.36),

$$
\begin{aligned}
g_{11}(\overline{e}'_2)&=1-\frac{(-1)^3 2^3+(-1)^3 2b^2}{((-1)^2+2^2+b^2)^3}+\frac{6(-1)2^3b^2+3(-1)2^5+3(-1)2b^4}{((-1)^2+2^2+b^2)^3}\\
g_{12}(\overline{e}'_2)&=\frac{1}{2}\frac{(-1)^6+6(-1)^4 2^2-3(-1)^2 2^4+4(-1)^4 b^2+3(-1)^2 b^4}{((-1)^2+2^2+b^2)^3}\\
g_{13}(\overline{e}'_2)&=\frac{(-1)^4 2b-3(-1)^2 2^3 b-3(-1)^2 2b^3}{((-1)^2+2^2+b^2)^3}\\
g_{22}(\overline{e}'_2)&=1+\frac{(-1)^3 2^3-3(-1)^5 2-3(-1)^3 2b^2}{((-1)^2+2^2+b^2)^3}\\
g_{23}(\overline{e}'_2)&=\frac{3(-1)^3 2^2 b-(-1)^5 b-(-1)^3 b^3}{((-1)^2+2^2+b^2)^3}\\
g_{33}(\overline{e}'_2)&=1+\frac{3(-1)^3 2b^2-(-1)^5 2-(-1)^3 2^3}{((-1)^2+2^2+b^2)^3}
\end{aligned}
\tag{3.4.37}
$$

From equations (3.4.37), it follows that

$$g(\overline{e}'_2)=\begin{pmatrix}1-\frac{88+46b^2+6b^4}{(5+b^2)^3}&\frac{-23+4b^2+3b^4}{2(5+b^2)^3}&\frac{-22b-6b^3}{(5+b^2)^3}\\&1+\frac{-2+6b^2}{(5+b^2)^3}&\frac{-11b+b^3}{(5+b^2)^3}\\&&1+\frac{-6b^2+10}{(5+b^2)^3}\end{pmatrix}\tag{3.4.38}$$

From equations (3.4.20), (3.4.38), it follows that

$$\begin{aligned}&g_{ij}(\overline{e}'_2)e'^i_2e'^j_1]=g_{i1}(\overline{e}'_2)e'^i_2\\
=&(1-\frac{88+46b^2+6b^4}{(5+b^2)^3})(-1)+\frac{-23+4b^2+3b^4}{2(5+b^2)^3}2+\frac{-22b-6b^3}{(5+b^2)^3}b\\
=&-1+\frac{88+46b^2+6b^4}{(5+b^2)^3}+\frac{-23+4b^2+3b^4}{(5+b^2)^3}+\frac{-22b^2-6b^4}{(5+b^2)^3}\\
=&-1+\frac{65+28b^2+3b^4}{(5+b^2)^3}=-1+\frac{13+3b^2}{(5+b^2)^2}=\frac{-b^4-10b^2-25+13+3b^2}{(5+b^2)^2}\\
=&\frac{-b^4-7b^2-12}{(5+b^2)^2}\end{aligned}$$

This expression is negative for any value of $b$. Therefore, the vector $\overline{e}'_1$ is not orthogonal to vector $\overline{e}'_2$.

Since vector $\overline{e}'_3$ is orthogonal to vectors $\overline{e}'_1$, $\overline{e}'_2$, then according to the definition 2.3.1

$$\begin{aligned}&g_{ij}(\overline{e}'_1)e'^i_1e'^j_3=0 &&g_{ij}(\overline{e}'_2)e'^i_2e'^j_3=0\\
&g_{1j}(\overline{e}'_1)e'^j_3=0 &&-g_{1j}(\overline{e}'_2)e'^j_3+2g_{2j}(\overline{e}'_2)e'^j_3+bg_{3j}(\overline{e}'_2)e'^j_3=0\end{aligned}\tag{3.4.39}$$

From equations (3.4.18) (3.4.38), (3.4.39), it follows that

$$e'^1_3+\frac{1}{2}e'^2_3=0\tag{3.4.40}$$

$$\begin{aligned}&-(1-\frac{88+46b^2+6b^4}{(5+b^2)^3})e'^1_3-\frac{-23+4b^2+3b^4}{2(5+b^2)^3}e'^2_3-\frac{-22b-6b^3}{(5+b^2)^3}e'^3_3\\
&+2\frac{-23+4b^2+3b^4}{2(5+b^2)^3}e'^1_3+2(1+\frac{-2+6b^2}{(5+b^2)^3})e'^2_3+2\frac{-11b+b^3}{(5+b^2)^3}e'^3_3\\
&+b\frac{-22b-6b^3}{(5+b^2)^3}e'^1_3+b\frac{-11b+b^3}{(5+b^2)^3}e'^2_3+b(1+\frac{-6b^2+10}{(5+b^2)^3})e'^3_3=0\end{aligned}\tag{3.4.41}$$

From the equation (3.4.40), it follows that we can assume

$$e'^1_3=-c\quad e'^2_3=2c\tag{3.4.42}$$

From equations (3.4.41), (3.4.42), it follows that

$$(3.4.43)\quad \begin{aligned}&(1-\frac{88}{5^3})c+\frac{-23+4b^2+3b^4}{2(5+b^2)^3}2c-2\frac{-23+4b^2+3b^4}{2(5+b^2)^3}c\\&+4(1+\frac{-2+6b^2}{(5+b^2)^3})c-b\frac{-22b-6b^3}{(5+b^2)^3}c+2b\frac{-11b+b^3}{(5+b^2)^3}c\\&=\frac{-22b-6b^3}{(5+b^2)^3}e'^3_3-2\frac{-11b+b^3}{(5+b^2)^3}e'^3_3-b(1+\frac{-6b^2+10}{(5+b^2)^3})e'^3_3\end{aligned}$$

Let $b\neq 0$.[3.7] From the equation (3.4.43), it follows that

$$(5+\frac{-50-30b^2-4b^4}{(5+b^2)^3})c=(-b+\frac{-2b^3-10b}{(5+b^2)^3})e'^3_3$$

$$(5+\frac{-4(b^2+5)^2+10(b^2+5)}{(5+b^2)^3})c=-b(1+2\frac{b^2+5}{(5+b^2)^3})e'^3_3$$

$$(5+\frac{-4(b^2+5)+10}{(5+b^2)^2})\frac{c}{b}=-(1+\frac{2}{(5+b^2)^2})e'^3_3$$

$$(3.4.45)\quad (5(5+b^2)^2-4(b^2+5)+10)\frac{c}{b}=-((5+b^2)^2+2)e'^3_3$$

From the equation (3.4.45), it follows that

$$(3.4.46)\quad e'^3_3=-\frac{c}{b}\frac{5(5+b^2)^2-4(b^2+5)+10}{(5+b^2)^2+2}=-\frac{c}{b}(5-\frac{4(b^2+5)}{(5+b^2)^2+2})$$

From equations (3.4.42), (3.4.46), it follows that $c\neq 0$ is arbitrary; we assume $c=1$

$$(3.4.47)\quad \overline{e}'_3=\begin{pmatrix}-1\\2\\-\frac{1}{b}\left(5-\frac{4(b^2+5)}{(5+b^2)^2+2}\right)\end{pmatrix}$$

[3.7] Let $b=0$. Then

$$e'_2=\begin{pmatrix}-1\\2\\0\end{pmatrix}$$

From equations (3.4.43), it follows that

$$(1-\frac{88}{5^3})c-\frac{23}{2*5^3}c-2*\frac{-23}{2*5^3}c+4*(1-\frac{2}{5^3})c=0$$

$$(3.4.44)\quad \left(5-\frac{96}{125}\right)c=0$$

$c=0$ follows from the equation (3.4.40). $e'^1_3=e'^2_3=0$ follows from the equation (3.4.42). Therefore, $e'^3_3\neq 0$ is arbitrary; without loss of generality, we assume

$$\overline{e}'_3=\begin{pmatrix}0\\0\\1\end{pmatrix}$$

In subsection 3.4.1, we proved that relation of orthogonality of vectors $\overline{e}'_2$ and $\overline{e}'_3$ is symmetric. So, the order of vectors $\overline{e}'_2$, $\overline{e}'_3$ in the basis is not important. Since we considered the basis $\overline{\overline{e}}'$ in the subsection 3.4.1, we do not consider the case $b=0$ in this subsection.

Let

$$d=-\frac{1}{b}\left(5-\frac{4(b^2+5)}{(5+b^2)^2+2}\right) \tag{3.4.48}$$

Then

$$\overline{e}'_3=\begin{pmatrix}-1\\2\\d\end{pmatrix} \tag{3.4.49}$$

Vectors $\overline{e}'_2$, $\overline{e}'_3$ have similar coordinates.[3.8] According to equations (3.4.36), (3.4.38), (3.4.49),

$$g(\overline{e}'_3)=\begin{pmatrix}1-\dfrac{88+46d^2+6d^4}{(5+d^2)^3} & \dfrac{-23+4d^2+3d^4}{2(5+d^2)^3} & \dfrac{-22d-6d^3}{(5+d^2)^3}\\ & 1+\dfrac{-2+6d^2}{(5+d^2)^3} & \dfrac{-11d+d^3}{(5+d^2)^3}\\ & & 1+\dfrac{-6d^2+10}{(5+d^2)^3}\end{pmatrix} \tag{3.4.53}$$

Since vectors $\overline{e}'_2$, $\overline{e}'_3$ have similar coordinates, then vector $\overline{e}'_1$ is not orthogonal to vector $\overline{e}'_3$.

Since

$$F^2(\overline{e}'_1)=1$$

then the length of the vector

$$\overline{e}''_1=\overline{e}'_1=\begin{pmatrix}1\\0\\0\end{pmatrix}$$

equal 1. Since

$$F^2(\overline{e}'_2)=(-1)^2+2^2+b^2+\frac{(-1)^3 2}{(-1)^2+2^2+b^2}=5+b^2-\frac{2}{5+b^2}$$

---

[3.8]The condition

$$d=b \tag{3.4.50}$$

means equality of vectors $\overline{e}'_3$, $\overline{e}'_3$. However from equations (3.4.48), (3.4.50), it follows that

$$b^2+5=\frac{4(b^2+5)}{(5+b^2)^2+2} \tag{3.4.51}$$

Since

$$\frac{4}{(5+b^2)^2+2}<1 \tag{3.4.52}$$

then the equation (3.4.50) is impossible.

then the length of the vector

$$\overline{e}''_2=\sqrt{\frac{5+b^2}{(5+b^2)^2-2}}\overline{e}'_2=\begin{pmatrix}-\sqrt{\frac{5+b^2}{(5+b^2)^2-2}}\\ 2\sqrt{\frac{5+b^2}{(5+b^2)^2-2}}\\ b\sqrt{\frac{5+b^2}{(5+b^2)^2-2}}\end{pmatrix}$$

equal 1. Since

$$F^2(\overline{e}'_3)=(-1)^2+2^2+d^2+\frac{(-1)^3 2}{(-1)^2+2^2+d^2}=5+d^2-\frac{2}{5+d^2}$$

then the length of the vector

$$\overline{e}''_3=\sqrt{\frac{5+d^2}{(5+d^2)^2-2}}\overline{e}'_3=\begin{pmatrix}-\sqrt{\frac{5+d^2}{(5+d^2)^2-2}}\\ 2\sqrt{\frac{5+d^2}{(5+d^2)^2-2}}\\ d\sqrt{\frac{5+d^2}{(5+d^2)^2-2}}\end{pmatrix}$$

equal 1. Therefore, the basis $\overline{\overline{e}}''$ is orthonormal basis.

Let

$$f_b=\sqrt{\frac{5+b^2}{(5+b^2)^2-2}} \tag{3.4.54}$$

$$f_d=\sqrt{\frac{5+d^2}{(5+d^2)^2-2}} \tag{3.4.55}$$

There is passive transformation between bases $\overline{\overline{e}}''$ and $\overline{\overline{e}}$

$$\begin{pmatrix}\overline{e}''_1 & \overline{e}''_2 & \overline{e}''_3\end{pmatrix}=\begin{pmatrix}\overline{e}_1 & \overline{e}_2 & \overline{e}_3\end{pmatrix}\begin{pmatrix}1 & -f_b & -f_d\\ 0 & 2f_b & 2f_d\\ 0 & bf_b & df_d\end{pmatrix} \tag{3.4.56}$$

THEOREM 3.4.5. *Consider Minkowski space with the Finslerian metric defined by equation* (3.4.1). *Let vector $\overline{v}$ has coordinates*

$$\overline{v}=\begin{pmatrix}-\alpha\\ 2\alpha\\ \gamma\end{pmatrix} \tag{3.4.57}$$

*Coordinates of $g(\overline{v})$ have form*

(3.4.58)
$$\left\{\begin{aligned}
g_{11}(\overline{v})&=1-\frac{88\alpha^6+46\alpha^4\gamma^2+6\alpha^2\gamma^4}{(5\alpha^2+\gamma^2)^3}\\
g_{12}(\overline{v})&=\frac{-23\alpha^6+4\alpha^4\gamma^2+3\alpha^2\gamma^4}{2(5\alpha^2+\gamma^2)^3}\\
g_{13}(\overline{v})&=\frac{-22\alpha^5\gamma-6\alpha^3\gamma^3}{(5\alpha^2+\gamma^2)^3}\\
g_{22}(\overline{v})&=1+\frac{-2\alpha^6+6\alpha^4\gamma^2}{(5\alpha^2+\gamma^2)^3}\\
g_{23}(\overline{v})&=\frac{-11\alpha^5\gamma+\alpha^3\gamma^3}{(5\alpha^2+\gamma^2)^3}\\
g_{33}(\overline{v})&=1+\frac{-6\alpha^4\gamma^2+10\alpha^6}{(5\alpha^2+\gamma^2)^3}
\end{aligned}\right.$$

Proof. According to the theorem 3.4.1 and the equation (3.4.57),

(3.4.59)
$$\left\{\begin{aligned}
g_{11}(\overline{v})&=1+\frac{6(-\alpha)(2\alpha)^3\gamma^2+3(-\alpha)2\alpha\gamma^4}{(5\alpha^2+\gamma^2)^3}\\
&+\frac{3(-\alpha)(2\alpha)^5-(-\alpha)^3(2\alpha)^3-(-\alpha)^32\alpha\gamma^2}{(5\alpha^2+\gamma^2)^3}\\
&=1+\frac{-48\alpha^4\gamma^2-6\alpha^2\gamma^4-96\alpha^6+8\alpha^6+2\alpha^4\gamma^2}{(5\alpha^2+\gamma^2)^3}\\
g_{12}(\overline{v})&=\frac{3(-\alpha)^4(2\alpha)^2+2(-\alpha)^4\gamma^2}{(5\alpha^2+\gamma^2)^3}\\
&+\frac{1}{2}\frac{\alpha^6-3(-\alpha)^2(2\alpha)^4+3(-\alpha)^2\gamma^4}{(5\alpha^2+\gamma^2)^3}\\
g_{13}(\overline{v})&=\frac{(-\alpha)^42\alpha\gamma-3(-\alpha)^2(2\alpha)^3\gamma-3(-\alpha)^22\alpha\gamma^3}{(5\alpha^2+\gamma^2)^3}\\
&=\frac{2\alpha^5\gamma-24\alpha^5\gamma-6\alpha^3\gamma^3}{(5\alpha^2+\gamma^2)^3}\\
g_{22}(\overline{v})&=1+\frac{(-\alpha)^3(2\alpha)^3-3(-\alpha)^52\alpha-3(-\alpha)^32\alpha\gamma^2}{(5\alpha^2+\gamma^2)^3}\\
&=1+\frac{-8\alpha^6+6\alpha^6+6\alpha^4\gamma^2}{(5\alpha^2+\gamma^2)^3}\\
g_{23}(\overline{v})&=\frac{3(-\alpha)^3(2\alpha)^2\gamma-(-\alpha)^5\gamma-(-\alpha)^3\gamma^3}{(5\alpha^2+\gamma^2)^3}\\
&=\frac{-12\alpha^5\gamma+\alpha^5\gamma+\alpha^3\gamma^3}{(5\alpha^2+\gamma^2)^3}\\
g_{33}(\overline{v})&=1+\frac{3(-\alpha)^32\alpha\gamma^2-(-\alpha)^52\alpha-(-\alpha)^3(2\alpha)^3}{(5\alpha^2+\gamma^2)^3}
\end{aligned}\right.$$

Equations (3.4.58) follows from equations (3.4.59). □

Theorem 3.4.6. *Coordinates of metric tensor relative to the basis $\overline{\overline{e}}''$ have form*

(3.4.60)
$$g''(\overline{e}''_1)=\begin{pmatrix}1&0&0\\0&f_b^2(3+b^2)&f_bf_d(3+bd)\\0&f_bf_d(3+bd)&f_d^2(3+d^2)\end{pmatrix}$$

$$g''(\overline{e}''_2) =$$

(3.4.61)
$$= \begin{pmatrix} 1-\frac{88+46b^2+6b^4}{(5+b^2)^3} & -\frac{(b^2+3)(b^2+4)}{(5+b^2)^2}f_b & \boxed{\left(\frac{-b^4-b^2+10}{(5+b^2)^2} - \frac{88+24b^2}{(5+b^2)^2((5+b^2)^2+2)}\right)f_d} \\ -\frac{(b^2+3)(b^2+4)}{(5+b^2)^2}f_b & 1 & 0 \\ \boxed{\left(\frac{-b^4-b^2+10}{(5+b^2)^2} - \frac{88+24b^2}{(5+b^2)^2((5+b^2)^2+2)}\right)f_d} & 0 & \boxed{\left(5+\left(1+\frac{-6b^2+10}{(5+b^2)^3}\right)d^2 - \frac{50+110b^2+12b^4}{(5+b^2)^3} + \frac{64b^2}{(5+b^2)^2\,((5+b^2)^2+2)}\right)f_d^2} \end{pmatrix}$$

$$g''(\overline{e}''_3) =$$

(3.4.62)
$$= \begin{pmatrix} 1-\frac{88+46d^2+6d^4}{(5+d^2)^3} & \boxed{f_b\left(-1+\frac{9d^2+35}{(5+d^2)^2} - \frac{8(11+3d^2)(5+b^2)}{(5+d^2)^3((5+b^2)^2+2)}\right)} & -\frac{(d^2+3)(d^2+4)}{(5+d^2)^2}f_d \\ \boxed{f_b\left(-1+\frac{9d^2+35}{(5+d^2)^2} - \frac{8(11+3d^2)(5+b^2)}{(5+d^2)^3((5+b^2)^2+2)}\right)} & \boxed{f_b^2\left(5+\left(1+\frac{-6d^2+10}{(5+d^2)^3}\right)b^2 - \frac{50+110d^2+12d^4}{(5+d^2)^3} + \frac{64d^2(5+b^2)}{(5+d^2)^3((5+b^2)^2+2)}\right)} & \boxed{f_bf_d\left(-\frac{4}{(5+d^2)^2} + \frac{4(5+b^2)}{(5+b^2)^2+2} * \left(1+\frac{2}{(5+d^2)^2}\right)\right)} \\ -\frac{(d^2+3)(d^2+4)}{(5+d^2)^2}f_d & \boxed{f_bf_d\left(-\frac{4}{(5+d^2)^2} + \frac{4(5+b^2)}{(5+b^2)^2+2} * \left(1+\frac{2}{(5+d^2)^2}\right)\right)} & 1 \end{pmatrix}$$

Proof. From equations (2.5.3), (3.4.18), (3.4.56), it follows that

$$
\left\{
\begin{aligned}
&g''_{11}(\overline{e}''_1) = g_{11}(\overline{e}''_1) = 1\\
&g''_{12}(\overline{e}''_1) = g_{11}(\overline{e}''_1)A^1_2 + g_{12}(\overline{e}''_1)A^2_2 = 0\\
&g''_{13}(\overline{e}''_1) = g_{11}(\overline{e}''_1)A^1_3 + g_{12}(\overline{e}''_1)A^2_3 = 0\\
&g''_{22}(\overline{e}''_1) = g_{11}(\overline{e}''_1)A^1_2A^1_2 + 2g_{12}(\overline{e}''_1)A^1_2A^2_2\\
&\qquad + g_{22}(\overline{e}''_1)A^2_2A^2_2 + g_{33}(\overline{e}''_1)A^3_2A^3_2\\
&\qquad = f_b^2 - 2f_b^2 + 4f_b^2 + b^2f_b^2 = f_b^2(3+b^2)\\
&g''_{23}(\overline{e}''_1) = g_{11}(\overline{e}''_1)A^1_2A^1_3 + g_{12}(\overline{e}''_1)A^1_2A^2_3 + g_{12}(\overline{e}''_1)A^1_3A^2_2\\
&\qquad + g_{22}(\overline{e}''_1)A^2_2A^2_3 + g_{33}(\overline{e}''_1)A^3_2A^3_3\\
&\qquad = (-f_b)(-f_d) + \frac{1}{2}(-f_b)2f_d + \frac{1}{2}2f_b(-f_d) + 4f_bf_d + bf_bdf_d\\
&\qquad = f_bf_d(3+bd)\\
&g''_{33}(\overline{e}''_1) = g_{11}(\overline{e}''_1)A^1_3A^1_3 + 2g_{12}(\overline{e}''_1)A^1_3A^2_3\\
&\qquad + g_{22}(\overline{e}''_1)A^2_3A^2_3 + g_{33}(\overline{e}''_1)A^3_3A^3_3\\
&\qquad = f_d^2 - 2f_d^2 + 4f_d^2 + d^2f_d^2 = f_d^2(3+d^2)
\end{aligned}
\right.
\tag{3.4.63}
$$

The equation (3.4.60) follows from (3.4.63). From equations (2.5.3), (3.4.38), (3.4.56), it follows that

$$
(3.4.64)\quad
\left\{
\begin{aligned}
g''_{11}(\overline{e}''_2) &= g_{11}(\overline{e}''_2) = 1-\frac{88+46b^2+6b^4}{(5+b^2)^3}\\
g''_{12}(\overline{e}''_2) &= g_{11}(\overline{e}''_2)A^1_2+g_{12}(\overline{e}''_2)A^2_2+g_{13}(\overline{e}''_2)A^3_2\\
&=\left(1-\frac{88+46b^2+6b^4}{(5+b^2)^3}\right)(-f_b)+\frac{-23+4b^2+3b^4}{2(5+b^2)^3}2f_b+\frac{-22b-6b^3}{(5+b^2)^3}bf_b\\
&=\left(-1+\frac{65+28b^2+3b^4}{(5+b^2)^3}\right)f_b=\left(-1+\frac{3b^2+13}{(5+b^2)^2}\right)f_b\\
&=\left(\frac{-b^4-7b^2-12}{(5+b^2)^2}\right)f_b=-\frac{(b^2+3)(b^2+4)}{(5+b^2)^2}f_b\\
g''_{13}(\overline{e}''_2) &= g_{11}(\overline{e}''_2)A^1_3+g_{12}(\overline{e}''_2)A^2_3+g_{13}(\overline{e}''_2)A^3_3\\
&=\left(1-\frac{88+46b^2+6b^4}{(5+b^2)^3}\right)(-f_d)+\frac{-23+4b^2+3b^4}{2(5+b^2)^3}2f_d+\frac{-22b-6b^3}{(5+b^2)^3}df_d\\
&=\left(-1+\frac{65+50b^2+9b^4}{(5+b^2)^3}-\frac{22+6b^2}{(5+b^2)^3}bd\right)f_d\\
&=\left(-1+\frac{65+50b^2+9b^4}{(5+b^2)^3}-\frac{22+6b^2}{(5+b^2)^3}\left(\frac{4(5+b^2)}{(5+b^2)^2+2}-5\right)\right)f_d\\
&=\left(-1+\frac{175+80b^2+9b^4}{(5+b^2)^3}-\frac{22+6b^2}{(5+b^2)^3}\frac{4(5+b^2)}{(5+b^2)^2+2}\right)f_d\\
&=\left(-1+\frac{9b^2+35}{(5+b^2)^2}-\frac{22+6b^2}{(5+b^2)^2}\frac{4}{(5+b^2)^2+2}\right)f_d\\
&=\left(\frac{-b^4-b^2+10}{(5+b^2)^2}-\frac{88+24b^2}{(5+b^2)^2((5+b^2)^2+2)}\right)f_d\\
g''_{22}(\overline{e}''_2) &= g_{11}(\overline{e}''_2)A^1_2A^1_2+2g_{12}(\overline{e}''_2)A^1_2A^2_2+2g_{13}(\overline{e}''_2)A^1_2A^3_2\\
&\quad+g_{22}(\overline{e}''_2)A^2_2A^2_2+2g_{23}(\overline{e}''_2)A^2_2A^3_2+g_{33}(\overline{e}''_2)A^3_2A^3_2\\
&=\left(1-\frac{88+46b^2+6b^4}{(5+b^2)^3}\right)(-f_b)(-f_b)\\
&\quad+2\frac{-23+4b^2+3b^4}{2(5+b^2)^3}(-f_b)2f_b\\
&\quad+2\frac{-22b-6b^3}{(5+b^2)^3}(-f_b)bf_b+\left(1+\frac{-2+6b^2}{(5+b^2)^3}\right)2f_b2f_b\\
&\quad+2\frac{-11b+b^3}{(5+b^2)^3}2f_bbf_b+\left(1+\frac{-6b^2+10}{(5+b^2)^3}\right)bf_bbf_b\\
&=(5+b^2-2\frac{25+10b^2+b^4}{(5+b^2)^3})f_b^2=1
\end{aligned}
\right.
$$

$$
\begin{aligned}
g''_{23}(\overline{e}''_2) &= g_{11}(\overline{e}''_2)A^1_2A^1_3 + g_{12}(\overline{e}''_2)A^1_2A^2_3 + g_{12}(\overline{e}''_2)A^1_3A^2_2 \\
&+g_{13}(\overline{e}''_2)A^1_2A^3_3 + g_{13}(\overline{e}''_2)A^1_3A^3_2 \\
&+g_{22}(\overline{e}''_2)A^2_2A^2_3 + g_{23}(\overline{e}''_2)A^2_2A^3_3 + g_{23}(\overline{e}''_2)A^2_3A^3_2 + g_{33}(\overline{e}''_2)A^3_2A^3_3 \\
&= \left(1 - \frac{88+46b^2+6b^4}{(5+b^2)^3}\right)(-f_b)(-f_d) \\
&+\frac{-23+4b^2+3b^4}{2(5+b^2)^3}(-f_b)2f_d + \frac{-23+4b^2+3b^4}{2(5+b^2)^3}(-f_d)2f_b \\
&+\frac{-22b-6b^3}{(5+b^2)^3}(-f_b)df_d + \frac{-22b-6b^3}{(5+b^2)^3}(-f_d)bf_b \\
&+\left(1+\frac{-2+6b^2}{(5+b^2)^3}\right)2f_b2f_d + \frac{-11b+b^3}{(5+b^2)^3}2f_bdf_d \\
&+\frac{-11b+b^3}{(5+b^2)^3}2f_dbf_b + \left(1+\frac{-6b^2+10}{(5+b^2)^3}\right)bf_ddf_d \\
&= \left(5+bd-\frac{50+30b^2+4b^4}{(5+b^2)^3}\right)f_bf_d + \frac{10+2b^2}{(5+b^2)^3}bdf_bf_d \qquad (3.4.65) \\
&= \left(\frac{4(5+b^2)}{(5+b^2)^2+2} - \frac{50+30b^2+4b^4}{(5+b^2)^3}\right. \\
&\left.+\frac{2}{(5+b^2)^2}\left(\frac{4(5+b^2)}{(5+b^2)^2+2}-5\right)\right)f_bf_d \\
&= \left(\frac{4(5+b^2)}{(5+b^2)^2+2} - \frac{50+30b^2+4b^4}{(5+b^2)^3} - \frac{10}{(5+b^2)^2}\right. \\
&\left.+\frac{8}{(5+b^2)\left((5+b^2)^2+2\right)}\right)f_bf_d \\
&= \left(\frac{4(5+b^2)}{(5+b^2)^2+2} - \frac{100+40b^2+4b^4}{(5+b^2)^3}\right. \\
&\left.+\frac{8}{(5+b^2)\left((5+b^2)^2+2\right)}\right)f_bf_d \\
&= \left(\frac{4(5+b^2)}{(5+b^2)^2+2} - \frac{4}{5+b^2} + \frac{8}{(5+b^2)\left((5+b^2)^2+2\right)}\right)f_bf_d = 0
\end{aligned}
$$

$$
\begin{aligned}
g''_{33}(\overline{e}''_2) &= g_{11}(\overline{e}''_2)A_3^1A_3^1 + 2g_{12}(\overline{e}''_2)A_3^1A_3^2 + 2g_{13}(\overline{e}''_2)A_3^1A_3^3 \\
&+g_{22}(\overline{e}''_2)A_3^2A_3^2 + 2g_{23}(\overline{e}''_2)A_3^2A_3^3 + g_{33}(\overline{e}''_2)A_3^3A_3^3 \\
&= \left(1 - \frac{88+46b^2+6b^4}{(5+b^2)^3}\right)(-f_d)(-f_d) \\
&+2\frac{-23+4b^2+3b^4}{2(5+b^2)^3}(-f_d)2f_d \\
&+2\frac{-22b-6b^3}{(5+b^2)^3}(-f_d)df_d + \left(1+\frac{-2+6b^2}{(5+b^2)^3}\right)2f_d2f_d \\
&+2\frac{-11b+b^3}{(5+b^2)^3}2f_ddf_d + \left(1+\frac{-6b^2+10}{(5+b^2)^3}\right)df_ddf_d \\
&= \left(5+\left(1+\frac{-6b^2+10}{(5+b^2)^3}\right)d^2 - \frac{50+30b^2+12b^4}{(5+b^2)^3}\right. \\
&\left.+\frac{16b^3d}{(5+b^2)^3}\right)f_d^2 \\
&= \left(5+\left(1+\frac{-6b^2+10}{(5+b^2)^3}\right)d^2 - \frac{50+30b^2+12b^4}{(5+b^2)^3}\right. \\
&\left.+\frac{16b^2}{(5+b^2)^3}\left(\frac{4(5+b^2)}{(5+b^2)^2+2}-5\right)\right)f_d^2 \\
&= \left(5+\left(1+\frac{-6b^2+10}{(5+b^2)^3}\right)d^2 - \frac{50+30b^2+12b^4}{(5+b^2)^3} - \frac{80b^2}{(5+b^2)^3}\right. \\
&\left.+\frac{64b^2}{(5+b^2)^2\left((5+b^2)^2+2\right)}\right)f_d^2 \\
&= \left(5+\left(1+\frac{-6b^2+10}{(5+b^2)^3}\right)d^2 - \frac{50+110b^2+12b^4}{(5+b^2)^3} + \frac{64b^2}{(5+b^2)^2\left((5+b^2)^2+2\right)}\right)f_d^2
\end{aligned}
\tag{3.4.66}
$$

The equation (3.4.61) follows from (3.4.64), (3.4.65), (3.4.66). From equations (2.5.3), (3.4.53), (3.4.56), it follows that

$$
(3.4.67)\quad
\left\{
\begin{aligned}
g''_{11}(\overline{e}''_3) &= g_{11}(\overline{e}''_3) = 1-\frac{88+46d^2+6d^4}{(5+d^2)^3}\\
g''_{12}(\overline{e}''_3) &= g_{11}(\overline{e}''_3)A_2^1+g_{12}(\overline{e}''_3)A_2^2+g_{13}(\overline{e}''_3)A_2^3\\
&=\left(1-\frac{88+46d^2+6d^4}{(5+d^2)^3}\right)(-f_b)+\frac{-23+4d^2+3d^4}{2(5+d^2)^3}2f_b+\frac{-22d-6d^3}{(5+d^2)^3}bf_b\\
&=\left(-1+\frac{65+50d^2+9d^4}{(5+d^2)^3}-\frac{22+6d^2}{(5+d^2)^3}bd\right)f_b\\
&=\left(-1+\frac{65+50d^2+9d^4}{(5+d^2)^3}-\frac{22+6d^2}{(5+d^2)^3}\left(\frac{4(5+b^2)}{(5+b^2)^2+2}-5\right)\right)f_b\\
&=\left(-1+\frac{175+80d^2+9d^4}{(5+d^2)^3}-\frac{8(11+3d^2)(5+b^2)}{(5+d^2)^3((5+b^2)^2+2)}\right)f_b\\
&=\left(-1+\frac{9d^2+35}{(5+d^2)^2}-\frac{8(11+3d^2)(5+b^2)}{(5+d^2)^3((5+b^2)^2+2)}\right)f_b\\
g''_{13}(\overline{e}''_3) &= g_{11}(\overline{e}''_3)A_3^1+g_{12}(\overline{e}''_3)A_3^2+g_{13}(\overline{e}''_3)A_3^3\\
&=\left(1-\frac{88+46d^2+6d^4}{(5+d^2)^3}\right)(-f_d)+\frac{-23+4d^2+3d^4}{2(5+d^2)^3}2f_d\\
&+\frac{-22d-6d^3}{(5+d^2)^3}df_d\\
&=\left(-1+\frac{65+28d^2+3d^4}{(5+d^2)^3}\right)f_d=\left(-1+\frac{3d^2+13}{(5+d^2)^2}\right)f_d\\
&=\left(\frac{-d^4-7d^2-12}{(5+d^2)^2}\right)f_d=-\frac{(d^2+3)(d^2+4)}{(5+d^2)^2}f_d\\
g''_{22}(\overline{e}''_3) &= g_{11}(\overline{e}''_3)A_2^1A_2^1+2g_{12}(\overline{e}''_3)A_2^1A_2^2+2g_{13}(\overline{e}''_3)A_2^1A_2^3\\
&+g_{22}(\overline{e}''_3)A_2^2A_2^2+2g_{23}(\overline{e}''_3)A_2^2A_2^3+g_{33}(\overline{e}''_3)A_2^3A_2^3\\
&=\left(1-\frac{88+46d^2+6d^4}{(5+d^2)^3}\right)(-f_b)(-f_b)\\
&+2\frac{-23+4d^2+3d^4}{2(5+d^2)^3}(-f_b)2f_b\\
&+2\frac{-22d-6d^3}{(5+d^2)^3}(-f_b)bf_b+\left(1+\frac{-2+6d^2}{(5+d^2)^3}\right)2f_b2f_b\\
&+2\frac{-11d+d^3}{(5+d^2)^3}2f_bbf_b+\left(1+\frac{-6d^2+10}{(5+d^2)^3}\right)bf_bbf_b\\
&=\left(5+\left(1+\frac{-6d^2+10}{(5+d^2)^3}\right)b^2-\frac{50+30d^2+12d^4}{(5+d^2)^3}+\frac{16d^2}{(5+d^2)^3}bd\right)f_b^2\\
&=\left(5+\left(1+\frac{-6d^2+10}{(5+d^2)^3}\right)b^2-\frac{50+30d^2+12d^4}{(5+d^2)^3}+\frac{16d^2}{(5+d^2)^3}\left(\frac{4(5+b^2)}{(5+b^2)^2+2}-5\right)\right)f_b^2\\
&=\left(5+\left(1+\frac{-6d^2+10}{(5+d^2)^3}\right)b^2-\frac{50+30d^2+12d^4}{(5+d^2)^3}-\frac{80d^2}{(5+d^2)^3}\right.\\
&\left.+\frac{64d^2(5+b^2)}{(5+d^2)^3((5+b^2)^2+2)}\right)f_b^2\\
&=\left(5+\left(1+\frac{-6d^2+10}{(5+d^2)^3}\right)b^2-\frac{50+110d^2+12d^4}{(5+d^2)^3}+\frac{64d^2(5+b^2)}{(5+d^2)^3((5+b^2)^2+2)}\right)f_b^2
\end{aligned}
\right.
$$

(3.4.68)
$$\left\{\begin{aligned}
g''_{23}(\overline{e}''_3) &= g_{11}(\overline{e}''_3)A_2^1A_3^1 + g_{12}(\overline{e}''_3)A_2^1A_3^2 + g_{12}(\overline{e}''_3)A_3^1A_2^2\\
&+g_{13}(\overline{e}''_3)A_2^1A_3^3 + g_{13}(\overline{e}''_3)A_3^1A_2^3\\
&+g_{22}(\overline{e}''_3)A_2^2A_3^2 + g_{23}(\overline{e}''_3)A_2^2A_3^3 + g_{23}(\overline{e}''_3)A_3^2A_2^3 + g_{33}(\overline{e}''_3)A_3^3A_3^3\\
&= \left(1 - \frac{88+46d^2+6d^4}{(5+d^2)^3}\right)(-f_b)(-f_d)\\
&+\frac{-23+4d^2+3d^4}{2(5+d^2)^3}(-f_b)2f_d + \frac{-23+4d^2+3d^4}{2(5+d^2)^3}(-f_d)2f_b\\
&+\frac{-22d-6d^3}{(5+d^2)^3}(-f_b)df_d + \frac{-22d-6d^3}{(5+d^2)^3}(-f_d)bf_b\\
&+\left(1+\frac{-2+6d^2}{(5+d^2)^3}\right)2f_b2f_d + \frac{-11d+d^3}{(5+d^2)^3}2f_bdf_d\\
&+\frac{-11d+d^3}{(5+d^2)^3}2f_dbf_b + \left(1+\frac{-6d^2+10}{(5+d^2)^3}\right)df_ddf_d\\
&= \left(\left(5+bd-\frac{50+30d^2+4d^4}{(5+d^2)^3}\right)+\frac{10+2d^2}{(5+d^2)^3}bd\right)f_bf_d\\
&= \left(-\frac{50+30d^2+4d^4}{(5+d^2)^3}+\frac{4(5+b^2)}{(5+b^2)^2+2}+\frac{2(5+d^2)}{(5+d^2)^3}\left(\frac{4(5+b^2)}{(5+b^2)^2+2}-5\right)\right)f_bf_d\\
&= \left(-\frac{50+30d^2+4d^4}{(5+d^2)^3}+\frac{4\,(5+b^2)}{(5+b^2)^2+2}+\frac{2}{(5+d^2)^2}\left(\frac{4\,(5+b^2)}{(5+b^2)^2+2}-5\right)\right)f_bf_d\\
&= \left(-\frac{50+30d^2+4d^4}{(5+d^2)^3}+\frac{4(5+b^2)}{(5+b^2)^2+2}-\frac{10}{(5+d^2)^2}+\frac{8(5+b^2)}{(5+d^2)^2((5+b^2)^2+2)}\right)f_bf_d\\
&= \left(-\frac{50+30d^2+4d^4}{(5+d^2)^3}-\frac{10(5+d^2)}{(5+d^2)^3}+\frac{4(5+b^2)}{(5+b^2)^2+2}\left(1+\frac{2}{(5+d^2)^2}\right)\right)f_bf_d\\
&= \left(-\frac{100+40d^2+4d^4}{(5+d^2)^3}+\frac{4(5+b^2)}{(5+b^2)^2+2}\left(1+\frac{2}{(5+d^2)^2}\right)\right)f_bf_d\\
&= \left(-\frac{4(5+d^2)^2}{(5+d^2)^3}+\frac{4(5+b^2)}{(5+b^2)^2+2}\left(1+\frac{2}{(5+d^2)^2}\right)\right)f_bf_d\\
&= \left(-\frac{4}{(5+d^2)^2}+\frac{4(5+b^2)}{(5+b^2)^2+2}\left(1+\frac{2}{(5+d^2)^2}\right)\right)f_bf_d\\
g''_{33}(\overline{e}''_3) &= g_{11}(\overline{e}''_3)A_3^1A_3^1 + 2g_{12}(\overline{e}''_3)A_3^1A_3^2 + 2g_{13}(\overline{e}''_3)A_3^1A_3^3\\
&+g_{22}(\overline{e}''_3)A_3^2A_3^2 + 2g_{23}(\overline{e}''_3)A_3^2A_3^3 + g_{33}(\overline{e}''_3)A_3^3A_3^3\\
&= \left(1-\frac{88+46d^2+6d^4}{(5+d^2)^3}\right)(-f_d)(-f_d)\\
&+2\frac{-23+4d^2+3d^4}{2(5+d^2)^3}(-f_d)2f_d\\
&+2\frac{-22d-6d^3}{(5+d^2)^3}(-f_d)df_d + \left(1+\frac{-2+6d^2}{(5+d^2)^3}\right)2f_d2f_d\\
&+2\frac{-11d+d^3}{(5+d^2)^3}2f_ddf_d + \left(1+\frac{-6d^2+10}{(5+d^2)^3}\right)df_ddf_d\\
&= \left(5+d^2-\frac{50+20d^2+2d^4}{(5+d^2)^3}\right)f_d^2 = \left(5+d^2-\frac{2(5+d^2)^2}{(5+d^2)^3}\right)f_d^2\\
&= \left(5+d^2-\frac{2}{5+d^2}\right)f_d^2 = \frac{1}{f_d^2}f_d^2 = 1
\end{aligned}\right.$$

The equation (3.4.62) follows from (3.4.67), (3.4.68). □

Corollary 3.4.7. *Relative to the basis $\overline{\overline{e}}''$ defined by the equation (3.4.56), the matrix (2.3.9) corresponding to the Finslerian metric $F^2$ defined by equation (3.4.1) has form*

$$\begin{pmatrix} 1 & 0 & 0 \\ -\frac{(b^2+3)(b^2+4)}{(5+b^2)^2} f_b & 1 & 0 \\ -\frac{(d^2+3)(d^2+4)}{(5+d^2)^2} f_d & \boxed{\begin{matrix} f_b f_d\left(-\frac{4}{(5+d^2)^2}\right. \\ +\frac{4(5+b^2)}{(5+b^2)^2+2} \\ \left.*\left(1+\frac{2}{(5+d^2)^2}\right)\right) \end{matrix}} & 1 \end{pmatrix} \tag{3.4.69}$$

□

### 3.4.3. Example of motion in Minkowski Space of Dimension 3.

Theorem 3.4.8. *Consider Minkowski space with the Finslerian metric defined by equation (3.4.1). According to construction in subsection 3.4.1, the basis $\overline{\overline{e}}_1$ defined by the equation (3.4.28),*

$$\begin{pmatrix} \overline{e}_{1\cdot 1} & \overline{e}_{1\cdot 2} & \overline{e}_{1\cdot 3} \end{pmatrix} = \begin{pmatrix} \overline{e}_1 & \overline{e}_2 & \overline{e}_3 \end{pmatrix} \begin{pmatrix} 1 & 0 & -\sqrt{\frac{5}{23}} \\ 0 & 0 & 2\sqrt{\frac{5}{23}} \\ 0 & 1 & 0 \end{pmatrix} \tag{3.4.70}$$

$$\overline{e}_{1\cdot k} = e_{1\cdot}{}^{i}_{k}\, \overline{e}_i \tag{3.4.71}$$

*is orthonormal basis.*

*According to construction in subsection 3.4.2, the basis $\overline{\overline{e}}_2$ defined by the equation (3.4.28),*

$$\begin{pmatrix} \overline{e}_{2\cdot 1} & \overline{e}_{2\cdot 2} & \overline{e}_{2\cdot 3} \end{pmatrix} = \begin{pmatrix} \overline{e}_1 & \overline{e}_2 & \overline{e}_3 \end{pmatrix} \begin{pmatrix} 1 & -f_b & -f_d \\ 0 & 2f_b & 2f_d \\ 0 & bf_b & df_d \end{pmatrix} \tag{3.4.72}$$

$$\overline{e}_{2\cdot k} = e_{2\cdot}{}^{i}_{k}\, \overline{e}_i \tag{3.4.73}$$

*is orthonormal basis.*

*Consider the motion $\overline{A}$ mapping the basis $\overline{\overline{e}}_1$ into the basis $\overline{\overline{e}}_2$*

$$e_{2\cdot}{}^{i}_{k} = A^i_j e_{1\cdot}{}^{j}_{k} \tag{3.4.74}$$

*Let the motion $\overline{A}$ maps the basis $\overline{\overline{e}}_2$ into the basis $\overline{\overline{e}}_3$*

$$e_{3\cdot}{}^{i}_{k} = A^i_j e_{2\cdot}{}^{j}_{k} \tag{3.4.75}$$

*The basis $\overline{\overline{e}}_3$ is not orthonormal basis.*

Proof. From the equation (3.4.74), it follows that

$$A^i_j = e_{2\cdot}{}^{i}_{k} e_1^{-1}{}_{\cdot j}^{k} \tag{3.4.76}$$

From the equation (3.4.70), it follows that

(3.4.77)
$$e_1^{-1} = \begin{pmatrix} 1 & \frac{1}{2} & 0 \\ 0 & 0 & 1 \\ 0 & \frac{1}{2}\sqrt{\frac{23}{5}} & 0 \end{pmatrix}$$

From equations (3.4.72), (3.4.76), (3.4.77), it follows that

(3.4.78)
$$A = \begin{pmatrix} 1 & -f_b & -f_d \\ 0 & 2f_b & 2f_d \\ 0 & bf_b & df_d \end{pmatrix} \begin{pmatrix} 1 & \frac{1}{2} & 0 \\ 0 & 0 & 1 \\ 0 & \frac{1}{2}\sqrt{\frac{23}{5}} & 0 \end{pmatrix} = \begin{pmatrix} 1 & \frac{1}{2} - \frac{1}{2}f_d\sqrt{\frac{23}{5}} & -f_b \\ 0 & f_d\sqrt{\frac{23}{5}} & 2f_b \\ 0 & \frac{1}{2}df_d\sqrt{\frac{23}{5}} & bf_b \end{pmatrix}$$

From equations (3.4.72), (3.4.78), (3.4.75), it follows that

(3.4.79)
$$\begin{aligned} \overline{\overline{e}}_3 &= \begin{pmatrix} 1 & \frac{1}{2} - \frac{1}{2}f_d\sqrt{\frac{23}{5}} & -f_b \\ 0 & f_d\sqrt{\frac{23}{5}} & 2f_b \\ 0 & \frac{1}{2}df_d\sqrt{\frac{23}{5}} & bf_b \end{pmatrix} \begin{pmatrix} 1 & -f_b & -f_d \\ 0 & 2f_b & 2f_d \\ 0 & bf_b & df_d \end{pmatrix} \\ &= \left( \begin{pmatrix} 1 \\ 0 \\ 0 \end{pmatrix} \quad \begin{pmatrix} -f_bf_d\sqrt{\frac{23}{5}} - bf_b^2 \\ 2f_bf_d\sqrt{\frac{23}{5}} + 2bf_b^2 \\ df_bf_d\sqrt{\frac{23}{5}} + b^2f_b^2 \end{pmatrix} \quad \begin{pmatrix} -f_d^2\sqrt{\frac{23}{5}} - df_bf_d \\ 2f_d^2\sqrt{\frac{23}{5}} + 2df_bf_d \\ df_d^2\sqrt{\frac{23}{5}} + bdf_bf_d \end{pmatrix} \right) \end{aligned}$$

Let

(3.4.80)
$$\begin{array}{ll} \alpha = f_bf_d\sqrt{\frac{23}{5}} + bf_b^2 & \gamma = df_bf_d\sqrt{\frac{23}{5}} + b^2f_b^2 \\ \beta = f_d^2\sqrt{\frac{23}{5}} + df_bf_d & \rho = df_d^2\sqrt{\frac{23}{5}} + bdf_bf_d \end{array}$$

Then

(3.4.81)
$$\overline{e}_{3\cdot 2} = \begin{pmatrix} -\alpha \\ 2\alpha \\ \gamma \end{pmatrix} \quad \overline{e}_{3.3} = \begin{pmatrix} -\beta \\ 2\beta \\ \rho \end{pmatrix}$$

Since the subsequent calculations are very complicated, then in order to get an answer for the specific value of $a$, we designed application on $C\#$. Calculations are based on the theorem 3.4.5 and equations (3.4.48), (3.4.54), (3.4.55), (3.4.80). The theorem follows from calculations considered in examples 3.4.9, 3.4.10. □

EXAMPLE 3.4.9. Let $b=1$. Then

$$d=-4.3684210526 \quad f_b=0.4200840252 \quad f_d=0.2041239040$$

$$\overline{e}_{3\cdot 2}=\begin{pmatrix}-0.3603821145\\0.7207642289\\-0.6269323948\end{pmatrix} \quad \overline{e}_{3\cdot 3}=\begin{pmatrix}0.2852237394\\-0.5704474788\\-0.7649717899\end{pmatrix}$$

$$g(\overline{e}_{3\cdot 2})=\begin{pmatrix}0.4543053517&0.0160337367&0.1351074203\\0.0160337367&1.031248994&0.0268266657\\0.1351074203&0.0268266657&0.9842227918\end{pmatrix}$$

$$g_{ij}(\overline{e}_{3\cdot 2})e_{3\cdot}{}^{i}_{2}e_{3\cdot}{}^{j}_{3}=0.0157973688$$

□

EXAMPLE 3.4.10. Let $b=100$. Then

$$d=-0.049996002 \quad f_b=0.009997501 \quad f_d=0.4661156563$$

$$\overline{e}_{3\cdot 2}=\begin{pmatrix}-0.019989571554\\0.0399791431\\0.9990005814\end{pmatrix} \quad \overline{e}_{3\cdot 3}=\begin{pmatrix}-0.4657459676\\0.9314919353\\-0.0465951802\end{pmatrix}$$

$$g(\overline{e}_{3.2})=\begin{pmatrix}0.9976047372&0.0005972998&-0.0000785154\\0.0005972998&1.000000956&0.0000079285\\-0.0000785154&0.0000079285&0.9999990446\end{pmatrix}$$

$$g_{ij}(\overline{e}_{3\cdot 2})e_{3\cdot}{}^{i}_{2}e_{3\cdot}{}^{j}_{3}=-0.000013181$$

□

## 3.5. Example of Minkowski Space of Dimension 4

We consider example of the Finslerian metric of Minkowski space of dimension 4 (considered metric is a special case of metric [12]-(3); we assume $c=1$)

(3.5.1)
$$\begin{aligned}F^2(\overline{v})&=-(v^0)^2+((v^1)^2+(v^2)^2+(v^3)^2)\left(1+\frac{(v^1)^3v^2}{((v^1)^2+(v^2)^2+(v^3)^2)^2}\right)\\&=-(v^0)^2+(v^1)^2+(v^2)^2+(v^3)^2+\frac{(v^1)^3v^2}{(v^1)^2+(v^2)^2+(v^3)^2}\end{aligned}$$

We can present the corresponding Minkowski space $\overline{V}$ as direct sum

$$\overline{V}=\overline{V}_3\oplus\overline{R}$$

where $\overline{V}_3$ is Minkowski Space considered in section 3.4.

Let

$$F_r^2(v^0)=-(v^0)^2$$

(3.5.2)
$$F_3^2(v^1,v^2,v^3)=(v^1)^2+(v^2)^2+(v^3)^2+\frac{(v^1)^3v^2}{(v^1)^2+(v^2)^2+(v^3)^2}$$

Accordingly, we can present a function $F^2$ as sum

$$F^2(v^0,v^1,v^2,v^3)=F_3^2(v^1,v^2,v^3)+F_r^2(v^0) \tag{3.5.3}$$

Function $F_3^2$ is the Finslerian metric of Minkowski space considered in the section 3.4.

THEOREM 3.5.1. *The function $F_3^2$ is continuous at $(0,0,0)$ and can be extended by continuity*[3.9]

$$F_3^2(0,0,0)=0 \tag{3.5.4}$$

PROOF. Since

$$(v^1)^2\le(v^1)^2+(v^2)^2+(v^3)^2$$

then it follows that $((v^1,v^2,v^3)\ne(0,0,0))$

$$\frac{1}{(v^1)^2+(v^2)^2+(v^3)^2}\le\frac{1}{(v^1)^2} \tag{3.5.5}$$

From the inequality (3.5.5) and from the equation (3.5.2), it follows that $((v^1,v^2,v^3)\ne(0,0,0))$

$$|F_3^2(v^1,v^2,v^3)|\le(v^1)^2+(v^2)^2+(v^3)^2+\left|\frac{(v^1)^3v^2}{(v^1)^2}\right|=(v^1)^2+(v^2)^2+(v^3)^2+|v^1v^2| \tag{3.5.6}$$

For given $\epsilon>0$, let $\delta=\frac{\sqrt{\epsilon}}{2}$. Then from the inequality (3.5.6) and requirement

$$|v^1|<\delta\quad|v^2|<\delta\quad|v^3|<\delta$$

it follows that $((v^1,v^2,v^3)\ne(0,0,0))$

$$|F_3^2(v^1,v^2,v^3)|\le\frac{\epsilon}{4}+\frac{\epsilon}{4}+\frac{\epsilon}{4}+\frac{\epsilon}{4}=\epsilon \tag{3.5.7}$$

Therefore, the function $F_3^2$ is continuous at $(0,0,0)$ and

$$\lim_{(v^1,v^2,v^3)\to(0,0,0)}F_3^2(v^1,v^2,v^3)=0$$

□

THEOREM 3.5.2.

$$F^2(a,0,0,0)=F_r^2(a)=-a^2 \tag{3.5.8}$$

PROOF. The equation (3.5.8) follows from the equation (3.5.3) and the theorem 3.5.1. □

[3.9] The value of function $F_3^2$ at $(0,0,0)$ is not good defined, since the expression involves the ratio of two infinitesimal.

Theorem 3.5.3. *The Finslerian metric $F^2$ defined by equation (3.5.1) generates metric tensor of Minkowski space*

$$g_{00}(\overline{v}) = -1 \tag{3.5.9}$$

$$g_{11}(\overline{v}) = 1 - \frac{(v^1)^3(v^2)^3 + (v^1)^3v^2(v^3)^2}{((v^1)^2+(v^2)^2+(v^3)^2)^3} + \frac{6v^1(v^2)^3(v^3)^2 + 3v^1(v^2)^5 + 3v^1v^2(v^3)^4}{((v^1)^2+(v^2)^2+(v^3)^2)^3} \tag{3.5.10}$$

$$g_{12}(\overline{v}) = \frac{1}{2}\frac{(v^1)^6 + 6(v^1)^4(v^2)^2 - 3(v^1)^2(v^2)^4 + 4(v^1)^4(v^3)^2 + 3(v^1)^2(v^3)^4}{((v^1)^2+(v^2)^2+(v^3)^2)^3} \tag{3.5.11}$$

$$g_{13}(\overline{v}) = \frac{(v^1)^4v^2v^3 - 3(v^1)^2(v^2)^3v^3 - 3(v^1)^2v^2(v^3)^3}{((v^1)^2+(v^2)^2+(v^3)^2)^3} \tag{3.5.12}$$

$$g_{22}(\overline{v}) = 1 + \frac{(v^1)^3(v^2)^3 - 3(v^1)^5v^2 - 3(v^1)^3v^2(v^3)^2}{((v^1)^2+(v^2)^2+(v^3)^2)^3} \tag{3.5.13}$$

$$g_{23}(\overline{v}) = \frac{3(v^1)^3(v^2)^2v^3 - (v^1)^5v^3 - (v^1)^3(v^3)^3}{((v^1)^2+(v^2)^2+(v^3)^2)^3} \tag{3.5.14}$$

$$g_{33}(\overline{v}) = 1 + \frac{3(v^1)^3v^2(v^3)^2 - (v^1)^5v^2 - (v^1)^3(v^2)^3}{((v^1)^2+(v^2)^2+(v^3)^2)^3} \tag{3.5.15}$$

Proof. Equations (3.5.10), (3.5.11), (3.5.12), (3.5.13), (3.5.14), (3.5.15) follow from the theorem 3.4.1 and the equation (3.5.3). From the equation (3.5.3), it follows that

$$\frac{\partial F^2(\overline{v})}{\partial v^i} = 2F_3(v^1,v^2,v^3)\frac{\partial F_3(v^1,v^2,v^3)}{\partial v^i} \quad i = 1,2,3 \tag{3.5.16}$$

From equations (2.2.2), (3.5.1), (3.5.16), it follows that

$$g_{0j}(\overline{v}) = 0 \quad i = 1,2,3 \tag{3.5.17}$$

□

Theorem 3.5.4. *Since vectors*

$$\overline{e}_1 = \begin{pmatrix}0\\e_1^1\\e_1^2\\e_1^3\end{pmatrix} \quad \overline{e}_2 = \begin{pmatrix}0\\e_2^1\\e_2^2\\e_2^3\end{pmatrix} \quad \overline{e}_3 = \begin{pmatrix}0\\e_3^1\\e_3^2\\e_3^3\end{pmatrix} \tag{3.5.18}$$

*generate orthonormal basis of Minkowski space $\overline{V}_3$, then vectors*

$$\overline{e}_1 = \begin{pmatrix}0\\e_1^1\\e_1^2\\e_1^3\end{pmatrix} \quad \overline{e}_2 = \begin{pmatrix}0\\e_2^1\\e_2^2\\e_2^3\end{pmatrix} \quad \overline{e}_3 = \begin{pmatrix}0\\e_3^1\\e_3^2\\e_3^3\end{pmatrix} \quad \overline{e}_0 = \begin{pmatrix}1\\0\\0\\0\end{pmatrix} \tag{3.5.19}$$

*generate orthonormal basis of Minkowski space $\overline{V}$.*

Proof. According to the equation (3.5.3), since vectors $\overline{a}$, $\overline{b} \in \overline{V}_3$ are orthogonal in Minkowski space $\overline{V}_3$, then these vectors are also orthogonal in Minkowski space $\overline{V}$. Therefore, vectors $\overline{e}_1$, $\overline{e}_2$, $\overline{e}_3$ belong to orthonormal basis of Minkowski space $\overline{V}$. From the theorem 3.5.3, it follows that

$$g_{kl}(\overline{e}_i)e_i^k e_0^l = g_{k0}(\overline{e}_i)e_i^k = 0 \quad i = 1, 2, 3$$

Therefore, vector $\overline{e}_0$ is orthogonal to vector $\overline{e}_i$, $i = 1, 2, 3$. According to the equation (3.5.8),

$$F^2(\overline{e}_0) = -1$$

□

CHAPTER 4

CHAPTER 5

# Index

# Ортонормированный базис в пространстве Минковского

Александр Клейн

Александр Ложье

*E-mail address*: Aleks_Kleyn@MailAPS.org

*E-mail address*: Laugier.Alexandre@orange.fr

Аннотация. Финслерово пространство - это дифференцируемое многообразие, для которого пространство Минковского является слоем касательного расслоения. Для того, чтобы понять строение системы отсчёта в Финслеровом пространстве, мы должны понять структуру ортонормированного базиса в пространстве Минковского.

В статье рассмотрено определение ортонормированного базиса в пространстве Минковского, структура метрического тензора относительно ортонормированного базиса, процедура ортогонализации. Линейное преобразование пространства Минковского отображающее по крайней мере один ортонормированный базис в ортонормированный базис называется движением. Множество движений пространства Минковского $\overline{V}$ порождает не полную группу $SO(\overline{V})$ действующую однотранзитивно на многообразии базисов.

Пассивное преобразование пространства Минковского отображающее по крайней мере один ортонормированный базис в ортонормированный базис называется квазидвижением пространства Минковского. Множество пассивных преобразований пространства Минковского порождает пассивное представление не полной группы $SO(\overline{V})$ на многообразии базисов. Так как парные представления (активное и пассивное) не полной группы $SO(\overline{V})$ на многообразии базисов однотранзитивны, то мы можем рассмотреть определение геометрического объекта.

# Оглавление

Глава 1

# Предисловие

## 1.1. Структура множества движений пространства Минковского

Эта книга появилась в результате совместной работы авторов над некоторыми вопросами, возникшими в результате анализа статьи [8]. Евклидово пространство является предельным случаем пространства Минковского, когда метрический тензор не зависит от направления. Однако, хотя матрица метрического тензора относительно ортогонального базиса евклидова пространства диагональна, аналог (2.3.9) этой матрицы в пространстве Минковского является треугольной матрицей. Это утверждение является крайне неудовлетворительным.

Поскольку теорема 2.3.8 не ограничивает значение элементов нижнего треугольника матрицы (2.3.9), то мы должны допустить, что это значение произвольно. Это приводит к утверждению, что многообразие базисов пространства Минковского больше чем многообразие базисов евклидова пространства. Когда мы завершили расчёты к секции 3.1, мы поняли, что хотя рассматриваемое значение неизвестно, но не произвольно.

Более существенной оказалась возможность упростить уравнение [8]-(4.3). Уравнение (2.4.8) напоминает аналогичное уравнение для евклидова пространства. Из уравнения (2.4.2) следует, что матрица движения пространства Минковского близка к ортогональной матрице.

Поскольку множество $SO(\overline{V})$ движений пространства Минковского размерности $n$ является подмножеством группы $GL(n)$ линейных преобразований, то мы можем рассматривать произведение движений. Однако множество $SO(\overline{V})$ не является замкнутым относительно операции произведения.[1.1] Другими словами, множество $SO(\overline{V})$ является не полной группой.[1.2] Сейчас трудно сказать какие физические следствия может иметь это утверждение. Однако так как не полная группа $SO(\overline{V})$ имеет однотранзитивное представление на многообразии базисов, и парное представление также однотранзитивно, то мы можем рассмотреть определение геометрического объекта.

## 1.2. Измерение угла

В евклидовом пространстве существует два способа измерения угла.

[1.1] Мы ожидаем статью, в которой покажем, что в пространстве Минковского с нормой (3.4.1) существуют движения, произведение которых не является движением.

[1.2] Не полная группа является частным случаем не полной $\Omega$-алгебры.

Определение 1.1.1. Если операция $\omega \in \Omega$ определена не для всякого кортежа элементов множества $A$, то на множестве $A$ определена структура **не полной $\Omega$-алгебры**. □

Например, пусть на множестве $A$ определено ассоциативное произведение, имеющее обратный элемент и единицу. Таким образом, для любого $a \in A$ существует $b \in A$, для которых определено произведение $ab$. Однако существуют $a$, $b \in A$, для которых произведение неопределенно. Тогда множестве $A$ называется **не полной группой**.

Рассмотрим угол, образованный лучами $a$ и $b$, которые имеют общую конечную точку называемую вершиной угла.

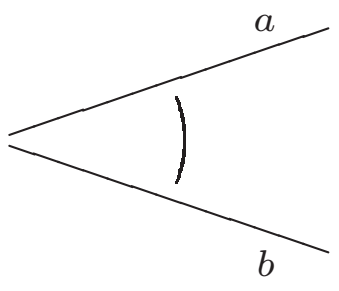

Рассмотрим единичный круг с центром в вершине угла. Мы можем определить величину угла как длину дуги, заключённой между лучами $a$ и $b$. Этот метод измерения угла используется в различных областях деятельности человека, включая астрономию (раздел [1]-1.9).

Когда мы выполняем аналитические вычисления нам удобнее определять косинус угла, опираясь на скалярное произведения векторов, параллельных лучам $a$ и $b$. Мы предполагаем, что измерения выполнены в псевдо евклидовом пространстве, где скалярное произведение симметрично. Поэтому результаты рассмотренных методов измерения совпадают.

Наблюдатель может выбрать произвольный метод измерения интервалов длины и времени. Однако для того, чтобы один наблюдатель мог сообщить результаты измерений другому наблюдателю, они должны договориться об эталонах длины и времени.

Поэтому каждый наблюдатель выбирает свой базис таким образом, что длина векторов базиса равна длине эталона. Соответствующий базис мы будем называть нормальным. Каждому базису сопутствует координатная система, построенная таким образом, что координата $x^i$ точки на оси координат $X^i$ равна расстоянию от начала координат. В евклидовом или псевдоевклидовом пространстве мы отождествляем каждую точку с радиус вектором этой точки. Длина этого вектора является расстоянием от начала координат до рассматриваемой точки. Квадрат длины этого вектора является квадратичной формой координат. В случае евклидова пространства эта квадратичная форма положительно определена.

Наблюдатель может выбрать произвольный нормальный базис. Соответственно, квадратичная форма, определяющая метрику, может иметь произвольный вид. Однако мы можем любую квадратичную форму привести к каноническому виду ([3], с. 169 - 172). Соответствующий базис является ортогональным базисом. Ортогональный нормальный базис называется ортонормированным базисом.

Таким образом, для измерения пространственных и временных интервалов наблюдатель пользуется ортонормированным базисом.[1.3]

В пространстве Минковского метрический тензор не является квадратичной формой. Поэтому процедура экспериментального определения метрики для заданного базиса является более сложной задачей. Выполняя измерение расстояния в разных направлениях, мы можем построить аналог единичной сферы в евклидовом пространстве, а именно поверхность

$$|F(\overline{v})| = 1$$

которая называется индикатрисой. В статье [10], Асанов рассмотрел, как опираясь на определённые физические гипотезы, можно определить форму финслеровой метрики.

Согласно замечанию 2.2.4, скалярное произведение в пространстве Минковского некоммутативно. Поэтому косинус угла зависит от направления измерения. Следовательно, в пространстве Минковского результаты рассмотренных выше методов измерения не совпадают. Соответственно, мы можем рассмотреть два различных определения ортонормированного базиса.

Согласно [15], вектор $\overline{w}$ ортогонален вектору $\overline{v}$, если вектор $\overline{w}$ касателен гиперповерхности

$$F^2(\overline{x}) = F^2(\overline{v}) \tag{1.2.1}$$

[1.3] Этот шаг в рассуждении кажется нелогичным. Наблюдая свободное движение твёрдого тела и отвлекаясь от сил трения, Галилей открывает закон инерции. Когда мы переходим от плоского к искривленному пространству, естественно вместо движения по прямой предположить движение по геодезической. Но также как принцип Лагранжа, требование ортонормированности не вытекает непосредственно из эксперимента. Тем более кажется удивительным, что эксперимент подтверждает эту точку зрения.

гомотетичной индикатрисе. Из равенства (2.2.9) следует, что вектор $\overline{w}$, касательный гиперповерхности (1.2.1), удовлетворяет уравнению

$$g_{ij}(\overline{v})v^i w^j = 0 \tag{1.2.2}$$

Из замечания 2.2.4 следует, что вектор $\overline{w}$ ортогонален вектору $\overline{v}$, если скалярное произведение

$$\overline{v}\cdot\overline{w} = 0$$

В этой статье мы рассматриваем определение 2.3.2 ортогонального базиса, основанное на определении (1.2.2) ортогональности.

Однако возможны и другие определения ортогональности, например [11, 13].

## 1.3. Измерение скорости света

Наблюдатель в общей теории относительности пользуется локальным базисом в касательном пространстве, которое является пространством событий специальной теории относительности. Поэтому переход наблюдателя от одного локального базиса к другому можно описать преобразованием Лоренца.

Когда геометрия пространства событий меняется настолько, что преобразование базиса меняет свою структуру, некоторые авторы воспринимают это изменение как нарушение инвариантности Лоренца ([14]). Если мы видим в эксперименте, что есть отклонение в структуре преобразования Лоренца, то это является аргументом в пользу того, что мы должны рассматривать новую геометрию. Однако если новая геометрия верна, мы не должны ограничивать себя инвариантностью Лоренца. Наша задача - понять какую алгебру порождают преобразования новой геометрии и как это влияет на принцип инвариантности.

Однако существует другая концепция нарушение инвариантности Лоренца. Эта концепция связана с зависимостью скорости света от направления ([12], с. 7, [13], с. 12). Эталоны длины и времени и скорость света неразрывно связаны в теории относительности. Если мы знаем два из этих параметров, мы можем определить третий. Выбор геометрии также важен для определения зависимости между рассматриваемыми параметрами.

Свет распространяется в финслеровом пространстве. Однако, так как операция определения скорости - локальная операция, то мы можем рассмотреть эту операцию в пространстве Минковского. Мы полагаем, что вектор $\overline{e}_0$ является времениподобным вектором, а остальные векторы базиса являются пространственноподобными векторами.

Для того, чтобы измерять скорость света, мы должны предположить, что нам заданы эталоны длины и времени. Пусть финслерова метрика пространства событий имеет вид ([12]-(3))

$$F^2(\overline{v}) = -c^2(v^0)^2 + F_3^2(v^1, v^2, v^3) \tag{1.3.1}$$

где $c$ - коэффициент, связывающий эталоны длины и времени, и $F_3^2$ - положительно определённая финслерова метрика 3-мерного пространства. Чтобы измерить скорость света в направлении в пространстве $(v^1, v^2, v^3)$ , наблюдатель откладывает в заданном направлении отрезок, длина которого равна эталону длины

$$F_3(v^1, v^2, v^3) = 1$$

и определяет время движения светового сигнала. Поскольку вектор распространения света является изотропным, то из равенства (1.3.1) следует

$$-c^2(v^0)^2 + F_3^2(v^1, v^2, v^3) = 0 \tag{1.3.2}$$

Для любого выбранного направления в пространстве $(v^1, v^2, v^3)$

$$v^0 = \frac{F_3(v^1, v^2, v^3)}{c} = \frac{1}{c} \tag{1.3.3}$$

Следовательно, время распространения света $v^0$ на заданное расстояние не зависит от направления в пространстве.

Скорость света в рассматриваемой модели не зависит от направления и равна $c$. В пространстве Финслера финслерова метрика зависит от точки многообразия. Поэтому процедура измерения скорости света должна быть выполнена в бесконечно малой области.

Попытка понять, как взаимодействие элементарных частиц и квантовых полей может изменить геометрию пространства событий, является одной из причин изучения финслеровой геометрии ([10, 12]). Без сомнения, только эксперимент может подтвердить, является ли геометрия пространства событий финслеровой и зависит ли скорость света от направления.[1.4] Наша задача - найти эффективный инструмент для экспериментальной проверки характера геометрии.

[1.4]Физики проводят эксперименты, проверяющие изотропность скорости света ([4, 5, 6]). Задача экспериментов найти предел нарушения Лоренц инвариантности.

Даже если окажется, что скорость света зависит от направления, это не противоречит утверждению, что скорость света в заданном направлении - это максимальная скорость передачи сигнала в заданном направлении.

Глава 2

# Пространство Минковского

## 2.1. Однородная функция

В этой статье мы будем рассматривать векторное пространство над полем действительных чисел $R$.

Определение 2.1.1. Пусть $\overline{V}$ - векторное пространство. Функция $f(\overline{x})$, $\overline{x} \in \overline{V}$, называется **однородная степени** $k$, если

$$f(a\overline{x}) = a^k f(\overline{x})$$

□

Теорема 2.1.2 (Теорема Эйлера). *Функция $f(\overline{x})$, однородная степени $k$, удовлетворяет дифференциальному уравнению*

$$\frac{\partial f(\overline{x})}{\partial x^i} x^i = k f(\overline{x}) \tag{2.1.1}$$

Доказательство. Продифференцируем равенство (2.1.1) по $a$

$$\frac{df(a\overline{x})}{da} = \frac{da^k}{da} f(\overline{x}) \tag{2.1.2}$$

Согласно правилу дифференцирования по частям, мы имеем

$$\frac{df(a\overline{x})}{da} = \frac{\partial f(a\overline{x})}{\partial ax^i}\frac{dax^i}{da} = \frac{\partial f(a\overline{x})}{\partial ax^i} x^i \tag{2.1.3}$$

Из равенств (2.1.2), (2.1.3) следует

$$\frac{\partial f(a\overline{x})}{\partial ax^i} x^i = ka^{k-1} f(\overline{x}) \tag{2.1.4}$$

Равенство (2.1.1) следует из равенства (2.1.4) если положить $a = 1$. □

Теорема 2.1.3. *Если $f(\overline{x})$ - функция, однородная степени $k$, то частные производные $\dfrac{\partial f(\overline{x})}{\partial x^i}$ являются функциями, однородными степени $k-1$.*

Доказательство. Рассмотрим равенство[2.1]

$$f(a\overline{x}) = a^k f(\overline{x}) \tag{2.1.5}$$

Продифференцируем равенство (2.1.5) по $x^i$

$$a\frac{\partial f(a\overline{x})}{\partial x^i} = a^k \frac{\partial f(\overline{x})}{\partial x^i} \tag{2.1.6}$$

Из равенства (2.1.6) следует

$$\frac{\partial f(a\overline{x})}{\partial x^i} = a^{k-1} \frac{\partial f(\overline{x})}{\partial x^i} \tag{2.1.7}$$

[2.1] Аналогичное доказательство смотри в [17], с. 265.

Следовательно, производные $\dfrac{\partial f(\overline{x})}{\partial x^i}$ являются однородными функциями степени $k-1$. □

Если отображение

$$f: R \to R$$

однородно степени 0, то $f(x) = \text{const}$. В общем случае, мы можем только утверждать, что

$$f(t\overline{x}) = f(\overline{x}) \quad t \in R$$

Мы будем также говорить, что отображение $f$ постоянно в направлении $\overline{x}$.

## 2.2. Финслерово пространство

Определения в этом разделе даны по аналогии с определениями в [12].

Определение 2.2.1. Векторное пространство $V$ называется **пространством Минковского**,[2.2] если в векторном пространстве $V$ определена **финслерова структура**

$$F: \overline{V} \to R$$

такая, что

2.2.1.1: **Финслерова метрика** $F^2$ не обязательно положительно определена[2.3]

2.2.1.2: Функция $F(\overline{x})$ однородна степени 1

$$F(a\overline{x}) = aF(\overline{x}) \quad a > 0 \tag{2.2.1}$$

2.2.1.3: Пусть $\overline{\overline{e}}$ - базис векторного пространства $A$. Координаты **метрического тензора**

$$g_{ij}(\overline{v}) = \frac{1}{2}\frac{\partial^2 F^2(\overline{v})}{\partial v^i \partial v^j} \tag{2.2.2}$$

порождают невырожденную симметричную матрицу. □

Пусть $\overline{V}$, $\overline{W}$ - пространства Минковского одинаковой размерности. Пусть

$$f: \overline{V} \to \overline{W}$$

невырожденное линейное отображение. Мы примем за общее правило ([2], с. 176), что если

$$F_V^2(\overline{v}) = F_W^2(f(\overline{v}))$$

либо

$$F_V^2(\overline{v}) = -F_W^2(f(\overline{v}))$$

где $F_V$ - финслерова структура в пространстве Минковского $\overline{V}$, $F_W$ - финслерова структура в пространстве Минковского $\overline{W}$, то мы не будем соответствующие геометрии считать существенно различными и будем изучать лишь одну из них. В частности, если $F(\overline{v}) < 0$ для любого вектора $\overline{v}$, то не нарушая общности мы можем рассматривать финслерову структуру

$$F'(\overline{v}) = -F(\overline{v})$$

вместо финслеровой структуры $F$.

---

[2.2] Я рассмотрел определение пространства Минковского согласно определению в [15], с. 28 - 32, [16], с. 44. Хотя этот термин вызывает некоторые ассоциации со специальной теорией относительности, обычно из контекста ясно о какой геометрии идёт речь.

[2.3] Это требование связано с тем, что мы рассматриваем приложения в общей теории относительности.

Определение 2.2.2. Если $F(\overline{v}) = 0$, то вектор $\overline{v}$ называется **изотропным**. Если $F^2(\overline{v}) < 0$, то вектор $\overline{v}$ называется **времениподобным**. В этом случае, $F(\overline{v})$ мнимое положительное число. Если $F^2(\overline{v}) > 0$, то вектор $\overline{v}$ называется **пространственноподобным**. В этом случае, $F(\overline{v})$ действительное положительное число. □

Определение 2.2.3. Для произвольного вектора $\overline{v} \in \overline{V}$, $F(\overline{v}) \neq 0$, вектор

$$\overline{u} = \frac{1}{|F(\overline{v})|}\overline{v}$$

называется **единичным вектором** □

Замечание 2.2.4. Согласно определению (2.2.2), метрический тензор зависит от вектора $\overline{v}$. Поэтому скалярное произведение, определённое равенством

$$\overline{v} \cdot \overline{w} = g_{ij}(\overline{v})v^i w^j$$

вообще говоря, некоммутативно. □

Пусть $\overline{\overline{e}}$ - базис пространства Минковского. Поскольку вектор $\overline{v}$ имеет разложение

$$\overline{v} = v^i \overline{e}_i$$

относительно базиса $\overline{\overline{e}}$, то отображение $F$ можно также представить в виде

$$F(v^1, ..., v^n) = F(\overline{v})$$

Теорема 2.2.5. *Финслерова структура пространства Минковского удовлетворяет дифференциальным уравнениям*

$$\frac{\partial F(\overline{a})}{\partial a^i}a^i = F(\overline{a}) \tag{2.2.3}$$

$$\frac{\partial^2 F(\overline{a})}{\partial a^i \partial a^j}a^i = 0 \tag{2.2.4}$$

$$\frac{1}{2}\frac{\partial^2 F^2(\overline{a})}{\partial a^i \partial a^j}a^i a^j = F^2(\overline{a}) \tag{2.2.5}$$

Доказательство. Равенство (2.2.3) следует из утверждения (2.2.1.2) определения 2.2.1 и теоремы 2.1.2. Согласно теореме 2.1.2 производная $\dfrac{\partial F(\overline{x})}{\partial x^i}$ является однородной функцией степени 0, откуда следует равенство (2.2.4).

Последовательно дифференцируя функцию $F^2$, мы получим[2.4]

$$\frac{\partial F^2(\overline{x})}{\partial x^i} = 2F(\overline{x})\frac{\partial F(\overline{x})}{\partial x^i}$$

$$\frac{1}{2}\frac{\partial F^2(\overline{x})}{\partial x^j \partial x^i} = \frac{\partial F(\overline{x})}{\partial x^j}\frac{\partial F(\overline{x})}{\partial x^i} + F(\overline{x})\frac{\partial^2 F(\overline{x})}{\partial x^j \partial x^i} \tag{2.2.6}$$

Из равенств (2.2.3), (2.2.4), (2.2.6) следует

$$\frac{1}{2}\frac{\partial F^2(\overline{x})}{\partial x^j \partial x^i}x^i = \frac{\partial F(\overline{x})}{\partial x^j}\frac{\partial F(\overline{x})}{\partial x^i}x^i + F(\overline{x})\frac{\partial^2 F(\overline{x})}{\partial x^j \partial x^i}x^i = \frac{\partial F(\overline{x})}{\partial x^j}F(\overline{x}) \tag{2.2.7}$$

Из равенства (2.2.7) следует

$$\frac{1}{2}\frac{\partial F^2(\overline{x})}{\partial x^j \partial x^i}x^i x^j = \frac{\partial F(\overline{x})}{\partial x^j}x^j F(\overline{x}) \tag{2.2.8}$$

Равенство (2.2.5) следует из равенств (2.2.3), (2.2.8). □

[2.4] см. также [15], с. 24.

Теорема 2.2.6.

$$g_{ij}(\overline{v})v^iv^j = F^2(\overline{v}) \tag{2.2.9}$$

Доказательство. Равенство (2.2.9) является следствием равенств (2.2.2), (2.2.5). □

Теорема 2.2.7. *Метрический тензор* $g_{ij}(\overline{a})$ *является однородной функцией степени* 0 *и удовлетворяет уравнению*

$$\frac{\partial g_{ij}(\overline{a})}{\partial a^k}a^k = 0 \tag{2.2.10}$$

Доказательство. Из утверждения (2.2.1.2) определения 2.2.1 следует, что отображение $F^2(\overline{v})$ однородно степени 2. Из теоремы 2.1.3 следует, что функция

$$\frac{\partial F^2(\overline{x})}{\partial x^i}$$

однородна степени 1. Из теоремы 2.1.3 и определения (2.2.2) следует, что функция $g_{ij}(\overline{x})$ однородна степени 0. Равенство (2.2.10) следует из теоремы 2.1.2. □

Замечание 2.2.8. Как отметил Рунд в [15] на с. 35, тензор

$$C_{ijk}(\overline{a}) = \frac{\partial g_{ij}(\overline{a})}{\partial a^k} = \frac{1}{2}\frac{\partial^3 F^2(\overline{a})}{\partial a^k \partial a^i \partial a^j}$$

симметричен по всем индексам. Его компоненты являются однородными функциями степени $-1$ и удовлетворяют следующим уравнениям

$$\begin{aligned} C_{ijk}(\overline{a})a^k &= C_{kij}(\overline{a})a^k = C_{ikj}(\overline{a})a^k = 0 \\ \frac{\partial C_{kij}(\overline{a})}{\partial a^h}a^k &= \frac{\partial C_{ikj}(\overline{a})}{\partial a^h}a^k = \frac{\partial C_{ijk}(\overline{a})}{\partial a^h}a^k = 0 \end{aligned}$$

□

Теорема 2.2.9. *Пусть координаты тензора* $g$ *определены относительно базиса* $\overline{\overline{e}}$. *Тогда*

$$\left.\frac{\partial g_{ij}(\overline{a})}{\partial a^l}\right|_{\overline{a}=\overline{e}_l} = 0 \tag{2.2.11}$$

Доказательство. Равенство (2.2.11) является следствием равенства (2.2.10), так как $e_i^m = \delta_i^m$. □

Рассмотрим бесконечно малое преобразование

$$\overline{a}' = \overline{a} + d\overline{a}$$

Тогда координаты метрического тензора испытывает бесконечно малое преобразование

$$g_{ij}(\overline{a}') = g_{ij}(\overline{a}) + \frac{\partial g_{ij}(\overline{a})}{\partial a^k}da^k \tag{2.2.12}$$

Определение 2.2.10. Многообразие $M$ называется **финслеровым пространством**, если его касательное пространство являются пространством Минковского и финслерова структура $F(x,\overline{v})$ непрерывно зависит от точки касания $x \in M$. □

Замечание 2.2.11. Следствием того, что финслерова структура в касательном пространстве непрерывно зависит от точки касания, является возможность определения дифференциала длины кривой на многообразии

$$dl = |F(x, d\overline{x})|$$

Обычно сперва определяют финслерово пространство, а потом рассматривают касательное к нему пространство Минковского. На самом деле порядок определений несущественен. В этой статье, основным объектом исследования является пространство Минковского. □

## 2.3. Ортогональность

Как отметил Рунд в [15], с. 47, существуют различные определения тригонометрических функций в пространстве Минковского. Нас интересует прежде всего понятие ортогональности.

Определение 2.3.1. Вектор $v_2$ **ортогонален** вектору $v_1$, если

$$g_{ij}(v_1)v_1^i v_2^j = 0$$

□

Как мы видим из определения 2.3.1, отношение ортогональности некоммутативно. Это утверждение имеет следующее следствие. Пусть $n$ - размерность пространства Минковского $\overline{V}$. Множество векторов $\overline{v}$, ортогональных вектору $\overline{a}$, удовлетворяет линейному уравнению

$$g_{ij}(\overline{a})a^i v^j = 0$$

и, следовательно, является векторным пространством размерности $n-1$. Множество векторов $\overline{v}$, которым вектор $\overline{a}$ ортогонален, удовлетворяет уравнению

$$g_{ij}(\overline{v})v^i a^j = 0$$

Это множество, вообще говоря, не является векторным пространством.

Определение 2.3.2. Множество векторов $\overline{e}_1$, ..., $\overline{e}_p$ называется ортогональным, если

$$\begin{aligned} &g_{ij}(\overline{e}_k)e_k^i e_k^j \neq 0 \\ &g_{ij}(\overline{e}_k)e_k^i e_l^j = 0 \quad k < l \end{aligned} \tag{2.3.1}$$

Базис $\overline{\overline{e}}$ называется **ортогональным**, если его векторы формируют ортогональное множество. □

Определение 2.3.3. Базис $\overline{\overline{e}}$ называется **ортонормированным**, если это ортогональный базис и его векторы имеют единичную длину. □

Поскольку отношение ортогональности некоммутативно, то важен порядок векторов при определении ортогонального базиса. Существуют различные процедуры ортогонализации в пространстве Минковского, например, [18], с. 39. Ниже мы рассмотрим процедуру ортогонализации, предложенную в [3], с. 213 - 214.

Теорема 2.3.4. *Пусть $\overline{e}_1$, ..., $\overline{e}_p$ - ортогональное множество не изотропных векторов. Тогда векторы $\overline{e}_1$, ..., $\overline{e}_p$ линейно независимы.*

Доказательство. Рассмотрим равенство

$$a_1\overline{e}_1 + ... + a_p\overline{e}_p = 0 \tag{2.3.2}$$

Из равенства (2.3.2) следует

$$a_1 g_{ij}(\overline{e}_1)e_1^i e_1^j + ... + a_p g_{ij}(\overline{e}_1)e_1^i e_p^j = 0 \tag{2.3.3}$$

Поскольку вектор $\overline{e}_1$ не изотропен, то

$$g_{ij}(\overline{e}_1)e_1^i e_1^j \neq 0 \tag{2.3.4}$$

Из условий (2.3.1), (2.3.4) и равенства (2.3.3) следует $a_1 = 0$.

Если мы доказали, что $a^1 = ... = a^{m-1} = 0$, то равенство (2.3.2) примет вид

$$a_m\overline{e}_m + ... + a_p\overline{e}_p = 0 \tag{2.3.5}$$

Из равенства (2.3.5) следует

$$a_m g_{ij}(\overline{e}_m)e_m^i e_m^j + ... + a_p g_{ij}(\overline{e}_m)e_m^i e_p^j = 0 \tag{2.3.6}$$

Поскольку вектор $\overline{e}_m$ не изотропен, то

(2.3.7)
$$g_{ij}(\overline{e}_m)e^i_m e^j_m \neq 0$$

Из условий (2.3.1), (2.3.7) и равенства (2.3.6) следует $a_m = 0$. □

Теорема 2.3.5. *В пространстве Минковского с положительно определённой финслеровой метрикой существует ортонормированный базис.*

Доказательство. Пусть $n$ - размерность пространства Минковского. Пусть $\overline{\overline{e}}'$ - базис в пространстве Минковского.

Мы положим
$$\overline{e}_1 = \overline{e}'_1$$

Допустим мы построили множество векторов $\overline{e}_1$, ..., $\overline{e}_m$. Дополнительно предположим, что для всякого $i$, $1 \le i \le m$, вектор $\overline{e}_i$ является линейной комбинацией векторов $\overline{e}'_1$, ..., $\overline{e}'_m$. Это предположение будет выполнено и для вектора $\overline{e}_{m+1}$, если мы этот вектор представим в виде
$$\overline{e}_{m+1} = a_1\overline{e}_1 + ... + a_m\overline{e}_m + \overline{e}'_{m+1}$$

$\overline{e}_{m+1} \neq \overline{0}$, так как $\overline{\overline{e}}'$ - базис, и вектор $\overline{e}'_{m+1}$ не входит в разложение векторов $\overline{e}_1$, ..., $\overline{e}_m$. Для выбора вектора $\overline{e}_{m+1}$ мы потребуем, чтобы вектор $\overline{e}_{m+1}$ был ортогонален векторам $\overline{e}_1$, ..., $\overline{e}_m$.

(2.3.8)
$$\begin{gathered} g_{ij}(\overline{e}_1)e^i_1 e^j_{m+1} = 0 \\ ... \\ g_{ij}(\overline{e}_m)e^i_m e^j_{m+1} = 0 \end{gathered}$$

Система линейных уравнений (2.3.8) имеет вид
$$\begin{aligned} & a_1 g_{ij}(\overline{e}_1)e^i_1 e^j_1 && = -g_{ij}(\overline{e}_1)e^i_1 e'^j_{m+1} \\ & a_1 g_{ij}(\overline{e}_2)e^i_2 e^j_1 + a_2 g_{ij}(\overline{e}_2)e^i_2 e^j_2 && = -g_{ij}(\overline{e}_2)e^i_2 e'^j_{m+1} \\ & ... \\ & a_1 g_{ij}(\overline{e}_m)e^i_m e^j_1 + a_2 g_{ij}(\overline{e}_m)e^i_m e^j_2 + ... + a_m g_{ij}(\overline{e}_m)e^i_m e^j_m && = -g_{ij}(\overline{e}_m)e^i_m e'^j_{m+1} \end{aligned}$$

Следовательно, решение системы линейных уравнений (2.3.8) имеет вид
$$\begin{aligned} a_1 &= -\frac{g_{ij}(\overline{e}_1)e^i_1 e'^j_{m+1}}{g_{ij}(\overline{e}_1)e^i_1 e^j_1} \\ a_2 &= -\frac{g_{ij}(\overline{e}_2)e^i_2 e'^j_{m+1} + a_1 g_{ij}(\overline{e}_2)e^i_2 e^j_1}{g_{ij}(\overline{e}_2)e^i_2 e^j_2} \\ & ... \\ a_m &= -\frac{g_{ij}(\overline{e}_m)e^i_m e'^j_{m+1} + \sum_{k=1}^{m-1} a_k g_{ij}(\overline{e}_m)e^i_m e^j_k}{g_{ij}(\overline{e}_m)e^i_m e^j_m} \end{aligned}$$

Продолжая этот процесс, мы получим ортогональный базис $\overline{\overline{e}}$. Мы можем нормировать векторы базиса $\overline{\overline{e}}$ согласно правилу
$$\begin{gathered} E_k = g_{ij}(\overline{e}_k)e^i_k e^j_k \\ \overline{e}_k \rightarrow \frac{1}{\sqrt{|E_k|}}\overline{e}_k \\ k = 1, ..., n \end{gathered}$$

□

Рассмотрим пространство Минковского $\overline{V}$ размерности $n$. Пусть финслерова метрика $F^2$ пространства Минковского $\overline{V}$ не является положительно определённой. Множество пространственноподобных векторов не является векторным пространством, однако это множество содержит максимальное множество линейно независимых векторов

$$\overline{\overline{e}}_+ = (\overline{e}_{+\cdot 1}, ..., \overline{e}_{+\cdot p})$$

такое, что любая линейная комбинация этих векторов является пространственноподобным вектором. Множество $\overline{\overline{e}}'_+$ является базисом векторного пространства $\overline{V}_+$. Множество времениподобных векторов не является векторным пространством, однако это множество содержит максимальное множество линейно независимых векторов

$$\overline{\overline{e}}_- = (\overline{e}_{-\cdot 1}, ..., \overline{e}_{-\cdot q})$$

такое, что любая линейная комбинация этих векторов является времениподобным вектором. Множество $\overline{\overline{e}}'_-$ является базисом векторного пространства $\overline{V}_-$. Выбор множеств $\overline{\overline{e}}'_+$, $\overline{\overline{e}}'_-$, а также векторных пространств $\overline{V}_+$, $\overline{V}_-$ неоднозначен.

Теорема 2.3.6. *В пространстве Минковского $\overline{V}_+$ с финслеровой структурой $F$ существует ортонормированный базис $\overline{\overline{e}}_+$. Множество векторов, ортогональных векторам базиса $\overline{\overline{e}}_+$, является векторным пространством $\overline{V}_+^\perp$ размерности $n-p$.*

Доказательство. Существование ортонормированного базиса $\overline{\overline{e}}_+$ в пространстве Минковского $\overline{V}_+$ следует из теоремы 2.3.5. Вектор $\overline{v}$, ортогональный векторам базиса $\overline{\overline{e}}_+$, удовлетворяет системе линейных уравнений

$$g_{ij}(\overline{e}_{+\cdot 1})e^i_{+\cdot 1}v^j = 0$$
$$...$$
$$g_{ij}(\overline{e}_{+\cdot p})e^i_{+\cdot p}v^j = 0$$

Следовательно, множество векторов $\overline{v}$ порождает векторное пространство. □

Теорема 2.3.7. *В пространстве Минковского $\overline{V}$, где финслерова метрика $F^2$ не является положительно определённой, существует ортонормированный базис $\overline{\overline{e}}$, если выполнено одно из следующих условий*

2.3.7.1: *Финслеровая метрика $F^2$ пространства Минковского $\overline{V}$ имеет вид*

$$F^2(\overline{v}) = F^2_+(v^1_+, ..., v^p_+) - F^2_-(v^1_-, ..., v^q_-)$$

*где $F^2_+$, $F^2_-$ - положительно определённые финслеровы метрики. Множество переменных $v^1_+$, ..., $v^p_+$, $v^1_-$, ..., $v^q_-$ разбивает множество координат $v^1$, ..., $v^n$ вектора $\overline{v}$ на два непересекающихся множества.*

2.3.7.2: *Существует векторное пространство $\overline{V}_+$ такое, что векторное пространство $\overline{V}_+^\perp$ не содержит изотропных векторов*[2.5]

2.3.7.3: *Существуют векторные пространства $\overline{V}_+$, $\overline{V}_-$ такие, что верно следующее равенство*

$$\overline{V}_- = \overline{V}_+^\perp$$

[2.5] Хотя отображение $F$ не является линейным, мы формально можем определить множество $\ker F$ как множество изотропных векторов. Тогда условие 2.3.7.2 можно сформулировать как следующее условие. Существует векторное пространство $\overline{V}_+$ такое, что

$$\overline{V}_+^\perp \cap \ker F = \emptyset$$

Доказательство. Если выполнено условие 2.3.7.1, то естественно положить, что $p+q=n$; в противном случае, существует изотропная плоскость.[2.6] Полагая

$$v_-^1=...=v_-^q=0$$

мы видим, что векторное пространство $\overline{V}_+$ имеет размерность $p$. Полагая

$$v_+^1=...=v_+^p=0$$

мы видим, что векторное пространство $\overline{V}_-$ имеет размерность $q$. Согласно теореме 2.3.6, в пространстве Минковского $\overline{V}_+$ с финслеровой структурой $F$ существует ортонормированный базис $\overline{\overline{e}}_+$ и размерность векторного пространства $\overline{V}_+^\perp$ равна $n-p=q$. Следовательно, условие 2.3.7.1 является частным случаем условия 2.3.7.3.

Пространственноподобный вектор не может быть вектором векторного пространства $\overline{V}_+^\perp$. В противном случае, множество $\overline{\overline{e}}_+$ не является максимальным множеством линейно независимых пространственноподобный векторов. Если условие 2.3.7.2 верно, то векторное пространство $\overline{V}_+^\perp$ содержит только времениподобные векторы. Следовательно, условие 2.3.7.2 эквивалентно условию 2.3.7.3.

Рассмотрим условие 2.3.7.3. Пусть $\overline{\overline{e}}_+$ - ортонормированный базис пространства Минковского $\overline{V}_+$. Пусть $\overline{\overline{e}}_-$ - ортонормированный базис пространства Минковского $\overline{V}_-$. Согласно условие 2.3.7.3, каждый

---

[2.6] Например, если переменная $v^1$ не принадлежит множеству $v_+^1, ..., v_+^p, v_-^1, ..., v_-^q$ то любой вектор вида $\overline{v}=\begin{pmatrix}v^1\\0\\...\\0\end{pmatrix}$ является изотропным вектором. Это противоречит требованию 2.2.1.3 невырожденности финслеровой метрики.

Однако вырожденность финслеровой метрики необязательно следует из существования изотропной плоскости. В псевдоевклидовом пространстве с метрикой

$$F^2(\overline{v})=(v^1)^2-(v^2)^2+(v^3)^2-(v^4)^2$$

векторы

$$z_1=\begin{pmatrix}1\\1\\0\\0\end{pmatrix}\qquad z_2=\begin{pmatrix}0\\0\\1\\1\end{pmatrix}$$

являются изотропными векторами. Вектор

$$z=a^1z_1+a^2z_2$$

также является изотропным вектором. Однако векторы

$$e_{+\cdot 1}=\begin{pmatrix}1\\0\\0\\0\end{pmatrix}\qquad e_{+\cdot 2}=\begin{pmatrix}0\\0\\1\\0\end{pmatrix}$$

порождают векторное пространство $\overline{V}_+$, а векторы

$$e_{-\cdot 1}=\begin{pmatrix}0\\1\\0\\0\end{pmatrix}\qquad e_{-\cdot 2}=\begin{pmatrix}0\\0\\0\\1\end{pmatrix}$$

порождают векторное пространство $\overline{V}_-$.

вектор базиса $\overline{\overline{e}}_-$ ортогонален всем векторам базиса $\overline{\overline{e}}_+$. Согласно определению 2.3.2, базис

$$\overline{\overline{e}} = (\overline{e}_{+\cdot 1}, ..., \overline{e}_{+\cdot p}, \overline{e}_{-\cdot 1}, ..., \overline{e}_{-\cdot q})$$

является ортогональным базисом. Поскольку каждый вектор базиса $\overline{\overline{e}}$ является единичным вектором, то согласно определению 2.3.3, базис $\overline{\overline{e}}$ является ортонормированным базисом. □

Теорема 2.3.8. *Пусть $\overline{\overline{e}}$ - ортонормированный базис пространства Минковского. Если мы запишем координаты метрического тензора $g_{ij}(\overline{e}_k)$ относительно базиса $\overline{\overline{e}}$ в виде матрицы*

$$\begin{pmatrix} g_{11}(\overline{e}_1) & ... & g_{1n}(\overline{e}_1) \\ ... & ... & ... \\ g_{n1}(\overline{e}_n) & ... & g_{nn}(\overline{e}_n) \end{pmatrix} \tag{2.3.9}$$

*то матрица (2.3.9) является треугольной матрицей, диагональные элементы которой равны 1 или −1.*

Доказательство. Если мы рассмотрим координаты базиса $\overline{\overline{e}}$ относительно базиса $\overline{\overline{e}}$, то $e^i_j = \delta^i_j$. Если вектор $\overline{e}_k$ - времениподобный вектор, то, согласно определениям 2.3.2, 2.3.3,

$$\begin{aligned} g_{ij}(\overline{e}_k)\delta^i_k\delta^j_k &= -1 \quad \overline{e}_k \text{ - времениподобный вектор} \\ g_{ij}(\overline{e}_k)\delta^i_k\delta^j_k &= \ \ 1 \quad \overline{e}_k \text{ - времениподобный вектор} \\ g_{ij}(\overline{e}_k)\delta^i_k\delta^j_l &= \ \ 0 \quad k < l \end{aligned} \tag{2.3.10}$$

Из равенств (2.3.10) следует[2.7]

$$\begin{aligned} g_{kk}(\overline{e}_k) &= -1 \quad \overline{e}_k \text{ - времениподобный вектор} \\ g_{kk}(\overline{e}_k) &= \ \ 1 \quad \overline{e}_k \text{ - времениподобный вектор} \\ g_{kl}(\overline{e}_k) &= \ \ 0 \quad k < l \end{aligned}$$

Следовательно, $g_{kl}(\overline{e}_k)$ произвольно, если $k > l$. □

Если в теореме 2.3.5 ослабить требование и предположить произвольную финслерову метрику, то мы встретим серьёзную проблему. Множество векторов $\overline{e}_{m+1}$, удовлетворяющих системе линейных уравнений (2.3.8), порождает векторное пространство $\overline{V}_{n-m}$ размерности $n - m$. Если отношение ортогональности симметрично, то система линейных уравнений (2.3.8) эквивалентна системе линейных уравнений

$$\begin{aligned} g_{ij}(\overline{e}_{m+1})e^i_{m+1}e^j_1 &= 0 \\ &... \\ g_{ij}(\overline{e}_{m+1})e^i_{m+1}e^j_m &= 0 \end{aligned} \tag{2.3.11}$$

Из требования 2.2.1.3 невырожденности финслеровой метрики и равенств (2.3.11) следует, что существует вектор

$$\overline{v} \in \overline{V}_{n-m} \quad F(\overline{v}) \neq 0 \tag{2.3.12}$$

[2.7]Очевидно, что равенство (2.3.10) ничего не говорит о значении $g_{ij}(\overline{e}_k)$, если $i \neq k$, $j \neq k$, так как коэффициент при $g_{ij}(\overline{e}_k)$ в рассматриваемой сумме равен 0. Так как $g$ - симметричный тензор, то

$$g_{kl}(\overline{e}_k) = g_{lk}(\overline{e}_k) \quad k \leq l$$

При этом, если $k < l$, то $g_{kl}(\overline{e}_k) = 0$. Однако, значение $g_{lk}(\overline{e}_l)$ не определено.

Утверждение (2.3.12) гарантирует следующий шаг процедуры ортогонализации. В случае произвольной финслеровой метрики, мы не можем гарантировать, что процедура ортогонализации может быть завершена при начальном выборе некоторого базиса $\overline{\overline{e}}'$, так как очередной вектор оказался изотропным.

Задачи нахождения ортогонального базиса и приведения квадратичной формы к сумме квадратов тесно связаны в геометрии евклидова пространства. В пространстве Минковского, теорема 2.3.8 делает невозможной применимость алгоритма приведения квадратичной формы к каноническому виду.

## 2.4. Движение пространства Минковского

Структура пространства Минковского близка структуре евклидова пространства. Автоморфизм евклидова пространства называется движением. Движение евклидова пространства является линейным преобразованием, отображающим ортонормированный базис в ортонормированный.

Автоморфизм пространства Минковского так же называется движением. Движение пространства Минковского является линейным преобразованием. Если $\overline{\overline{e}}_1$, $\overline{\overline{e}}_2$ - ортонормированные базисы, то существует единственное линейное преобразование, отображающее базис $\overline{\overline{e}}_1$ в базис $\overline{\overline{e}}_2$. Однако из теоремы 3.4.8 следует, что многообразие базисов пространства Минковского не замкнуто относительно линейного преобразования. Поэтому мы должны ослабить требования к определению движения пространства Минковского.

Определение 2.4.1. Линейное преобразование пространства Минковского отображающее по крайней мере один ортонормированный базис в ортонормированный базис называется **движением**. □

Пусть движение $\overline{A}$ пространства Минковского отображает базис $\overline{\overline{e}}_1$ в базис $\overline{\overline{e}}_2$. Пусть $e_i$ - матрица координат базиса $\overline{\overline{e}}_i$. Тогда справедливо равенство

$$e_2 = Ae_1$$

Теорема 2.4.2. *Пусть движение $\overline{A}$ пространства Минковского отображает ортонормированный базис $\overline{\overline{e}}_1$ в ортонормированный базис $\overline{\overline{e}}_2$. Пусть движение $\overline{B}$ пространства Минковского отображает ортонормированный базис $\overline{\overline{e}}_2$ в ортонормированный базис $\overline{\overline{e}}_3$. Тогда произведение движений $\overline{AB}$ также является движением.*

Доказательство. Поскольку произведение линейных преобразований является линейным преобразованием, то утверждение теоремы следует из определения 2.4.1 и утверждения, что линейное преобразование $\overline{AB}$ отображает ортонормированный базис $\overline{\overline{e}}_1$ в ортонормированный базис $\overline{\overline{e}}_3$. □

Теорема 2.4.3. *Пусть линейное преобразование $\overline{A}$ является движением. Тогда линейное преобразование $\overline{A}^{-1}$ является движением.*

Доказательство. Пусть движение $\overline{A}$ пространства Минковского отображает ортонормированный базис $\overline{\overline{e}}_1$ в ортонормированный базис $\overline{\overline{e}}_2$. Поскольку линейное преобразование $\overline{A}^{-1}$ отображает ортонормированный базис $\overline{\overline{e}}_2$ в ортонормированный базис $\overline{\overline{e}}_1$, то утверждение теоремы следует из определения 2.4.1. □

Согласно определению 2.4.1, множество движений $SO(\overline{V})$ пространства Минковского $\overline{V}$ является подмножеством группы линейных преобразований $GL$. Из теорем 2.4.2, reftheorem: inverse motion следует, что на множестве движений $SO(\overline{V})$ определено ассоциативное произведение, которое имеет единицу и обратный элемент.

Пусть движение $\overline{A}$ пространства Минковского отображает ортонормированный базис $\overline{\overline{e}}_1$ в ортонормированный базис $\overline{\overline{e}}_2$. Пусть движение $\overline{B}$ пространства Минковского отображает ортонормированный базис $\overline{\overline{e}}_2$ в базис $\overline{\overline{e}}_3$. Если базис $\overline{\overline{e}}_3$ не является ортонормированным, то мы не знаем, является ли произведение $\overline{AB}$ движением. Следовательно, произведение движений не всегда определено и множество движений $SO(\overline{V})$ пространства Минковского $\overline{V}$ является не полной группой.

Теорема 2.4.4. *Пусть координаты тензора $g$ определены относительно ортонормированного базиса $\overline{\overline{e}}$. Пусть движение пространства Минковского $\overline{V}$*

$$v'^i = A^i_j v^j \tag{2.4.1}$$

*отображает базис $\overline{\overline{e}}$ в ортонормированный базис $\overline{\overline{e}}'$. Тогда*[2.8]

$$g_{ij}(\overline{e}'_k)A^i_k A^j_l = g_{kl}(\overline{e}_k) \quad k \le l \tag{2.4.2}$$

Доказательство. Согласно определениям 2.3.2, 2.3.3, векторы базиса $\overline{\overline{e}}$ удовлетворяют равенству

$$\begin{aligned} g_{ij}(\overline{e}_k)e^i_k e^j_l &= 0 \quad k < l \\ g_{ij}(\overline{e}_k)e^i_k e^j_k &= 1 \end{aligned} \tag{2.4.3}$$

и векторы базиса $\overline{\overline{e}}'$ удовлетворяют равенству

$$\begin{aligned} g_{ij}(\overline{e}'_k)e'^i_k e'^j_l &= 0 \quad k < l \\ g_{ij}(\overline{e}'_k)e'^i_k e'^j_k &= 1 \end{aligned} \tag{2.4.4}$$

Из равенств (2.4.1), (2.4.3), (2.4.4) следует

$$g_{ij}(\overline{e}_k)e^i_k e^j_l = g_{ij}(\overline{e}'_k)e'^i_k e'^j_l = g_{ij}(\overline{e}'_k)A^i_p e^p_k A^j_q e^q_l \quad k \le l \tag{2.4.5}$$

Так как $e^i_j = \delta^i_j$, то равенство (2.4.2) следует из равенства (2.4.5). □

Следствие 2.4.5. *Пусть координаты тензора $g$ определены относительно ортонормированного базиса $\overline{\overline{e}}$. Пусть движение пространства Минковского $\overline{V}$*

$$v'^i = A^i_j v^j$$

*отображает базис $\overline{\overline{e}}$ в ортонормированный базис $\overline{\overline{e}}'$. Тогда*

$$\begin{aligned} g_{ij}(\overline{e}'_k)A^i_k A^j_k &= 1 \\ g_{ij}(\overline{e}'_k)A^i_k A^j_l &= 0 \quad k < l \end{aligned} \tag{2.4.6}$$

□

Теорема 2.4.4 также следует из утверждения, что координаты движения $\overline{A}$ относительно базиса $\overline{\overline{e}}$ совпадают с координатами базиса $\overline{\overline{e}}'$

$$e'^i_j = A^i_j$$

Следствие 2.4.6. *Пусть $\overline{\overline{e}}_1$, $\overline{\overline{e}}_2$ - ортонормированные базисы. Пусть $e_{21\cdot}{}^k_l$ координаты базиса $\overline{\overline{e}}_2$ относительно базиса $\overline{\overline{e}}_1$. Тогда*

$$g_{1\cdot ij}(\overline{e}'_{2\cdot k})e_{21\cdot}{}^i_k e_{21\cdot}{}^j_l = g_{1\cdot kl}(\overline{e}_{1\cdot k}) \quad k \le l$$

*где $g_{1\cdot ij}(\overline{v})$ координаты метрического тензора относительно базиса $\overline{\overline{e}}_1$.* □

Теорема 2.4.7. *Пусть бесконечно малое движение*

$$a'^i = a^j(\delta^i_j + A^i_j dt)$$

*отображает ортонормированный базис $\overline{\overline{e}}$ в ортонормированный базис $\overline{\overline{e}}'$*

$$e'^i_k = e^j_k(\delta^i_j + A^i_j dt) \tag{2.4.7}$$

---

[2.8] Поскольку при движении базис пространства не меняется, то отображения $g_{ij}$ также не меняются.

*Тогда* $(k \le l)$

(2.4.8) $$g_{kp}(\overline{e}_k)A_l^p + g_{pl}(\overline{e}_k)A_k^p = 0$$

Доказательство. Для $k \le l$ из равенств (2.4.7), (2.2.12), следует

(2.4.9) $$\begin{aligned} g_{ij}(\overline{e}'_k)e'^i_k e'^j_l &= \left(g_{ij}(\overline{e}_k) + \frac{\partial g_{ij}(\overline{a})}{\partial a^r}\bigg|_{\overline{a}=\overline{e}_k} e_k^m A_m^r dt\right)(e_k^i + e_k^m A_m^i dt)(e_l^j + e_l^p A_p^j dt) \\ &= g_{ij}(\overline{e}_k)e_k^i e_l^j + g_{ij}(\overline{e}_k)e_k^i e_p^j A_l^p dt + g_{ij}(\overline{e}_k)e_m^i A_k^m dt e_l^j + \frac{\partial g_{ij}(\overline{a})}{\partial a^r}\bigg|_{\overline{a}=\overline{e}_k} e_k^m A_m^r dt e_k^i e_l^j \end{aligned}$$

Из равенств (2.4.3), (2.4.4), (2.4.9), следует $(k \le l)$

(2.4.10) $$0 = g_{ij}(\overline{e}_k)e_k^i e_p^j A_l^p dt + g_{ij}(\overline{e}_k)e_m^i A_k^m dt e_l^j + \frac{\partial g_{ij}(\overline{a})}{\partial a^r}\bigg|_{\overline{a}=\overline{e}_k} e_k^m A_m^r dt e_k^i e_l^j$$

Так как $e_k^i = \delta_k^i$, то из равенства (2.4.10) следует $(k \le l)$

(2.4.11) $$g_{kp}(\overline{e}_k)A_l^p + g_{pl}(\overline{e}_k)A_k^p + \frac{\partial g_{kl}(\overline{a})}{\partial a^r}\bigg|_{\overline{a}=\overline{e}_k} A_k^r dt = 0$$

Равенство (2.4.8) следует из равенств (2.2.11), (2.4.11). □

Теорема 2.4.8. *Произведение инфинитезимальных движений пространства Минковского является инфинитезимальным движением пространства Минковского.*

Доказательство. Пусть

$$\begin{aligned} f(\overline{a})^i &= a^k(\delta_k^i + A_k^i dt) \\ g(\overline{a})^i &= a^k(\delta_k^i + B_k^i dt) \end{aligned}$$

инфинитезимальные движения пространства Минковского. Преобразование $fg$ имеет вид

$$\begin{aligned} f(g(\overline{a}))^i &= (g(\overline{a}))^l(\delta_l^i + A_l^i dt) \\ &= a^p(\delta_p^l + B_p^l dt)(\delta_l^i + A_l^i dt) \\ &= a^p(\delta_p^i + A_p^i dt + B_p^i dt) \\ &= a^p(\delta_p^i + (A_p^i + B_p^i)dt) \end{aligned}$$

Из теоремы 2.4.7, следует, что координаты отображений $f$ и $g$ удовлетворяют равенствам

(2.4.12) $$\begin{aligned} g_{kp}(\overline{e}_k)A_l^p + g_{pl}(\overline{e}_k)A_k^p &= 0 \\ g_{kp}(\overline{e}_k)B_l^p + g_{pl}(\overline{e}_k)B_k^p &= 0 \end{aligned} \qquad k \le l$$

Из равенства (2.4.12) следует $(k \le l)$

$$\begin{aligned} &g_{kp}(\overline{e}_k)(A_l^p + B_l^p) + g_{pl}(\overline{e}_k)(A_k^p + B_k^p) \\ =&g_{kp}(\overline{e}_k)A_l^p + g_{pl}(\overline{e}_k)A_k^p + g_{kp}(\overline{e}_k)B_l^p + g_{pl}(\overline{e}_k)B_k^p \\ =&0 \end{aligned}$$

Следовательно, отображение $fg$ является инфинитезимальным движением пространства Минковского. □

Теорема 2.4.9. *Пусть*

$$a'^i = a^j(\delta^i_j + A^i_j dt)$$

*бесконечно малое движение. Для любого $r \in R$, преобразование*

$$a'^i = a^j(\delta^i_j + rA^i_j dt)$$

*является бесконечно малым движением.*

Доказательство. Из равенства (2.4.8) следует $(k \le l)$

$$g_{kp}(\overline{e}_k)(rA^p_l) + g_{pl}(\overline{e}_k)(rA^p_k) = 0 \tag{2.4.13}$$

Утверждение теоремы следует из равенства (2.4.13) и теоремы 2.4.7. □

Теорема 2.4.10. *Множество инфинитезимальных движений является векторным пространством над полем действительных чисел.*

Доказательство. Согласно теореме 2.4.8, множество бесконечно малых движений порождает абелевую группу $g$. Согласно теореме 2.4.9, определено представление поля действительных чисел в группе $g$. Следовательно, группа $g$ является векторным пространством. □

Так как отношение ортогональности не симметрично, то это приводит к изменению структуры метрического тензора в ортонормированном базисе. В частности, так как скалярное произведение

$$g_{ij}(\overline{e}_k)e^i_k e^j_l \quad k > l$$

не определено, то мы не можем требовать, что автоморфизм пространства Минковского сохраняет скалярное произведение.

Теорема 2.4.11. *Инфинитезимальное движение пространства Минковского*

$$a'^i = a^j(\delta^i_j + A^i_j dt) \tag{2.4.14}$$

*не сохраняет скалярное произведение. Для заданных векторов $\overline{a}$, $\overline{b}$, величину нарушения скалярного произведения можно оценить выражением*

$$g_{ij}(\overline{a}')a'^i b'^j - g_{ij}(\overline{a})a^i b^j = \left(g_{pj}(\overline{a})A^p_i + g_{ip}(\overline{a})A^p_j + \frac{\partial g_{ij}(\overline{a})}{\partial a^p}A^p_l a^l\right)a^i b^j dt \tag{2.4.15}$$

Доказательство. Пусть инфинитезимальное движение пространства Минковского отображает векторы $\overline{a}$, $\overline{b}$ в векторы $\overline{a}'$, $\overline{b}'$

$$\begin{aligned} a'^i &= a^j(\delta^i_j + A^i_j dt) \\ b'^i &= b^j(\delta^i_j + A^i_j dt) \end{aligned} \tag{2.4.16}$$

Скалярное произведение векторов $\overline{a}'$, $\overline{b}'$ имеет вид

$$\begin{aligned} g_{ij}(\overline{a}')a'^i b'^j &= g_{ij}((a^k + A^k_l a^l dt)\overline{e}_k)(a^i + A^i_p a^p dt)(b^j + A^j_q b^q dt) \\ &= (g_{ij}(\overline{a}) + \frac{\partial g_{ij}(\overline{a})}{\partial a^r}A^r_l a^l dt)(a^i + A^i_p a^p dt)(b^j + A^j_q b^q dt) \\ &= g_{ij}(\overline{a})a^i b^j + g_{ij}(\overline{a})A^i_p a^p b^j dt + g_{ij}(\overline{a})a^i A^j_p b^p dt + \frac{\partial g_{ij}(\overline{a})}{\partial a^r}A^r_l a^l dt a^i b^j \\ &= g_{ij}(\overline{a})a^i b^j + \left(g_{pj}(\overline{a})A^p_i + g_{ip}(\overline{a})A^p_j + \frac{\partial g_{ij}(\overline{a})}{\partial a^p}A^p_l a^l\right)a^i b^j dt \end{aligned} \tag{2.4.17}$$

Равенство (2.4.15) следует из равенства (2.4.17). □

## 2.5. Квазидвижение пространства Минковского

Хотя структура множества $SO(\overline{V})$ движений пространства Минковского $\overline{V}$ не является полностью определённой группой, множество $SO(\overline{V})$ является подмножеством группы $GL$ и действует однотранзитивно на многообразии базисов пространства Минковского $\overline{V}$. Поэтому в множестве пассивных преобразований векторного пространства $\overline{V}$ мы можем выделить преобразования, сохраняющие структуру пространства Минковского $\overline{V}$.

Определение 2.5.1. Пассивное преобразование пространства Минковского

$$\overline{e}'_i = A_i^j \overline{e}_j$$

отображающее по крайней мере один ортонормированный базис в ортонормированный базис называется **квазидвижением пространства Минковского**. □

Пусть квазидвижение $\overline{A}$ пространства Минковского отображает базис $\overline{\overline{e}}_1$ в базис $\overline{\overline{e}}_2$. Пусть $e_i$ - матрица координат базиса $\overline{\overline{e}}_i$. Тогда справедливо равенство

$$e_2 = e_1 A$$

Координаты вектора

$$\overline{a} = a^i \overline{e}_i$$

преобразуются согласно правилу

$$a'^j = A^{-1\cdot j}{}_i a^i \tag{2.5.1}$$

Из равенств (2.2.9), (2.5.1) следует

$$g_{ij}(\overline{a}) a^i a^j = g'_{kl}(\overline{a}) a'^k a'^l = g'_{kl}(\overline{a}) A^{-1\cdot k}{}_i a^i A^{-1\cdot l}{}_j a^j \tag{2.5.2}$$

Из равенства (2.5.2) следует

$$g'_{kl}(\overline{a}) = g_{ij}(\overline{a}) A_k^i A_l^j \tag{2.5.3}$$

Следовательно, $g_{ij}(\overline{a})$ является тензором.

Замечание 2.5.2. До сих пор мы рассматривали метрику пространства Минковского $\overline{V}$ как набор функций $g_{ij}(\overline{v})$ вектора $\overline{v} \in \overline{V}$. Поскольку базис $\overline{\overline{e}}$ задан, то мы имеем однозначное разложение вектора $\overline{v}$ относительно базиса $\overline{\overline{e}}$

$$\overline{v} = v^i \overline{e}_i$$

Следовательно, мы можем рассматривать отображение $g_{ij}(\overline{v})$ как функцию $g_{ij}(v^1, ..., v^n)$ переменных $v^1$, ..., $v^n$

$$g_{ij}(v^1, ..., v^n) = g_{ij}(\overline{v}) \tag{2.5.4}$$

□

Теорема 2.5.3. *Пусть координаты тензора $g$ определены относительно ортонормированного базиса $\overline{\overline{e}}$. Пусть квазидвижение пространства Минковского $\overline{V}$*

$$\overline{e}'_i = A_i^j \overline{e}_j \tag{2.5.5}$$

*отображает базис $\overline{\overline{e}}$ в ортонормированный базис $\overline{\overline{e}}'$. Тогда*

$$g_{ij}(\overline{e}'_k) A_k^i A_l^j = g_{kl}(\overline{e}_k) \quad k \le l \tag{2.5.6}$$

Доказательство. Согласно определениям 2.3.2, 2.3.3, векторы базиса $\overline{\overline{e}}$ удовлетворяют равенству

(2.5.7)
$$g_{ij}(\overline{e}_k)e^i_k e^j_l = 0 \quad k < l$$
$$g_{ij}(\overline{e}_k)e^i_k e^j_k = 1$$

и векторы базиса $\overline{\overline{e}}'$ удовлетворяют равенству

(2.5.8)
$$g_{ij}(\overline{e}'_k)e'^i_k e'^j_l = 0 \quad k < l$$
$$g_{ij}(\overline{e}'_k)e'^i_k e'^j_k = 1$$

Из равенств (2.5.5), (2.5.7), (2.5.8) следует

(2.5.9)
$$g_{ij}(\overline{e}_k)e^i_k e^j_l = g_{ij}(\overline{e}'_k)e'^i_k e'^j_l = g_{ij}(\overline{e}'_k)A^p_k e^i_p A^q_l e^j_q \quad k \le l$$

Так как $e^i_j = \delta^i_j$, то равенство (2.5.6) следует из равенства (2.5.9). □

Из теорем 2.4.4, 2.5.3 следует, что множество квазидвижений пространства Минковского и множество движений пространства Минковского совпадают.

Замечание 2.5.4. Множество квазидвижений порождает пассивное представление не полной группы $SO(\overline{V})$.

Пусть $\overline{\overline{e}}_1$, $\overline{\overline{e}}_2$ - ортонормированные базисы пространства Минковского $\overline{V}$. Согласно теореме 2.5.3, существует квазидвижение отображающее базис $\overline{\overline{e}}_1$ в базис $\overline{\overline{e}}_2$. Поскольку не полная группа $SO(\overline{V})$ является подгруппой группы $GL$, то это квазидвижение определено однозначно. Следовательно, представление не полной группы $SO(\overline{V})$ на многообразии базисов пространства Минковского $\overline{V}$ однотранзитивно.

Поэтому определение геометрического объекта (раздел [7]-5.3) остаётся верным также в пространстве Минковского. □

Теорема 2.5.5. *Пусть $\overline{A}$ - квазидвижение, отображающее базис $\overline{\overline{e}}$ в базис $\overline{\overline{e}}'$*

(2.5.10)
$$\overline{e}'_k = A^i_k \overline{e}_i$$

*Квазидвижение $\overline{A}$ отображает отображение $g_{ij}(v^1, ..., v^n)$ в отображение $g'_{ij}(v'^1, ..., v'^n)$ согласно правилу*

(2.5.11)
$$g'_{kl}(v'^1, ..., v'^n) = g_{ij}(v'^k A^1_k, ..., v'^k A^n_k)A^i_k A^j_l$$

Доказательство. Из равенств (2.5.3), (2.5.4), следует

(2.5.12)
$$g'_{kl}(v'^1, ..., v'^n) = g'_{kl}(\overline{v}) = g_{ij}(\overline{v})A^i_k A^j_l = g_{ij}(v^1, ..., v^n)A^i_k A^j_l$$

Поскольку квазидвижение $\overline{A}$ не изменяет вектор $\overline{v}$, то

(2.5.13)
$$\overline{v} = v^k \overline{e}_k = v'^i \overline{e}'_i = v'^k A^i_k \overline{e}_i$$

Из равенства (2.5.13) следует

(2.5.14)
$$v^i = v'^k A^i_k$$

Равенство (2.5.11) следует из равенств (2.5.12), (2.5.14). □

Теорема 2.5.6. *Пусть бесконечно малое квазидвижение*

$$\overline{e}'_k = \overline{e}_l(\delta^l_k + A^l_k dt)$$

*отображает ортонормированный базис $\overline{\overline{e}}$ в ортонормированный базис $\overline{\overline{e}}'$*

(2.5.15)
$$e'^i_k = e^i_l(\delta^l_k + A^l_k dt)$$

*Пусть координаты тензора $g$ определены относительно ортонормированного базиса $\overline{\overline{e}}$. Пусть координаты тензора $g'$ определены относительно ортонормированного базиса $\overline{\overline{e}}'$. Тогда $(k \le l)$*

$$g_{kp}(\overline{e}_k)A_l^p + g_{lp}(\overline{e}_k)A_k^p = 0 \tag{2.5.16}$$

Доказательство. Согласно (2.5.3), (2.2.12), следует, что метрический тензор испытывает бесконечно малое преобразование[2.9]

$$\begin{aligned} g'_{kl}(\overline{e}'_p) &= g_{ij}(\overline{e}'_p)(\delta^i_k\delta^j_l + \delta^i_k A^j_l dt + A^i_k\delta^j_l dt) \\ &= \left(g_{ij}(\overline{e}_p) + \left.\frac{\partial g_{ij}(\overline{a})}{\partial a^r}\right|_{\overline{a}=\overline{e}_p} \delta^r_m A^m_p dt\right)(\delta^i_k\delta^j_l + \delta^i_k A^j_l dt + A^i_k\delta^j_l dt) \\ &= g_{kl}(\overline{e}_p) + g_{kj}(\overline{e}_p)A^j_l dt + g_{il}(\overline{e}_p)A^i_k dt + \left.\frac{\partial g_{kl}(\overline{a})}{\partial a^m}\right|_{\overline{a}=\overline{e}_p} A^m_p dt \end{aligned} \tag{2.5.17}$$

Если $p = k$, то из равенства (2.5.17) следует

$$g'_{kl}(\overline{e}'_k) = g_{kl}(\overline{e}_k) + g_{kj}(\overline{e}_k)A^j_l dt + g_{il}(\overline{e}_k)A^i_k dt + \left.\frac{\partial g_{kl}(\overline{a})}{\partial a^m}\right|_{\overline{a}=\overline{e}_k} A^m_p dt \tag{2.5.18}$$

При $k \le l$, согласно теореме 2.3.8,

$$g_{kl}(\overline{e}_k) = g'_{kl}(\overline{e}'_k) \tag{2.5.19}$$

Равенство (2.5.16) следует из равенств (2.2.11), (2.5.18), (2.5.19). □

Из теорем 2.4.7, 2.5.6 следует, что множество бесконечно малых квазидвижений пространства Минковского и множество бесконечно малых движений пространства Минковского совпадают.

[2.9] Мы используем равенство $e^m_i = \delta^m_i$.

Глава 3

# Примеры пространства Минковского

## 3.1. Пример пространства Минковского размерности 2

Мы рассмотрим пример финслеровой метрики пространства Минковского размерности 2 (рассматриваемая метрика является частным случаем метрики [12]-(3))

$$
\begin{aligned}
F^2(\overline{v}) &= ((v^1)^2+(v^2)^2)\left(1+\frac{(v^1)^3v^2}{((v^1)^2+(v^2)^2)^2}\right)\\
&= (v^1)^2+(v^2)^2+\frac{(v^1)^3v^2}{(v^1)^2+(v^2)^2}
\end{aligned}
\tag{3.1.1}
$$

Теорема 3.1.1. *Финслеровая метрика, определённая равенством (3.1.1) порождает метрику пространства Минковского*

$$
g_{11}(\overline{v}) = 1+\frac{3v^1(v^2)^5-(v^1)^3(v^2)^3}{((v^1)^2+(v^2)^2)^3}
\tag{3.1.2}
$$

$$
g_{12}(\overline{v}) = \frac{1}{2}\frac{(v^1)^6+6(v^1)^4(v^2)^2-3(v^1)^2(v^2)^4}{((v^1)^2+(v^2)^2)^3}
\tag{3.1.3}
$$

$$
g_{22}(\overline{v}) = 1+\frac{(v^1)^3(v^2)^3-3(v^1)^5v^2}{((v^1)^2+(v^2)^2)^3}
\tag{3.1.4}
$$

Доказательство. Из равенства (3.1.1) следует

$$
\begin{aligned}
\frac{\partial F^2(\overline{v})}{\partial v^1} &= 2v^1+\frac{3(v^1)^2v^2((v^1)^2+(v^2)^2)-(v^1)^3v^2 2v^1}{((v^1)^2+(v^2)^2)^2}\\
&= 2v^1+\frac{3(v^1)^4v^2+3(v^1)^2(v^2)^3-2(v^1)^4v^2}{((v^1)^2+(v^2)^2)^2}\\
&= 2v^1+\frac{(v^1)^4v^2+3(v^1)^2(v^2)^3}{((v^1)^2+(v^2)^2)^2}
\end{aligned}
\tag{3.1.5}
$$

$$
\begin{aligned}
\frac{\partial F^2(\overline{v})}{\partial v^2} &= 2v^2+\frac{(v^1)^3((v^1)^2+(v^2)^2)-(v^1)^3v^2 2v^2}{((v^1)^2+(v^2)^2)^2}\\
&= 2v^2+\frac{(v^1)^5-(v^1)^3(v^2)^2}{((v^1)^2+(v^2)^2)^2}
\end{aligned}
\tag{3.1.6}
$$

Из равенств (2.2.2), (3.1.5), (3.1.6) следует

$$\begin{aligned}
g_{11}(\overline{v}) &= 1 + \frac{1}{2}\frac{(4(v^1)^3v^2 + 6v^1(v^2)^3)((v^1)^2 + (v^2)^2)^2}{((v^1)^2 + (v^2)^2)^4}\\
&- \frac{1}{2}\frac{((v^1)^4v^2 + 3(v^1)^2(v^2)^3)2((v^1)^2 + (v^2)^2)2v^1}{((v^1)^2 + (v^2)^2)^4}\\
&= 1 + \frac{2(v^1)^5v^2 + 3(v^1)^3(v^2)^3 + 2(v^1)^3(v^2)^3 + 3v^1(v^2)^5}{((v^1)^2 + (v^2)^2)^3} - \frac{2(v^1)^5v^2 + 6(v^1)^3(v^2)^3}{((v^1)^2 + (v^2)^2)^3}
\end{aligned} \tag{3.1.7}$$

$$\begin{aligned}
g_{12}(\overline{v}) &= \frac{1}{2}\frac{(5(v^1)^4 - 3(v^1)^2(v^2)^2)((v^1)^2 + (v^2)^2)^2}{((v^1)^2 + (v^2)^2)^4}\\
&- \frac{1}{2}\frac{((v^1)^5 - (v^1)^3(v^2)^2)2((v^1)^2 + (v^2)^2)2v^1}{((v^1)^2 + (v^2)^2)^4}\\
&= \frac{1}{2}\frac{5(v^1)^6 - 3(v^1)^4(v^2)^2 + 5(v^1)^4(v^2)^2 - 3(v^1)^2(v^2)^4}{((v^1)^2 + (v^2)^2)^3} - \frac{1}{2}\frac{4(v^1)^6 - 4(v^1)^4(v^2)^2}{((v^1)^2 + (v^2)^2)^3}
\end{aligned} \tag{3.1.8}$$

$$\begin{aligned}
g_{22}(\overline{v}) &= 1 + \frac{1}{2}\frac{-2(v^1)^3v^2((v^1)^2 + (v^2)^2)^2}{((v^1)^2 + (v^2)^2)^4}\\
&- \frac{1}{2}\frac{((v^1)^5 - (v^1)^3(v^2)^2)2((v^1)^2 + (v^2)^2)2v^2}{((v^1)^2 + (v^2)^2)^4}\\
&= 1 - \frac{(v^1)^5v^2 + (v^1)^3(v^2)^3}{((v^1)^2 + (v^2)^2)^3} - \frac{2(v^1)^5v^2 - 2(v^1)^3(v^2)^3}{((v^1)^2 + (v^2)^2)^3}
\end{aligned} \tag{3.1.9}$$

Равенство (3.1.2) следует из равенства (3.1.7). Равенство (3.1.3) следует из равенства (3.1.8). Равенство (3.1.4) следует из равенства (3.1.9). □

Мы начнём процедуру ортогонализации[3.1] с базиса

$$\overline{e}_1 = \begin{pmatrix}1\\0\end{pmatrix} \quad \overline{e}_2 = \begin{pmatrix}0\\1\end{pmatrix}$$

Пусть

$$\overline{e}'_1 = \overline{e}_1 \quad \overline{e}'_2 = \begin{pmatrix}e'^1_2\\e'^2_2\end{pmatrix}$$

ортогональный базис. Так как вектор $\overline{e}'_2$ ортогонален вектору $\overline{e}'_1$, то согласно определению 2.3.1

$$g_{ij}(\overline{e}'_1)e'^i_1e'^j_2 = 0$$

$$g_{1j}(\overline{e}'_1)e'^j_2 = 0 \tag{3.1.10}$$

Согласно равенствам (3.1.2), (3.1.3),

$$g_{11}(\overline{e}'_1) = 1 \quad g_{12}(\overline{e}'_1) = \frac{1}{2} \quad g_{22}(\overline{e}'_1) = 1 \tag{3.1.11}$$

Из равенств (3.1.10), (3.1.11) следует

$$e'^1_2 + \frac{1}{2}e'^2_2 = 0 \tag{3.1.12}$$

[3.1] Смотри раздел [8]-3.

Из равенства (3.1.12) следует, что мы можем положить

$$\overline{e}'_2=\begin{pmatrix}-1\\2\end{pmatrix}\tag{3.1.13}$$

Согласно равенствам (3.1.2), (3.1.3), (3.1.4), (3.1.13),

$$\left\{\begin{aligned}g_{11}(\overline{e}'_2)&=1+\frac{3(-1)(2)^5-(-1)^3(2)^3}{((-1)^2+(2)^2)^3}\\&=1+\frac{-96+8}{(1+4)^3}=\frac{37}{125}\\g_{12}(\overline{e}'_2)&=\frac{1}{2}\frac{(-1)^6+6(-1)^4(2)^2-3(-1)^2(2)^4}{((-1)^2+(2)^2)^3}\\&=\frac{1}{2}\frac{1+24-48}{(1+4)^3}=-\frac{23}{250}\\g_{22}(\overline{e}'_2)&=1+\frac{(-1)^3(2)^3-3(-1)^5 2}{((-1)^2+(2)^2)^3}\\&=1+\frac{-8+6}{(1+4)^3}=\frac{123}{125}\end{aligned}\right.\tag{3.1.14}$$

Из равенств (3.1.13), (3.1.14) следует

$$\begin{aligned}g_{ij}(\overline{e}'_2)e'^i_2e'^j_1&=g_{i1}(\overline{e}'_2)e'^i_2\\&=\frac{37}{125}(-1)+(-\frac{23}{250})2=-\frac{37+23}{125}\\&=-\frac{12}{25}\end{aligned}$$

Следовательно, вектор $\overline{e}'_1$ не ортогонален вектору $\overline{e}'_2$.

Поскольку

$$F^2(\overline{e}'_1)=1$$

то длина вектора

$$\overline{e}''_1=\overline{e}'_1=\begin{pmatrix}1\\0\end{pmatrix}$$

равна 1. Поскольку

$$F^2(\overline{e}'_2)=(-1)^2+(2)^2+\frac{(-1)^3 2}{(-1)^2+(2)^2}=1+4-\frac{2}{1+4}=\frac{23}{5}$$

то длина вектора

$$\overline{e}''_2=\sqrt{\frac{5}{23}}\overline{e}'_2=\begin{pmatrix}-\sqrt{\frac{5}{23}}\\2\sqrt{\frac{5}{23}}\end{pmatrix}$$

равна 1. Следовательно, базис $\overline{\overline{e}}''$ является ортонормированным базисом. Базисы $\overline{\overline{e}}''$ и $\overline{\overline{e}}$ связаны пассивным преобразованием

$$\begin{pmatrix}\overline{e}''_1 & \overline{e}''_2\end{pmatrix}=\begin{pmatrix}\overline{e}_1 & \overline{e}_2\end{pmatrix}\begin{pmatrix}1 & -\sqrt{\frac{5}{23}}\\ 0 & 2\sqrt{\frac{5}{23}}\end{pmatrix}\tag{3.1.15}$$

Теорема 3.1.2. *Компоненты метрического тензора относительно базиса $\overline{\overline{e}}''$ имеют вид*

$$g''(\overline{e}''_1)=\begin{pmatrix}1 & 0\\ 0 & \frac{15}{23}\end{pmatrix}\tag{3.1.16}$$

$$g''(\overline{e}''_2)=\begin{pmatrix}\frac{37}{125} & -\frac{12}{25}\sqrt{\frac{5}{23}}\\ -\frac{12}{25}\sqrt{\frac{5}{23}} & 1\end{pmatrix}\tag{3.1.17}$$

Доказательство. Из равенств (2.5.3), (3.1.11), (3.1.15) следует

$$\left\{\begin{aligned}
&g''_{11}(\overline{e}''_1)=g_{11}(\overline{e}''_1)=1\\
&g''_{12}(\overline{e}''_1)=g_{11}(\overline{e}''_1)A^1_2+g_{12}(\overline{e}''_1)A^2_2=-\sqrt{\frac{5}{23}}+\frac{1}{2}2\sqrt{\frac{5}{23}}=0\\
&g''_{22}(\overline{e}''_1)=g_{11}(\overline{e}''_1)A^1_2A^1_2+2g_{12}(\overline{e}''_1)A^1_2A^2_2+g_{22}(\overline{e}''_1)A^2_2A^2_2\\
&\qquad=\left(-\sqrt{\frac{5}{23}}\right)\left(-\sqrt{\frac{5}{23}}\right)+2\frac{1}{2}\left(-\sqrt{\frac{5}{23}}\right)2\sqrt{\frac{5}{23}}+2\sqrt{\frac{5}{23}}2\sqrt{\frac{5}{23}}=\frac{15}{23}
\end{aligned}\right.\tag{3.1.18}$$

Равенство (3.1.17) следует из (3.1.18). Из равенств (2.5.3), (3.1.14), (3.1.15) следует

$$\left\{\begin{aligned}
&g''_{11}(\overline{e}''_2)=g_{11}(\overline{e}''_2)=\frac{37}{125}\\
&g''_{12}(\overline{e}''_2)=g_{11}(\overline{e}''_2)A^1_2+g_{12}(\overline{e}''_2)A^2_2\\
&\qquad=-\frac{37}{125}\sqrt{\frac{5}{23}}+\left(-\frac{23}{250}\right)2\sqrt{\frac{5}{23}}=-\frac{12}{25}\sqrt{\frac{5}{23}}\\
&g''_{22}(\overline{e}''_2)=g_{11}(\overline{e}''_2)A^1_2A^1_2+2g_{12}(\overline{e}''_2)A^1_2A^2_2+g_{22}(\overline{e}''_2)A^2_2A^2_2\\
&\qquad=\frac{37}{125}\left(-\sqrt{\frac{5}{23}}\right)\left(-\sqrt{\frac{5}{23}}\right)+2\left(-\frac{23}{250}\right)\left(-\sqrt{\frac{5}{23}}\right)2\sqrt{\frac{5}{23}}\\
&\qquad+\frac{123}{125}2\sqrt{\frac{5}{23}}2\sqrt{\frac{5}{23}}\\
&\qquad=\left(\frac{37}{125}+\frac{46}{125}\right)\frac{5}{23}+\frac{492}{125}\frac{5}{23}=\frac{575}{125}\frac{5}{23}=1
\end{aligned}\right.\tag{3.1.19}$$

Равенство (3.1.2) следует из (3.1.19). □

Следствие 3.1.3. *Относительно базиса (3.1.15), матрица (2.3.9), соответствующая финслеровой метрике, определённой равенством (3.1.1), имеет вид*

$$\begin{pmatrix} 1 & 0 \\ -\frac{12}{25}\sqrt{\frac{5}{23}} & 1 \end{pmatrix} \tag{3.1.20}$$

□

## 3.2. Структура метрики пространства Минковского

Теорема 2.3.8 определяет значение $g_{kl}(\overline{e}_k)$, если $k \geq l$. Однако теорема 2.3.8 ничего не говорит о значении $g_{kl}(\overline{e}_k)$, если $k < l$. Впечатление, что это значение произвольно, порождает серьёзную проблему. Например, если метрика не зависит от направления, то пространство Минковского становится евклидовым пространством и матрица (2.3.9) становится диагональной матрицей; поэтому непонятно каким образом последовательность произвольных значений сходится к 0.

Из равенства (3.1.20) следует, что значение $g_{kl}(\overline{e}_k)$, $k < l$, неизвестно, но не произвольно. В разделе 3.1, мы выполнили построение для финслеровой метрики, определённой равенством (3.1.1). Наша задача - воспроизвести эти вычисления для произвольной нормы $F$ пространства Минковского размерности 2.

Пусть $F(\overline{v})$ - финслерова структура пространства Минковского размерности 2. Пусть $g(\overline{v})$ - метрика, определённая равенством (2.2.2).

Мы начнём процедуру ортогонализации[3.2] с базиса

$$\overline{e}_1 = \begin{pmatrix} 1 \\ 0 \end{pmatrix} \quad \overline{e}_2 = \begin{pmatrix} 0 \\ 1 \end{pmatrix}$$

Пусть

$$\overline{e}'_1 = \overline{e}_1 \quad \overline{e}'_2 = \begin{pmatrix} e'^1_2 \\ e'^2_2 \end{pmatrix}$$

ортогональный базис. Так как вектор $\overline{e}'_2$ ортогонален вектору $\overline{e}'_1$, то согласно определению 2.3.1

$$g_{ij}(\overline{e}'_1)e'^i_1 e'^j_2 = 0$$

$$g_{1j}(\overline{e}'_1)e'^j_2 = g_{11}(\overline{e}'_1)e'^1_2 + g_{12}(\overline{e}'_1)e'^2_2 = 0 \tag{3.2.1}$$

Из равенства (3.2.1) следует, что мы можем положить

$$\overline{e}'_2 = \begin{pmatrix} -g_{12}(\overline{e}'_1) \\ g_{11}(\overline{e}'_1) \end{pmatrix} \tag{3.2.2}$$

Из равенства (3.2.2) следует

$$\begin{aligned} g_{ij}(\overline{e}'_2)e'^i_2 e'^j_1 = g_{i1}(\overline{e}'_2)e'^i_2 &= -g_{11}(\overline{e}'_2)g_{12}(\overline{e}'_1) + g_{21}(\overline{e}'_2)g_{11}(\overline{e}'_1) \\ &= \begin{vmatrix} g_{11}(\overline{e}'_1) & g_{11}(\overline{e}'_2) \\ g_{12}(\overline{e}'_1) & g_{12}(\overline{e}'_2) \end{vmatrix} \end{aligned}$$

Если мы положим

$$F(v^1, v^2) = F(\overline{v})$$

[3.2] Смотри раздел [8]-3.

то ортонормированный базис $\overline{\overline{e}}_2''$ имеет вид

$$\begin{pmatrix}\overline{e}_1'' & \overline{e}_2''\end{pmatrix} = \begin{pmatrix}\overline{e}_1 & \overline{e}_2\end{pmatrix}\begin{pmatrix}\dfrac{1}{F(1,0)} & -\dfrac{g_{12}(1,0)}{F(-g_{12}(1,0),g_{11}(1,0))}\\ 0 & \dfrac{g_{11}(1,0)}{F(-g_{12}(1,0),g_{11}(1,0))}\end{pmatrix} \tag{3.2.3}$$

## 3.3. Пример 1 пространства Минковского размерности 3

Мы рассмотрим пример финслеровой метрики пространства Минковского размерности 3 (рассматриваемая метрика является частным случаем метрики [12]-(3))

$$\begin{aligned}F^2(\overline{v}) &= ((v^1)^2+(v^2)^2+(v^3)^2)\left(1+\frac{(v^1)^2v^2v^3}{((v^1)^2+(v^2)^2+(v^3)^2)^2}\right)\\ &= (v^1)^2+(v^2)^2+(v^3)^2+\frac{(v^1)^2v^2v^3}{(v^1)^2+(v^2)^2+(v^3)^2}\end{aligned} \tag{3.3.1}$$

Теорема 3.3.1. *Финслеровая метрика, определённая равенством (3.3.1) порождает метрику пространства Минковского*

$$g_{11}(\overline{v}) = 1+\frac{-3(v^1)^2(v^2)^3v^3-3(v^1)^2v^2(v^3)^3+(v^2)^5v^3+2(v^2)^3(v^3)^3+v^2(v^3)^5}{((v^1)^2+(v^2)^2+(v^3)^2)^3} \tag{3.3.2}$$

$$g_{12}(\overline{v}) = \frac{3(v^1)^3(v^2)^2v^3-v^1(v^2)^4v^3+(v^1)^3(v^3)^3+v^1(v^3)^5}{((v^1)^2+(v^2)^2+(v^3)^2)^3} \tag{3.3.3}$$

$$g_{13}(\overline{v}) = \frac{(v^1)^3(v^2)^3+v^1(v^2)^5+3(v^1)^3v^2(v^3)^2-v^1v^2(v^3)^4}{((v^1)^2+(v^2)^2+(v^3)^2)^3} \tag{3.3.4}$$

$$g_{22}(\overline{v}) = 1+\frac{-3(v^1)^4v^2v^3+(v^1)^2(v^2)^3v^3-3(v^1)^2v^2(v^3)^3}{((v^1)^2+(v^2)^2+(v^3)^2)^3} \tag{3.3.5}$$

$$g_{23}(\overline{v}) = \frac{1}{2}\frac{(v^1)^6-(v^1)^2(v^3)^4-(v^1)^2(v^2)^4+6(v^1)^2(v^2)^2(v^3)^2}{((v^1)^2+(v^2)^2+(v^3)^2)^3} \tag{3.3.6}$$

$$g_{33}(\overline{v}) = 1+\frac{(v^1)^2v^2(v^3)^3-3(v^1)^4v^2v^3-3(v^1)^2(v^2)^3v^3}{((v^1)^2+(v^2)^2+(v^3)^2)^3} \tag{3.3.7}$$

Доказательство. Из равенства (3.3.1) следует

$$\begin{aligned}\frac{\partial F^2(\overline{v})}{\partial v^1} &= 2v^1+\frac{2v^1v^2v^3((v^1)^2+(v^2)^2+(v^3)^2)-(v^1)^2v^2v^3 2v^1}{((v^1)^2+(v^2)^2+(v^3)^2)^2}\\ &= 2v^1+\frac{2v^1(v^2)^3v^3+2v^1v^2(v^3)^3}{((v^1)^2+(v^2)^2+(v^3)^2)^2}\end{aligned} \tag{3.3.8}$$

$$\begin{aligned}\frac{\partial F^2(\overline{v})}{\partial v^2} &= 2v^2+\frac{(v^1)^2v^3((v^1)^2+(v^2)^2+(v^3)^2)-(v^1)^2v^2v^3 2v^2}{((v^1)^2+(v^2)^2+(v^3)^2)^2}\\ &= 2v^2+\frac{(v^1)^4v^3-(v^1)^2(v^2)^2v^3+(v^1)^2(v^3)^3}{((v^1)^2+(v^2)^2+(v^3)^2)^2}\end{aligned} \tag{3.3.9}$$

$$\begin{aligned}\frac{\partial F^2(\overline{v})}{\partial v^3} &= 2v^3+\frac{(v^1)^2v^2((v^1)^2+(v^2)^2+(v^3)^2)-(v^1)^2v^2v^3 2v^3}{((v^1)^2+(v^2)^2+(v^3)^2)^2}\\ &= 2v^3+\frac{(v^1)^4v^2+(v^1)^2(v^2)^3-(v^1)^2v^2(v^3)^2}{((v^1)^2+(v^2)^2+(v^3)^2)^2}\end{aligned} \tag{3.3.10}$$

Из равенств (2.2.2), (3.3.8), (3.3.9), (3.3.10) следует

$$
\left\{
\begin{aligned}
g_{11}(\overline{v}) &= 1+\frac{1}{2}\frac{(2(v^2)^3v^3+2v^2(v^3)^3)((v^1)^2+(v^2)^2+(v^3)^2)^2}{((v^1)^2+(v^2)^2+(v^3)^2)^4}\\
&-\frac{1}{2}\frac{(2v^1(v^2)^3v^3+2v^1v^2(v^3)^3)2((v^1)^2+(v^2)^2+(v^3)^2)2v^1}{((v^1)^2+(v^2)^2+(v^3)^2)^4}\\
&=1+\frac{((v^2)^3v^3+v^2(v^3)^3)((v^1)^2+(v^2)^2+(v^3)^2)}{((v^1)^2+(v^2)^2+(v^3)^2)^3}\\
&-\frac{(v^1(v^2)^3v^3+v^1v^2(v^3)^3)4v^1}{((v^1)^2+(v^2)^2+(v^3)^2)^3}\\
&=1+\frac{(v^1)^2(v^2)^3v^3+(v^1)^2v^2(v^3)^3+(v^2)^5v^3+(v^2)^3(v^3)^3}{((v^1)^2+(v^2)^2+(v^3)^2)^3}\\
&+\frac{(v^2)^3(v^3)^3+v^2(v^3)^3(v^3)^2}{((v^1)^2+(v^2)^2+(v^3)^2)^3}-\frac{4(v^1)^2(v^2)^3v^3+4(v^1)^2v^2(v^3)^3}{((v^1)^2+(v^2)^2+(v^3)^2)^3}
\end{aligned}
\right.
\tag{3.3.11}
$$

$$
\left\{
\begin{aligned}
&g_{12}(\overline{v})\\
&=\frac{1}{2}\frac{(4(v^1)^3v^3-2v^1(v^2)^2v^3+2v^1(v^3)^3)((v^1)^2+(v^2)^2+(v^3)^2)^2}{((v^1)^2+(v^2)^2+(v^3)^2)^4}\\
&-\frac{1}{2}\frac{((v^1)^4v^3-(v^1)^2(v^2)^2v^3+(v^1)^2(v^3)^3)2((v^1)^2+(v^2)^2+(v^3)^2)2v^1}{((v^1)^2+(v^2)^2+(v^3)^2)^4}\\
&=\frac{(2(v^1)^3v^3-v^1(v^2)^2v^3+v^1(v^3)^3)((v^1)^2+(v^2)^2+(v^3)^2)}{((v^1)^2+(v^2)^2+(v^3)^2)^3}\\
&-\frac{((v^1)^4v^3-(v^1)^2(v^2)^2v^3+(v^1)^2(v^3)^3)2v^1}{((v^1)^2+(v^2)^2+(v^3)^2)^3}\\
&=\frac{2(v^1)^5v^3-(v^1)^3(v^2)^2v^3+(v^1)^3(v^3)^3}{((v^1)^2+(v^2)^2+(v^3)^2)^3}\\
&+\frac{2(v^1)^3(v^2)^2v^3-v^1(v^2)^4v^3+v^1(v^2)^2(v^3)^3}{((v^1)^2+(v^2)^2+(v^3)^2)^3}\\
&+\frac{2(v^1)^3(v^3)^3-v^1(v^2)^2(v^3)^3+v^1(v^3)^5}{((v^1)^2+(v^2)^2+(v^3)^2)^3}\\
&-\frac{2(v^1)^5v^3-2(v^1)^3(v^2)^2v^3+2(v^1)^3(v^3)^3}{((v^1)^2+(v^2)^2+(v^3)^2)^3}
\end{aligned}
\right.
\tag{3.3.12}
$$

$$
(3.3.13)\quad
\left\{
\begin{aligned}
g_{13}(\overline{v}) &= \frac{1}{2}\frac{(4(v^1)^3v^2+2v^1(v^2)^3-2v^1v^2(v^3)^2)((v^1)^2+(v^2)^2+(v^3)^2)^2}{((v^1)^2+(v^2)^2+(v^3)^2)^4}\\
&-\frac{1}{2}\frac{((v^1)^4v^2+(v^1)^2(v^2)^3-(v^1)^2v^2(v^3)^2)2((v^1)^2+(v^2)^2+(v^3)^2)2v^1}{((v^1)^2+(v^2)^2+(v^3)^2)^4}\\
&=\frac{(2(v^1)^3v^2+v^1(v^2)^3-v^1v^2(v^3)^2)((v^1)^2+(v^2)^2+(v^3)^2)}{((v^1)^2+(v^2)^2+(v^3)^2)^3}\\
&-\frac{((v^1)^4v^2+(v^1)^2(v^2)^3-(v^1)^2v^2(v^3)^2)2v^1}{((v^1)^2+(v^2)^2+(v^3)^2)^3}\\
&=\frac{2(v^1)^5v^2+(v^1)^3(v^2)^3-(v^1)^3v^2(v^3)^2}{((v^1)^2+(v^2)^2+(v^3)^2)^3}+\frac{2(v^1)^3(v^2)^3+v^1(v^2)^5-v^1(v^2)^3(v^3)^2}{((v^1)^2+(v^2)^2+(v^3)^2)^3}\\
&+\frac{2(v^1)^3v^2(v^3)^2+v^1(v^2)^3(v^3)^2-v^1v^2(v^3)^4}{((v^1)^2+(v^2)^2+(v^3)^2)^3}-\frac{2(v^1)^5v^2+2(v^1)^3(v^2)^3-2(v^1)^3v^2(v^3)^2}{((v^1)^2+(v^2)^2+(v^3)^2)^3}
\end{aligned}
\right.
$$

$$
(3.3.14)\quad
\left\{
\begin{aligned}
g_{22}(\overline{v}) &= 1+\frac{1}{2}\frac{(-2(v^1)^2v^2v^3)((v^1)^2+(v^2)^2+(v^3)^2)^2}{((v^1)^2+(v^2)^2+(v^3)^2)^4}\\
&-\frac{1}{2}\frac{((v^1)^4v^3-(v^1)^2(v^2)^2v^3+(v^1)^2(v^3)^3)2((v^1)^2+(v^2)^2+(v^3)^2)2v^2}{((v^1)^2+(v^2)^2+(v^3)^2)^4}\\
&=1+\frac{-(v^1)^2v^2v^3((v^1)^2+(v^2)^2+(v^3)^2)}{((v^1)^2+(v^2)^2+(v^3)^2)^3}\\
&+\frac{(-(v^1)^4v^3+(v^1)^2(v^2)^2v^3-(v^1)^2(v^3)^3)2v^2}{((v^1)^2+(v^2)^2+(v^3)^2)^3}\\
&=1+\frac{-(v^1)^4v^2v^3-(v^1)^2(v^2)^3v^3-(v^1)^2v^2(v^3)^3}{((v^1)^2+(v^2)^2+(v^3)^2)^3}\\
&+\frac{-2(v^1)^4v^2v^3+2(v^1)^2(v^2)^3v^3-2(v^1)^2v^2(v^3)^3}{((v^1)^2+(v^2)^2+(v^3)^2)^3}
\end{aligned}
\right.
$$

$$
(3.3.15)\quad \left\{\begin{aligned}
g_{23}(\overline{v}) &= \frac{1}{2}\frac{((v^1)^4+3(v^1)^2(v^2)^2-(v^1)^2(v^3)^2)((v^1)^2+(v^2)^2+(v^3)^2)^2}{((v^1)^2+(v^2)^2+(v^3)^2)^4}\\
&-\frac{1}{2}\frac{((v^1)^4v^2+(v^1)^2(v^2)^3-(v^1)^2v^2(v^3)^2)2((v^1)^2+(v^2)^2+(v^3)^2)2v^2}{((v^1)^2+(v^2)^2+(v^3)^2)^4}\\
&=\frac{1}{2}\frac{((v^1)^4+3(v^1)^2(v^2)^2-(v^1)^2(v^3)^2)((v^1)^2+(v^2)^2+(v^3)^2)}{((v^1)^2+(v^2)^2+(v^3)^2)^3}\\
&+\frac{1}{2}\frac{(-(v^1)^4v^2-(v^1)^2(v^2)^3+(v^1)^2v^2(v^3)^2)4v^2}{((v^1)^2+(v^2)^2+(v^3)^2)^3}\\
&=\frac{1}{2}\frac{(v^1)^6+3(v^1)^4(v^2)^2-(v^1)^4(v^3)^2}{((v^1)^2+(v^2)^2+(v^3)^2)^3}\\
&+\frac{1}{2}\frac{(v^1)^4(v^2)^2+3(v^1)^2(v^2)^4-(v^1)^2(v^2)^2(v^3)^2}{((v^1)^2+(v^2)^2+(v^3)^2)^3}\\
&+\frac{1}{2}\frac{(v^1)^4(v^3)^2+3(v^1)^2(v^2)^2(v^3)^2-(v^1)^2(v^3)^4}{((v^1)^2+(v^2)^2+(v^3)^2)^3}\\
&+\frac{1}{2}\frac{-4(v^1)^4(v^2)^2-4(v^1)^2(v^2)^4+4(v^1)^2(v^2)^2(v^3)^2}{((v^1)^2+(v^2)^2+(v^3)^2)^3}
\end{aligned}\right.
$$

$$
(3.3.16)\quad \left\{\begin{aligned}
&g_{33}(\overline{v})\\
&=1+\frac{1}{2}\frac{(-2(v^1)^2v^2v^3)((v^1)^2+(v^2)^2+(v^3)^2)^2}{((v^1)^2+(v^2)^2+(v^3)^2)^4}\\
&-\frac{1}{2}\frac{((v^1)^4v^2+(v^1)^2(v^2)^3-(v^1)^2v^2(v^3)^2)2((v^1)^2+(v^2)^2+(v^3)^2)2v^3}{((v^1)^2+(v^2)^2+(v^3)^2)^4}\\
&=1+\frac{-(v^1)^2v^2v^3((v^1)^2+(v^2)^2+(v^3)^2)}{((v^1)^2+(v^2)^2+(v^3)^2)^4}\\
&+\frac{(-(v^1)^4v^2-(v^1)^2(v^2)^3+(v^1)^2v^2(v^3)^2)2v^3}{((v^1)^2+(v^2)^2+(v^3)^2)^3}\\
&=1+\frac{-(v^1)^4v^2v^3-(v^1)^2(v^2)^3v^3-(v^1)^2v^2(v^3)^3}{((v^1)^2+(v^2)^2+(v^3)^2)^4}\\
&+\frac{-2(v^1)^4v^2v^3-2(v^1)^2(v^2)^3v^3+2(v^1)^2v^2(v^3)^3}{((v^1)^2+(v^2)^2+(v^3)^2)^3}
\end{aligned}\right.
$$

Равенство (3.3.2) следует из равенства (3.3.11). Равенство (3.3.3) следует из равенства (3.3.12). Равенство (3.3.4) следует из равенства (3.3.13). Равенство (3.3.5) следует из равенства (3.3.14). Равенство (3.3.6) следует из равенства (3.3.15). Равенство (3.3.7) следует из равенства (3.3.16). □

Теорема 3.3.2. *Финслерова метрика, определённая равенством (3.3.1), порождает симметричное отношение ортогональности.*

Доказательство. Мы начнём процедуру ортогонализации[3.3] с базиса

$$\overline{e}_1=\begin{pmatrix}1\\0\\0\end{pmatrix}\quad\overline{e}_2=\begin{pmatrix}0\\1\\0\end{pmatrix}\quad\overline{e}_3=\begin{pmatrix}0\\0\\1\end{pmatrix}$$

Пусть

$$\overline{e}'_1=\overline{e}_1\quad\overline{e}'_2=\begin{pmatrix}e'^1_2\\e'^2_2\\e'^3_2\end{pmatrix}\quad\overline{e}'_3=\begin{pmatrix}e'^1_3\\e'^2_3\\e'^3_3\end{pmatrix}$$

ортогональный базис. Так как вектор $\overline{e}'_2$ ортогонален вектору $\overline{e}'_1$, то согласно определению 2.3.1

$$g_{ij}(\overline{e}'_1)e'^i_1e'^j_2=0$$

$$g_{1j}(\overline{e}'_1)e'^j_2=0\tag{3.3.17}$$

Согласно равенствам (3.3.2), (3.3.3), (3.3.4),

$$\begin{array}{lll}g_{11}(\overline{e}'_1)=1 & g_{12}(\overline{e}'_1)=0 & g_{13}(\overline{e}'_1)=0\\ & g_{22}(\overline{e}'_1)=1 & g_{23}(\overline{e}'_1)=\frac{1}{2}\\ & & g_{33}(\overline{e}'_1)=1\end{array}\tag{3.3.18}$$

Из равенств (3.3.17), (3.3.18) следует

$$e'^1_2=0\tag{3.3.19}$$

Поскольку значения $e'^2_2$, $e'^3_2$ произвольны, мы положим[3.4]

$$e'^2_2=a\quad e'^3_2=b\tag{3.3.20}$$

Согласно теореме 3.3.1 и равенствам (3.3.19), (3.3.20),

$$\begin{array}{lll}g_{11}(\overline{e}'_2)=1+\frac{a^5b+2a^3b^3+ab^5}{(a^2+b^2)^3}=1+\frac{ab}{a^2+b^2} & g_{12}(\overline{e}'_2)=0 & g_{13}(\overline{e}'_2)=0\\ & g_{22}(\overline{e}'_2)=1 & g_{23}(\overline{e}'_2)=0\\ & & g_{33}(\overline{e}'_2)=1\end{array}\tag{3.3.21}$$

Из равенств (3.3.20), (3.3.21) следует

$$g_{ij}(\overline{e}'_2)e'^i_2e'^j_1=g_{i1}(\overline{e}'_2)e'^i_2=0$$

Следовательно, вектор $\overline{e}'_1$ ортогонален вектору $\overline{e}'_2$.

Так как вектор $\overline{e}'_3$ ортогонален векторам $\overline{e}'_1$, $\overline{e}'_2$, то согласно определению 2.3.1

$$g_{ij}(\overline{e}'_1)e'^i_1e'^j_3=0\qquad g_{ij}(\overline{e}'_2)e'^i_2e'^j_3=0$$

$$g_{1j}(\overline{e}'_1)e'^j_3=0\qquad g_{2j}(\overline{e}'_2)ae'^j_3+g_{3j}(\overline{e}'_2)be'^j_3=0\tag{3.3.22}$$

---

[3.3] Смотри раздел [8]-3.

[3.4] Если мы положим $a=1$, $b=0$, то легко видеть, что вектор $\overline{e}_1$ ортогонален вектору $\overline{e}_2$ и вектор $\overline{e}_2$ ортогонален вектору $\overline{e}_1$.

Из равенств (3.3.18), (3.3.21), (3.3.22), следует

(3.3.23) $$e'^1_3 = 0 \quad ae'^2_3 + be'^3_3 = 0$$

Из равенства (3.3.23) следует, что мы можем положить

(3.3.24) $$e'^1_3 = 0 \quad e'^2_3 = -b \quad e'^3_3 = a$$

Согласно теореме 3.3.1 и равенству (3.3.24),

(3.3.25) $$\begin{aligned} & g_{11}(\overline{e}'_3) \\ & = 1 + \frac{(-b)^5 a + 2(-b)^3(a)^3 + (-b)(a)^5}{(a^2+b^2)^3} \\ & = 1 - \frac{ab}{a^2+b^2} \end{aligned} \quad g_{12}(\overline{e}'_3) = 0 \quad g_{13}(\overline{e}'_3) = 0$$
$$g_{22}(\overline{e}'_3) = 1 \quad g_{23}(\overline{e}'_3) = 0$$
$$g_{33}(\overline{e}'_3) = 1$$

Из равенств (3.3.20), (3.3.24), (3.3.25) следует

$$\begin{aligned} g_{ij}(\overline{e}'_3)e'^i_3 e'^j_1 &= g_{i1}(\overline{e}'_3)e'^i_3 = 0 \\ g_{ij}(\overline{e}'_3)e'^i_3 e'^j_2 &= g_{i2}(\overline{e}'_3)e'^i_3 a + g_{i3}(\overline{e}'_3)e'^i_3 b \\ &= g_{22}(\overline{e}'_3)e'^2_3 a + g_{33}(\overline{e}'_3)e'^3_3 b \\ &= (-b)a + ab = 0 \end{aligned}$$

Следовательно, векторы $\overline{e}'_1$, $\overline{e}'_2$ ортогональны вектору $\overline{e}'_3$. □

Поскольку

$$F^2(\overline{e}'_1) = 1$$

то длина вектора

$$\overline{e}''_1 = \overline{e}'_1 = \begin{pmatrix} 1 \\ 0 \\ 0 \end{pmatrix}$$

равна 1. Поскольку

$$F^2(\overline{e}'_2) = a^2 + b^2$$

то длина вектора

$$\overline{e}''_2 = \frac{1}{\sqrt{a^2+b^2}}\overline{e}'_2 = \begin{pmatrix} 0 \\ \frac{a}{\sqrt{a^2+b^2}} \\ \frac{b}{\sqrt{a^2+b^2}} \end{pmatrix}$$

равна 1. Поскольку

$$F^2(\overline{e}'_3) = (-b)^2 + (a)^2 = b^2 + a^2$$

то длина вектора

$$\overline{e}''_3=\frac{1}{\sqrt{a^2+b^2}}\overline{e}'_3=\begin{pmatrix}0\\ -\frac{b}{\sqrt{a^2+b^2}}\\ \frac{a}{\sqrt{a^2+b^2}}\end{pmatrix}$$

равна 1. Следовательно, базис $\overline{\overline{e}}''$ является ортонормированным базисом. Базисы $\overline{\overline{e}}''$ и $\overline{\overline{e}}$ связаны пассивным преобразованием

$$\begin{pmatrix}\overline{e}''_1 & \overline{e}''_2 & \overline{e}''_3\end{pmatrix}=\begin{pmatrix}\overline{e}_1 & \overline{e}_2 & \overline{e}_3\end{pmatrix}\begin{pmatrix}1 & 0 & 0\\ 0 & \frac{a}{\sqrt{a^2+b^2}} & -\frac{b}{\sqrt{a^2+b^2}}\\ 0 & \frac{b}{\sqrt{a^2+b^2}} & \frac{a}{\sqrt{a^2+b^2}}\end{pmatrix} \tag{3.3.26}$$

Теорема 3.3.3. *Компоненты метрического тензора относительно базиса $\overline{\overline{e}}''$ имеют вид*

$$g''(\overline{e}''_1)=\begin{pmatrix}1 & 0 & 0\\ 0 & 1+\frac{ab}{a^2+b^2} & \frac{1}{2}\frac{a^2-b^2}{a^2+b^2}\\ 0 & \frac{1}{2}\frac{a^2-b^2}{a^2+b^2} & 1-\frac{ab}{a^2+b^2}\end{pmatrix} \tag{3.3.27}$$

$$g''(\overline{e}''_2)=\begin{pmatrix}1+\frac{ab}{a^2+b^2} & 0 & 0\\ 0 & 1 & 0\\ 0 & 0 & 1\end{pmatrix} \tag{3.3.28}$$

$$g''(\overline{e}''_3)=\begin{pmatrix}1-\frac{ab}{a^2+b^2} & 0 & 0\\ 0 & 1 & 0\\ 0 & 0 & 1\end{pmatrix} \tag{3.3.29}$$

Доказательство. Из равенств (2.5.3), (3.3.18), (3.3.26), следует

$$
(3.3.30)\quad
\begin{cases}
g''_{11}(\overline{e}''_1)=g_{11}(\overline{e}''_1)=1\\
g''_{12}(\overline{e}''_1)=g_{12}(\overline{e}''_1)A_2^2+g_{13}(\overline{e}''_1)A_2^3=0\\
g''_{13}(\overline{e}''_1)=g_{12}(\overline{e}''_1)A_3^2+g_{13}(\overline{e}''_1)A_3^3=0\\
g''_{22}(\overline{e}''_1)=g_{22}(\overline{e}''_1)A_2^2A_2^2+2g_{23}(\overline{e}''_1)A_2^2A_2^3+g_{33}(\overline{e}''_1)A_2^3A_2^3\\
\qquad=\dfrac{a^2}{a^2+b^2}+\dfrac{ab}{a^2+b^2}+\dfrac{b^2}{a^2+b^2}=1+\dfrac{ab}{a^2+b^2}\\
g''_{23}(\overline{e}''_1)=g_{22}(\overline{e}''_1)A_2^2A_3^2+g_{23}(\overline{e}''_1)A_2^2A_3^3+g_{23}(\overline{e}''_1)A_3^2A_2^3+g_{33}(\overline{e}''_1)A_2^3A_3^3\\
\qquad=\dfrac{-ab}{a^2+b^2}+\dfrac{1}{2}\dfrac{a^2}{a^2+b^2}-\dfrac{1}{2}\dfrac{1}{2}\dfrac{b^2}{a^2+b^2}+\dfrac{ab}{a^2+b^2}=\dfrac{1}{2}\dfrac{a^2-b^2}{a^2+b^2}\\
g''_{33}(\overline{e}''_1)=g_{22}(\overline{e}''_1)A_3^2A_3^2+2g_{23}(\overline{e}''_1)A_3^2A_3^3+g_{33}(\overline{e}''_1)A_3^3A_3^3\\
\qquad=\dfrac{b^2}{a^2+b^2}-\dfrac{ab}{a^2+b^2}+\dfrac{a^2}{a^2+b^2}=1-\dfrac{ab}{a^2+b^2}
\end{cases}
$$

Равенство (3.3.27) следует из (3.3.30). Из равенств (2.5.3), (3.3.21), (3.3.26), следует

$$
(3.3.31)\quad
\begin{cases}
g''_{11}(\overline{e}''_2)=g_{11}(\overline{e}''_2)=1+\dfrac{ab}{a^2+b^2}\\
g''_{12}(\overline{e}''_2)=g_{11}(\overline{e}''_2)A_2^1+g_{12}(\overline{e}''_2)A_2^2+g_{13}(\overline{e}''_2)A_2^3=0\\
g''_{13}(\overline{e}''_2)=g_{11}(\overline{e}''_2)A_3^1+g_{12}(\overline{e}''_2)A_3^2+g_{13}(\overline{e}''_2)A_3^3=0\\
g''_{22}(\overline{e}''_2)=g_{22}(\overline{e}''_2)A_2^2A_2^2+2g_{23}(\overline{e}''_2)A_2^2A_2^3+g_{33}(\overline{e}''_2)A_2^3A_2^3\\
\qquad=\dfrac{a^2}{a^2+b^2}+\dfrac{b^2}{a^2+b^2}=1\\
g''_{23}(\overline{e}''_2)=g_{22}(\overline{e}''_2)A_2^2A_3^2+g_{33}(\overline{e}''_2)A_2^3A_3^3=0\\
g''_{33}(\overline{e}''_2)=g_{22}(\overline{e}''_2)A_3^2A_3^2+2g_{23}(\overline{e}''_2)A_3^2A_3^3+g_{33}(\overline{e}''_2)A_3^3A_3^3\\
\qquad=\dfrac{b^2}{a^2+b^2}+\dfrac{a^2}{a^2+b^2}=1
\end{cases}
$$

Равенство (3.3.28) следует из (3.3.31). Из равенств (2.5.3), (3.3.25), (3.3.26), следует

$$
(3.3.32)\quad
\begin{cases}
g''_{11}(\overline{e}''_3)=g_{11}(\overline{e}''_3)=1-\dfrac{ab}{a^2+b^2}\\
g''_{12}(\overline{e}''_3)=g_{11}(\overline{e}''_3)A_2^1+g_{12}(\overline{e}''_3)A_2^2+g_{13}(\overline{e}''_3)A_2^3=0\\
g''_{13}(\overline{e}''_3)=g_{11}(\overline{e}''_3)A_3^1+g_{12}(\overline{e}''_3)A_3^2+g_{13}(\overline{e}''_3)A_3^3=0\\
g''_{22}(\overline{e}''_3)=g_{22}(\overline{e}''_3)A_2^2A_2^2+2g_{23}(\overline{e}''_3)A_2^2A_2^3+g_{33}(\overline{e}''_3)A_2^3A_2^3\\
\qquad=\dfrac{a^2}{a^2+b^2}+\dfrac{b^2}{a^2+b^2}=1\\
g''_{23}(\overline{e}''_3)=g_{22}(\overline{e}''_3)A_2^2A_3^2+g_{33}(\overline{e}''_3)A_2^3A_3^3=0\\
g''_{33}(\overline{e}''_3)=g_{22}(\overline{e}''_3)A_3^2A_3^2+2g_{23}(\overline{e}''_3)A_3^2A_3^3+g_{33}(\overline{e}''_3)A_3^3A_3^3\\
\qquad=\dfrac{b^2}{a^2+b^2}+\dfrac{a^2}{a^2+b^2}=1
\end{cases}
$$

Равенство (3.3.29) следует из (3.3.32). □

Следствие 3.3.4. *Матрица (2.3.9), соответствующая финслеровой метрике, определённой равенством (3.3.1), является диагональной матрицей.* □

## 3.4. Пример 2 пространства Минковского размерности 3

Мы рассмотрим пример финслеровой метрики пространства Минковского размерности 3 (рассматриваемая метрика является частным случаем метрики [12]-(3))

$$\begin{aligned} F^2(\overline{v}) &= ((v^1)^2+(v^2)^2+(v^3)^2)\left(1+\frac{(v^1)^3v^2}{((v^1)^2+(v^2)^2+(v^3)^2)^2}\right) \\ &= (v^1)^2+(v^2)^2+(v^3)^2+\frac{(v^1)^3v^2}{(v^1)^2+(v^2)^2+(v^3)^2} \end{aligned} \tag{3.4.1}$$

Теорема 3.4.1. *Финслерова метрика $F^2$, определённая равенством (3.4.1) порождает метрику пространства Минковского*

$$\begin{aligned} g_{11}(\overline{v}) &= 1-\frac{(v^1)^3(v^2)^3+(v^1)^3v^2(v^3)^2}{((v^1)^2+(v^2)^2+(v^3)^2)^3} \\ &+\frac{6v^1(v^2)^3(v^3)^2+3v^1(v^2)^5+3v^1v^2(v^3)^4}{((v^1)^2+(v^2)^2+(v^3)^2)^3} \end{aligned} \tag{3.4.2}$$

$$g_{12}(\overline{v}) = \frac{1}{2}\frac{(v^1)^6+6(v^1)^4(v^2)^2-3(v^1)^2(v^2)^4+4(v^1)^4(v^3)^2+3(v^1)^2(v^3)^4}{((v^1)^2+(v^2)^2+(v^3)^2)^3} \tag{3.4.3}$$

$$g_{13}(\overline{v}) = \frac{(v^1)^4v^2v^3-3(v^1)^2(v^2)^3v^3-3(v^1)^2v^2(v^3)^3}{((v^1)^2+(v^2)^2+(v^3)^2)^3} \tag{3.4.4}$$

$$g_{22}(\overline{v}) = 1+\frac{(v^1)^3(v^2)^3-3(v^1)^5v^2-3(v^1)^3v^2(v^3)^2}{((v^1)^2+(v^2)^2+(v^3)^2)^3} \tag{3.4.5}$$

$$g_{23}(\overline{v}) = \frac{3(v^1)^3(v^2)^2v^3-(v^1)^5v^3-(v^1)^3(v^3)^3}{((v^1)^2+(v^2)^2+(v^3)^2)^3} \tag{3.4.6}$$

$$g_{33}(\overline{v}) = 1+\frac{3(v^1)^3v^2(v^3)^2-(v^1)^5v^2-(v^1)^3(v^2)^3}{((v^1)^2+(v^2)^2+(v^3)^2)^3} \tag{3.4.7}$$

Доказательство. Из равенства (3.4.1) следует

$$\begin{aligned} \frac{\partial F^2(\overline{v})}{\partial v^1} &= 2v^1+\frac{3(v^1)^2v^2((v^1)^2+(v^2)^2+(v^3)^2)-(v^1)^3v^22v^1}{((v^1)^2+(v^2)^2+(v^3)^2)^2} \\ &= 2v^1+\frac{3(v^1)^4v^2+3(v^1)^2(v^2)^3+3(v^1)^2v^2(v^3)^2-2(v^1)^4v^2}{((v^1)^2+(v^2)^2+(v^3)^2)^2} \\ &= 2v^1+\frac{(v^1)^4v^2+3(v^1)^2(v^2)^3+3(v^1)^2v^2(v^3)^2}{((v^1)^2+(v^2)^2+(v^3)^2)^2} \end{aligned} \tag{3.4.8}$$

$$\begin{aligned} \frac{\partial F^2(\overline{v})}{\partial v^2} &= 2v^2+\frac{(v^1)^3((v^1)^2+(v^2)^2+(v^3)^2)-(v^1)^3v^22v^2}{((v^1)^2+(v^2)^2+(v^3)^2)^2} \\ &= 2v^2+\frac{(v^1)^5+(v^1)^3(v^2)^2+(v^1)^3(v^3)^2-2(v^1)^3(v^2)^2}{((v^1)^2+(v^2)^2+(v^3)^2)^2} \\ &= 2v^2+\frac{(v^1)^5+(v^1)^3(v^3)^2-(v^1)^3(v^2)^2}{((v^1)^2+(v^2)^2+(v^3)^2)^2} \end{aligned} \tag{3.4.9}$$

$$\frac{\partial F^2(\overline{v})}{\partial v^3} = 2v^3 + \frac{-2(v^1)^3v^2v^3}{((v^1)^2+(v^2)^2+(v^3)^2)^2} \tag{3.4.10}$$

Из равенств (2.2.2), (3.4.8), (3.4.9), (3.4.10) следует

$$\left\{\begin{aligned}
&g_{11}(\overline{v})\\
&=1+\frac{1}{2}\frac{(4(v^1)^3v^2+6v^1(v^2)^3+6v^1v^2(v^3)^2)((v^1)^2+(v^2)^2+(v^3)^2)^2}{((v^1)^2+(v^2)^2+(v^3)^2)^4}\\
&-\frac{1}{2}\frac{(v^1)^4v^2+3(v^1)^2(v^2)^3+3(v^1)^2v^2(v^3)^2}{((v^1)^2+(v^2)^2+(v^3)^2)^4}\\
&*2((v^1)^2+(v^2)^2+(v^3)^2)2v^1\\
&=1+\frac{(2(v^1)^3v^2+3v^1(v^2)^3+3v^1v^2(v^3)^2)((v^1)^2+(v^2)^2+(v^3)^2)}{((v^1)^2+(v^2)^2+(v^3)^2)^3}\\
&-\frac{((v^1)^4v^2+3(v^1)^2(v^2)^3+3(v^1)^2v^2(v^3)^2)2v^1}{((v^1)^2+(v^2)^2+(v^3)^2)^3}\\
&=1+\frac{2(v^1)^5v^2+3(v^1)^3(v^2)^3+3(v^1)^3v^2(v^3)^2}{((v^1)^2+(v^2)^2+(v^3)^2)^3}\\
&+\frac{2(v^1)^3(v^2)^3+3v^1(v^2)^5+3v^1(v^2)^3(v^3)^2}{((v^1)^2+(v^2)^2+(v^3)^2)^3}\\
&+\frac{2(v^1)^3v^2(v^3)^2+3v^1(v^2)^3(v^3)^2+3v^1v^2(v^3)^4}{((v^1)^2+(v^2)^2+(v^3)^2)^3}\\
&-\frac{2(v^1)^5v^2+6(v^1)^3(v^2)^3+6(v^1)^3v^2(v^3)^2}{((v^1)^2+(v^2)^2+(v^3)^2)^3}
\end{aligned}\right. \tag{3.4.11}$$

$$\left\{\begin{aligned}
&g_{12}(\overline{v})\\
&=\frac{1}{2}\frac{(5(v^1)^4+3(v^1)^2(v^3)^2-3(v^1)^2(v^2)^2)((v^1)^2+(v^2)^2+(v^3)^2)^2}{((v^1)^2+(v^2)^2+(v^3)^2)^4}\\
&-\frac{1}{2}\frac{((v^1)^5+(v^1)^3(v^3)^2-(v^1)^3(v^2)^2)2((v^1)^2+(v^2)^2+(v^3)^2)2v^1}{((v^1)^2+(v^2)^2+(v^3)^2)^4}\\
&=\frac{1}{2}\frac{(5(v^1)^4+3(v^1)^2(v^3)^2-3(v^1)^2(v^2)^2)((v^1)^2+(v^2)^2+(v^3)^2)}{((v^1)^2+(v^2)^2+(v^3)^2)^3}\\
&-\frac{1}{2}\frac{((v^1)^5+(v^1)^3(v^3)^2-(v^1)^3(v^2)^2)4v^1}{((v^1)^2+(v^2)^2+(v^3)^2)^3}\\
&=\frac{1}{2}\frac{5(v^1)^6+3(v^1)^4(v^3)^2-3(v^1)^4(v^2)^2}{((v^1)^2+(v^2)^2+(v^3)^2)^3}\\
&+\frac{1}{2}\frac{5(v^1)^4(v^2)^2+3(v^1)^2(v^2)^2(v^3)^2-3(v^1)^2(v^2)^4}{((v^1)^2+(v^2)^2+(v^3)^2)^3}\\
&+\frac{1}{2}\frac{5(v^1)^4(v^3)^2+3(v^1)^2(v^3)^4-3(v^1)^2(v^2)^2(v^3)^2}{((v^1)^2+(v^2)^2+(v^3)^2)^3}\\
&+\frac{1}{2}\frac{-4(v^1)^6-4(v^1)^4(v^3)^2+4(v^1)^4(v^2)^2}{((v^1)^2+(v^2)^2+(v^3)^2)^3}
\end{aligned}\right. \tag{3.4.12}$$

$$
(3.4.13)\quad \begin{cases}
g_{13}(\overline{v}) = \dfrac{1}{2}\dfrac{(-6(v^1)^2v^2v^3)((v^1)^2+(v^2)^2+(v^3)^2)^2}{((v^1)^2+(v^2)^2+(v^3)^2)^4} \\
\quad -\dfrac{1}{2}\dfrac{(-2(v^1)^3v^2v^3)2((v^1)^2+(v^2)^2+(v^3)^2)2v^1}{((v^1)^2+(v^2)^2+(v^3)^2)^4} \\
\quad = \dfrac{(-3(v^1)^2v^2v^3)((v^1)^2+(v^2)^2+(v^3)^2)+((v^1)^3v^2v^3)4v^1}{((v^1)^2+(v^2)^2+(v^3)^2)^3} \\
\quad = \dfrac{-3(v^1)^4v^2v^3-3(v^1)^2(v^2)^3v^3-3(v^1)^2v^2(v^3)^3)+4(v^1)^4v^2v^3}{((v^1)^2+(v^2)^2+(v^3)^2)^3}
\end{cases}
$$

$$
(3.4.14)\quad \begin{cases}
g_{22}(\overline{v}) = 1+\dfrac{1}{2}\dfrac{(-2(v^1)^3v^2)((v^1)^2+(v^2)^2+(v^3)^2)^2}{((v^1)^2+(v^2)^2+(v^3)^2)^4} \\
\quad -\dfrac{1}{2}\dfrac{((v^1)^5+(v^1)^3(v^3)^2-(v^1)^3(v^2)^2)2((v^1)^2+(v^2)^2+(v^3)^2)2v^2}{((v^1)^2+(v^2)^2+(v^3)^2)^4} \\
\quad = 1+\dfrac{(-(v^1)^3v^2)((v^1)^2+(v^2)^2+(v^3)^2)}{((v^1)^2+(v^2)^2+(v^3)^2)^3} \\
\quad -\dfrac{((v^1)^5+(v^1)^3(v^3)^2-(v^1)^3(v^2)^2)2v^2}{((v^1)^2+(v^2)^2+(v^3)^2)^3} \\
\quad = 1+\dfrac{-(v^1)^5v^2-(v^1)^3(v^2)^3-(v^1)^3v^2(v^3)^2}{((v^1)^2+(v^2)^2+(v^3)^2)^3} \\
\quad +\dfrac{(-2(v^1)^5v^2-2(v^1)^3v^2(v^3)^2+2(v^1)^3(v^2)^3}{((v^1)^2+(v^2)^2+(v^3)^2)^3}
\end{cases}
$$

$$
(3.4.15)\quad \begin{cases}
g_{23}(\overline{v}) \\
= \dfrac{1}{2}\dfrac{(-2(v^1)^3v^3)((v^1)^2+(v^2)^2+(v^3)^2)^2}{((v^1)^2+(v^2)^2+(v^3)^2)^4} \\
-\dfrac{1}{2}\dfrac{(-2(v^1)^3v^2v^3)2((v^1)^2+(v^2)^2+(v^3)^2)2v^2}{((v^1)^2+(v^2)^2+(v^3)^2)^4} \\
= \dfrac{-(v^1)^3v^3((v^1)^2+(v^2)^2+(v^3)^2)+2(v^1)^3v^2v^32v^2}{((v^1)^2+(v^2)^2+(v^3)^2)^3} \\
= \dfrac{-(v^1)^5v^3-(v^1)^3(v^2)^2v^3-(v^1)^3(v^3)^3+4(v^1)^3(v^2)^2v^3}{((v^1)^2+(v^2)^2+(v^3)^2)^3}
\end{cases}
$$

$$
(3.4.16)\quad \begin{cases}
g_{33}(\overline{v}) \\
= 1+\dfrac{1}{2}\dfrac{(-2(v^1)^3v^2)((v^1)^2+(v^2)^2+(v^3)^2)^2}{((v^1)^2+(v^2)^2+(v^3)^2)^4} \\
-\dfrac{1}{2}\dfrac{(-2(v^1)^3v^2v^3)2((v^1)^2+(v^2)^2+(v^3)^2)2v^3}{((v^1)^2+(v^2)^2+(v^3)^2)^4} \\
= 1+\dfrac{(-(v^1)^3v^2)((v^1)^2+(v^2)^2+(v^3)^2)+(2(v^1)^3v^2v^3)2v^3}{((v^1)^2+(v^2)^2+(v^3)^2)^3} \\
= 1+\dfrac{-(v^1)^5v^2-(v^1)^3(v^2)^3-(v^1)^3v^2(v^3)^2+4(v^1)^3v^2(v^3)^2}{((v^1)^2+(v^2)^2+(v^3)^2)^3}
\end{cases}
$$

Равенство (3.4.2) следует из равенства (3.4.11). Равенство (3.4.3) следует из равенства (3.4.12). Равенство (3.4.4) следует из равенства (3.4.13). Равенство (3.4.5) следует из равенства (3.4.14). Равенство (3.4.6) следует из равенства (3.4.15). Равенство (3.4.7) следует из равенства (3.4.16). □

Теорема 3.4.2. *Финслерова метрика $F^2$, определённая равенством* (3.4.1)*, порождает несимметричное отношение ортогональности. Однако существуют векторы $\overline{v}$, $\overline{w}$, для которых отношение ортогональности симметрично.*

Доказательство. Мы начнём процедуру ортогонализации[3.5] с базиса

$$\overline{e}_1=\begin{pmatrix}1\\0\\0\end{pmatrix}\quad \overline{e}_2=\begin{pmatrix}0\\1\\0\end{pmatrix}\quad \overline{e}_3=\begin{pmatrix}0\\0\\1\end{pmatrix}$$

Пусть

$$\overline{e}'_1=\overline{e}_1\quad \overline{e}'_2=\begin{pmatrix}e'^1_2\\e'^2_2\\e'^3_2\end{pmatrix}\quad \overline{e}'_3=\begin{pmatrix}e'^1_3\\e'^2_3\\e'^3_3\end{pmatrix}$$

ортогональный базис. Так как вектор $\overline{e}'_2$ ортогонален вектору $\overline{e}'_1$, то согласно определению 2.3.1

$$g_{ij}(\overline{e}'_1)e'^i_1e'^j_2=0$$

$$g_{1j}(\overline{e}'_1)e'^j_2=0 \tag{3.4.17}$$

Согласно равенствам (3.4.2), (3.4.3), (3.4.4),

$$\begin{array}{ccc} g_{11}(\overline{e}'_1)=1 & g_{12}(\overline{e}'_1)=\frac{1}{2} & g_{13}(\overline{e}'_1)=0\\ & g_{22}(\overline{e}'_1)=1 & g_{23}(\overline{e}'_1)=0\\ & & g_{33}(\overline{e}'_1)=1 \end{array} \tag{3.4.18}$$

Из равенств (3.4.17), (3.4.18) следует

$$e'^1_2+\frac{1}{2}e'^2_2=0 \tag{3.4.19}$$

Поскольку значение $e'^3_2$ произвольно, мы положим[3.6]

$$e'^1_2=-a\quad e'^2_2=2a\quad e'^3_2=b \tag{3.4.20}$$

Координаты вектора $\overline{e}'_2$ определены с точностью до произвольного множителя. Поэтому, не нарушая общности, мы можем поделить координаты вектора $\overline{e}'_2$ на $a$. При этом мы должны рассмотреть два случая. Мы рассмотрим процедуру ортогонализации, когда $a=0$, в подразделе 3.4.1. Мы увидим, что в этом случае отношение ортогонализации коммутативно. Мы рассмотрим процедуру ортогонализации, когда $a\neq 0$, в подразделе 3.4.2. Мы увидим, что в этом случае отношение ортогонализации некоммутативно. □

[3.5] Смотри раздел [8]-3.

[3.6] Если мы положим $a=1$, $b=0$, то легко видеть, что вектор $\overline{e}_1$ ортогонален вектору $\overline{e}_2$ и вектор $\overline{e}_2$ ортогонален вектору $\overline{e}_1$.

**3.4.1. Процедура ортогонализации базиса, $(a = 0)$.** Поскольку $b$ произвольно, мы можем положить $b = 1$. Следовательно, вектор $\overline{e}'_2$ имеет вид

$$\overline{e}'_2 = \begin{pmatrix} 0 \\ 0 \\ 1 \end{pmatrix} \tag{3.4.21}$$

Согласно теореме 3.4.1 и равенству (3.4.21),

$$\begin{array}{lll} g_{11}(\overline{e}'_2) = 1 & g_{12}(\overline{e}'_2) = 0 & g_{13}(\overline{e}'_2) = 0 \\ & g_{22}(\overline{e}'_2) = 1 & g_{23}(\overline{e}'_2) = 0 \\ & & g_{33}(\overline{e}'_2) = 1 \end{array} \tag{3.4.22}$$

Из равенств (3.4.20), (3.4.22) следует

$$g_{ij}(\overline{e}'_2){e'}^i_2{e'}^j_1 = g_{22}(\overline{e}'_2){e'}^2_1 = 0$$

Следовательно, вектор $\overline{e}'_1$ ортогонален вектору $\overline{e}'_2$.

Так как вектор $\overline{e}'_3$ ортогонален векторам $\overline{e}'_1$, $\overline{e}'_2$, то согласно определению 2.3.1

$$\begin{array}{ll} g_{ij}(\overline{e}'_1){e'}^i_1{e'}^j_3 = 0 & g_{ij}(\overline{e}'_2){e'}^i_2{e'}^j_3 = 0 \\ g_{1j}(\overline{e}'_1){e'}^j_3 = 0 & g_{3j}(\overline{e}'_2){e'}^j_3 = 0 \end{array} \tag{3.4.23}$$

Из равенств (3.4.18) (3.4.22), (3.4.23), следует

$${e'}^1_3 + \frac{1}{2}{e'}^2_3 = 0 \quad {e'}^3_3 = 0 \tag{3.4.24}$$

Из равенства (3.4.24) следует, что мы можем положить

$${e'}^1_3 = -1 \quad {e'}^2_3 = 2 \quad {e'}^3_3 = 0 \tag{3.4.25}$$

Согласно теореме 3.4.1 и равенству (3.4.25),

$$\begin{aligned} g_{11}(\overline{e}'_3) &= 1 + \frac{-(-1)^3 2^3 + 3(-1)2^5}{((-1)^2 + 2^2)^3} \\ g_{12}(\overline{e}'_3) &= \frac{1}{2}\frac{(-1)^6 + 6(-1)^4 2^2 - 3(-1)^2 2^4}{((-1)^2 + 2^2)^3} \\ g_{13}(\overline{e}'_3) &= 0 \\ g_{22}(\overline{e}'_3) &= 1 + \frac{(-1)^3(2)^3 - 3(-1)^5 2}{((-1)^2 + 2^2)^3} \\ g_{23}(\overline{e}'_3) &= 0 \\ g_{33}(\overline{e}'_3) &= 1 + \frac{-(-1)^5 2 - (-1)^3 2^3}{((-1)^2 + 2^2)^3} \end{aligned} \tag{3.4.26}$$

Из равенства (3.4.26) следует

$$g(\overline{e}'_3) = \begin{pmatrix} \frac{37}{125} & -\frac{23}{250} & 0 \\ & \frac{123}{125} & 0 \\ & & \frac{27}{25} \end{pmatrix} \tag{3.4.27}$$

Из равенств (3.4.20), (3.4.25), (3.4.27) следует

$$\begin{aligned} g_{ij}(\overline{e}'_3)e'^i_3 e'^j_1 &= g_{11}(\overline{e}'_3)e'^1_3 + g_{21}(\overline{e}'_3)e'^2_3 \\ &= \frac{37}{125}(-1) - \frac{23}{250}2 = -\frac{12}{25} \\ g_{ij}(\overline{e}'_3)e'^i_3 e'^j_2 &= g_{33}(\overline{e}'_3)e'^3_3 = 0 \end{aligned}$$

Следовательно, вектор $\overline{e}'_2$ ортогонален вектору $\overline{e}'_3$ и вектор $\overline{e}'_1$ не ортогонален вектору $\overline{e}'_3$.

Поскольку

$$F^2(\overline{e}'_1) = 1$$

то длина вектора

$$\overline{e}''_1 = \overline{e}'_1 = \begin{pmatrix} 1 \\ 0 \\ 0 \end{pmatrix}$$

равна 1. Поскольку

$$F^2(\overline{e}'_2) = 1$$

то длина вектора

$$\overline{e}''_2 = \overline{e}'_2 = \begin{pmatrix} 0 \\ 0 \\ 1 \end{pmatrix}$$

равна 1. Поскольку

$$F^2(\overline{e}'_3) = (-1)^2 + 2^2 + \frac{(-1)^3 2}{(-1)^2 + 2^2} = \frac{23}{5}$$

то длина вектора

$$\overline{e}''_3 = \sqrt{\frac{5}{23}}\overline{e}'_3 = \begin{pmatrix} -\sqrt{\frac{5}{23}} \\ 2\sqrt{\frac{5}{23}} \\ 0 \end{pmatrix}$$

равна 1. Следовательно, базис $\overline{\overline{e}}''$ является ортонормированным базисом. Базисы $\overline{\overline{e}}''$ и $\overline{\overline{e}}$ связаны пассивным преобразованием

$$\begin{pmatrix} \overline{e}''_1 & \overline{e}''_2 & \overline{e}''_3 \end{pmatrix} = \begin{pmatrix} \overline{e}_1 & \overline{e}_2 & \overline{e}_3 \end{pmatrix} \begin{pmatrix} 1 & 0 & -\sqrt{\frac{5}{23}} \\ 0 & 0 & 2\sqrt{\frac{5}{23}} \\ 0 & 1 & 0 \end{pmatrix} \tag{3.4.28}$$

Теорема 3.4.3. *Компоненты метрического тензора относительно базиса $\overline{\overline{e}}''$ имеют вид*

$$g''(\overline{e}''_1)=\begin{pmatrix}1&0&0\\0&1&0\\0&0&\frac{20}{23}\end{pmatrix}\tag{3.4.29}$$

$$g''(\overline{e}''_2)=\begin{pmatrix}1&0&-\sqrt{\frac{5}{23}}\\0&1&0\\-\sqrt{\frac{5}{23}}&0&\frac{25}{23}\end{pmatrix}\tag{3.4.30}$$

$$g''(\overline{e}''_3)=\begin{pmatrix}\frac{37}{125}&\frac{123}{125}&-\frac{12}{25}\sqrt{\frac{5}{23}}\\\frac{123}{125}&\frac{27}{25}&0\\-\frac{12}{25}\sqrt{\frac{5}{23}}&0&1\end{pmatrix}\tag{3.4.31}$$

Доказательство. Из равенств (2.5.3), (3.4.18), (3.4.28), следует

$$\left\{\begin{array}{l}g''_{11}(\overline{e}''_1)=g_{11}(\overline{e}''_1)=1\\g''_{12}(\overline{e}''_1)=g_{13}(\overline{e}''_1)A^1_3=0\\g''_{13}(\overline{e}''_1)=g_{11}(\overline{e}''_1)A^1_3+g_{12}(\overline{e}''_1)A^2_3=-\sqrt{\frac{5}{23}}+\frac{1}{2}2\sqrt{\frac{5}{23}}=0\\g''_{22}(\overline{e}''_1)=g_{33}(\overline{e}''_1)A^3_2A^3_2=1\\g''_{23}(\overline{e}''_1)=g_{13}(\overline{e}''_1)A^3_2A^1_3+g_{23}(\overline{e}''_1)A^3_2A^2_3=0\\g''_{33}(\overline{e}''_1)=g_{11}(\overline{e}''_1)A^1_3A^1_3+2g_{12}(\overline{e}''_1)A^1_3A^2_3+g_{22}(\overline{e}''_1)A^2_3A^2_3\\\qquad=\frac{5}{23}+2*\frac{1}{2}\left(-\sqrt{\frac{5}{23}}\right)2\sqrt{\frac{5}{23}}+4*\frac{5}{23}=\frac{20}{23}\end{array}\right.\tag{3.4.32}$$

Равенство (3.4.29) следует из (3.4.32). Из равенств (2.5.3), (3.4.22), (3.4.28), следует

$$\left\{\begin{array}{l}g''_{11}(\overline{e}''_2)=1\\g''_{12}(\overline{e}''_2)=g_{13}(\overline{e}''_2)A^3_2=0\\g''_{13}(\overline{e}''_2)=g_{11}(\overline{e}''_2)A^1_3+g_{12}(\overline{e}''_2)A^2_3=-\sqrt{\frac{5}{23}}\\g''_{22}(\overline{e}''_2)=g_{33}(\overline{e}''_2)A^3_2A^3_2=1\\g''_{23}(\overline{e}''_2)=g_{31}(\overline{e}''_2)A^3_2A^1_3+g_{32}(\overline{e}''_2)A^3_2A^2_3=0\\g''_{33}(\overline{e}''_2)=g_{11}(\overline{e}''_2)A^1_3A^1_3+2g_{12}(\overline{e}''_2)A^1_3A^2_3+g_{22}(\overline{e}''_2)A^2_3A^2_3\\\qquad=\frac{5}{23}+4*\frac{5}{23}=\frac{25}{23}\end{array}\right.\tag{3.4.33}$$

Равенство (3.4.30) следует из (3.4.33). Из равенств (2.5.3), (3.4.27), (3.4.28), следует

$$
\begin{cases}
g''_{11}(\overline{e}''_3) = g_{11}(\overline{e}''_3) = \frac{37}{125}\\
g''_{12}(\overline{e}''_3) = g_{13}(\overline{e}''_3)A^3_2 = \frac{123}{125}\\
g''_{13}(\overline{e}''_3) = g_{11}(\overline{e}''_3)A^1_3 + g_{12}(\overline{e}''_3)A^2_3 = -\frac{37}{125}\sqrt{\frac{5}{23}} - \frac{23}{250}2\sqrt{\frac{5}{23}}\\
g''_{22}(\overline{e}''_3) = g_{33}(\overline{e}''_3)A^3_2A^3_2\\
\qquad = \frac{27}{25}\\
g''_{23}(\overline{e}''_3) = g_{31}(\overline{e}''_3)A^3_2A^1_3 + g_{32}(\overline{e}''_3)A^3_2A^2_3 = 0\\
g''_{33}(\overline{e}''_3) = g_{11}(\overline{e}''_3)A^1_3A^1_3 + 2g_{12}(\overline{e}''_3)A^1_3A^2_3 + g_{22}(\overline{e}''_3)A^2_3A^2_3\\
\qquad = \frac{37}{125}\frac{5}{23} - 2(-\frac{23}{250})2 * \frac{5}{23} + \frac{123}{125}4 * \frac{5}{23} = 1
\end{cases}
\tag{3.4.34}
$$

Равенство (3.4.31) следует из (3.4.34). □

Следствие 3.4.4. *Относительно базиса $\overline{\overline{e}}''$, определённого равенством* (3.4.28), *матрица* (2.3.9), *соответствующая финслеровой метрике $F^2$, определённой равенством* (3.4.1), *имеет вид*

$$
\begin{pmatrix}
1 & 0 & 0\\
0 & 1 & 0\\
-\frac{12}{25}\sqrt{\frac{5}{23}} & 0 & 1
\end{pmatrix}
\tag{3.4.35}
$$

□

**3.4.2. Процедура ортогонализации базиса, ($a \ne 0$).** Поскольку $b$ произвольно, мы можем положить $a = 1$. Следовательно, вектор $\overline{e}'_2$ имеет вид

$$
\overline{e}'_2 = \begin{pmatrix} -1\\ 2\\ b \end{pmatrix}
\tag{3.4.36}
$$

Согласно теореме 3.4.1 и равенству (3.4.36),

$$
\begin{aligned}
g_{11}(\overline{e}'_2) &= 1 - \frac{(-1)^3 2^3 + (-1)^3 2b^2}{((-1)^2 + 2^2 + b^2)^3} + \frac{6(-1)2^3b^2 + 3(-1)2^5 + 3(-1)2b^4}{((-1)^2 + 2^2 + b^2)^3}\\
g_{12}(\overline{e}'_2) &= \frac{1}{2}\frac{(-1)^6 + 6(-1)^4 2^2 - 3(-1)^2 2^4 + 4(-1)^4 b^2 + 3(-1)^2 b^4}{((-1)^2 + 2^2 + b^2)^3}\\
g_{13}(\overline{e}'_2) &= \frac{(-1)^4 2b - 3(-1)^2 2^3 b - 3(-1)^2 2b^3}{((-1)^2 + 2^2 + b^2)^3}\\
g_{22}(\overline{e}'_2) &= 1 + \frac{(-1)^3 2^3 - 3(-1)^5 2 - 3(-1)^3 2b^2}{((-1)^2 + 2^2 + b^2)^3}\\
g_{23}(\overline{e}'_2) &= \frac{3(-1)^3 2^2 b - (-1)^5 b - (-1)^3 b^3}{((-1)^2 + 2^2 + b^2)^3}\\
g_{33}(\overline{e}'_2) &= 1 + \frac{3(-1)^3 2b^2 - (-1)^5 2 - (-1)^3 2^3}{((-1)^2 + 2^2 + b^2)^3}
\end{aligned}
\tag{3.4.37}
$$

Из равенств (3.4.37) следует

(3.4.38)
$$g(\overline{e}'_2)=\begin{pmatrix}1-\frac{88+46b^2+6b^4}{(5+b^2)^3}&\frac{-23+4b^2+3b^4}{2(5+b^2)^3}&\frac{-22b-6b^3}{(5+b^2)^3}\\&1+\frac{-2+6b^2}{(5+b^2)^3}&\frac{-11b+b^3}{(5+b^2)^3}\\&&1+\frac{-6b^2+10}{(5+b^2)^3}\end{pmatrix}$$

Из равенств (3.4.20), (3.4.38) следует

$$\begin{aligned}&g_{ij}(\overline{e}'_2){e'}^i_2{e'}^j_1]=g_{i1}(\overline{e}'_2){e'}^i_2\\
=&(1-\frac{88+46b^2+6b^4}{(5+b^2)^3})(-1)+\frac{-23+4b^2+3b^4}{2(5+b^2)^3}2+\frac{-22b-6b^3}{(5+b^2)^3}b\\
=&-1+\frac{88+46b^2+6b^4}{(5+b^2)^3}+\frac{-23+4b^2+3b^4}{(5+b^2)^3}+\frac{-22b^2-6b^4}{(5+b^2)^3}\\
=&-1+\frac{65+28b^2+3b^4}{(5+b^2)^3}=-1+\frac{13+3b^2}{(5+b^2)^2}=\frac{-b^4-10b^2-25+13+3b^2}{(5+b^2)^2}\\
=&\frac{-b^4-7b^2-12}{(5+b^2)^2}\end{aligned}$$

Полученное выражение отрицательно для любого значения $b$. Следовательно, вектор $\overline{e}'_1$ не ортогонален вектору $\overline{e}'_2$.

Так как вектор $\overline{e}'_3$ ортогонален векторам $\overline{e}'_1$, $\overline{e}'_2$, то согласно определению 2.3.1

(3.4.39)
$$\begin{aligned}&g_{ij}(\overline{e}'_1){e'}^i_1{e'}^j_3=0 && g_{ij}(\overline{e}'_2){e'}^i_2{e'}^j_3=0\\
&g_{1j}(\overline{e}'_1){e'}^j_3=0 && -g_{1j}(\overline{e}'_2){e'}^j_3+2g_{2j}(\overline{e}'_2){e'}^j_3+bg_{3j}(\overline{e}'_2){e'}^j_3=0\end{aligned}$$

Из равенств (3.4.18) (3.4.38), (3.4.39), следует

(3.4.40)
$${e'}^1_3+\frac{1}{2}{e'}^2_3=0$$

(3.4.41)
$$\begin{aligned}&-(1-\frac{88+46b^2+6b^4}{(5+b^2)^3}){e'}^1_3-\frac{-23+4b^2+3b^4}{2(5+b^2)^3}{e'}^2_3-\frac{-22b-6b^3}{(5+b^2)^3}{e'}^3_3\\
&+2\frac{-23+4b^2+3b^4}{2(5+b^2)^3}{e'}^1_3+2(1+\frac{-2+6b^2}{(5+b^2)^3}){e'}^2_3+2\frac{-11b+b^3}{(5+b^2)^3}{e'}^3_3\\
&+b\frac{-22b-6b^3}{(5+b^2)^3}{e'}^1_3+b\frac{-11b+b^3}{(5+b^2)^3}{e'}^2_3+b(1+\frac{-6b^2+10}{(5+b^2)^3}){e'}^3_3=0\end{aligned}$$

Из равенства (3.4.40) следует, что мы можем положить

(3.4.42)
$${e'}^1_3=-c\quad {e'}^2_3=2c$$

Из равенств (3.4.41), (3.4.42) следует

$$
\begin{aligned}
&(1-\frac{88}{5^3})c+\frac{-23+4b^2+3b^4}{2(5+b^2)^3}2c-2\frac{-23+4b^2+3b^4}{2(5+b^2)^3}c\\
&+4(1+\frac{-2+6b^2}{(5+b^2)^3})c-b\frac{-22b-6b^3}{(5+b^2)^3}c+2b\frac{-11b+b^3}{(5+b^2)^3}c\\
&=\frac{-22b-6b^3}{(5+b^2)^3}e'^3_3-2\frac{-11b+b^3}{(5+b^2)^3}e'^3_3-b(1+\frac{-6b^2+10}{(5+b^2)^3})e'^3_3
\end{aligned}
\tag{3.4.43}
$$

Пусть $b\ne 0$.[3.7] Из равенства (3.4.43) следует

$$(5+\frac{-50-30b^2-4b^4}{(5+b^2)^3})c=(-b+\frac{-2b^3-10b}{(5+b^2)^3})e'^3_3$$

$$(5+\frac{-4(b^2+5)^2+10(b^2+5)}{(5+b^2)^3})c=-b(1+2\frac{b^2+5}{(5+b^2)^3})e'^3_3$$

$$(5+\frac{-4(b^2+5)+10}{(5+b^2)^2})\frac{c}{b}=-(1+\frac{2}{(5+b^2)^2})e'^3_3$$

$$(5(5+b^2)^2-4(b^2+5)+10)\frac{c}{b}=-((5+b^2)^2+2)e'^3_3 \tag{3.4.45}$$

Из равенства (3.4.45) следует

$$e'^3_3=-\frac{c}{b}\frac{5(5+b^2)^2-4(b^2+5)+10}{(5+b^2)^2+2}=-\frac{c}{b}(5-\frac{4(b^2+5)}{(5+b^2)^2+2}) \tag{3.4.46}$$

Из равенств (3.4.42), (3.4.46) следует, что $c\ne 0$ произвольно и мы можем положить $c=1$

$$\overline{e}'_3=\begin{pmatrix}-1\\ 2\\ -\frac{1}{b}\left(5-\frac{4(b^2+5)}{(5+b^2)^2+2}\right)\end{pmatrix} \tag{3.4.47}$$

---

[3.7] Если $b=0$, то

$$e'_2=\begin{pmatrix}-1\\ 2\\ 0\end{pmatrix}$$

Из равенств (3.4.43) следует

$$(1-\frac{88}{5^3})c-\frac{23}{2*5^3}c-2*\frac{-23}{2*5^3}c+4*(1-\frac{2}{5^3})c=0$$

$$\left(5-\frac{96}{125}\right)c=0 \tag{3.4.44}$$

Из равенства (3.4.44) следует, что $c=0$. Из равенства (3.4.42) следует, что $e'^1_3=e'^2_3=0$. Следовательно, $e'^3_3\ne 0$ произвольно, и не нарушая общности, мы положим

$$\overline{e}'_3=\begin{pmatrix}0\\ 0\\ 1\end{pmatrix}$$

В подразделе 3.4.1, мы показали, что отношение ортогональности векторов $\overline{e}'_2$ и $\overline{e}'_3$ симметрично. Следовательно, для нас не имеет значения порядок векторов $\overline{e}'_2$, $\overline{e}'_3$ в базисе. Поскольку базис $\overline{\overline{e}}'$ рассмотрен в подразделе 3.4.1, мы не рассматриваем в этом подразделе случай $b=0$.

Пусть

$$d=-\frac{1}{b}\left(5-\frac{4(b^2+5)}{(5+b^2)^2+2}\right) \tag{3.4.48}$$

Тогда

$$\overline{e}'_3=\begin{pmatrix}-1\\2\\d\end{pmatrix} \tag{3.4.49}$$

Векторы $\overline{e}'_2$, $\overline{e}'_3$ имеют похожие координаты.[3.8] Согласно равенствам (3.4.36), (3.4.38), (3.4.49),

$$g(\overline{e}'_3)=\begin{pmatrix}1-\frac{88+46d^2+6d^4}{(5+d^2)^3} & \frac{-23+4d^2+3d^4}{2(5+d^2)^3} & \frac{-22d-6d^3}{(5+d^2)^3}\\ & 1+\frac{-2+6d^2}{(5+d^2)^3} & \frac{-11d+d^3}{(5+d^2)^3}\\ & & 1+\frac{-6d^2+10}{(5+d^2)^3}\end{pmatrix} \tag{3.4.53}$$

Так как векторы $\overline{e}'_2$, $\overline{e}'_3$ имеют похожие координаты, то вектор $\overline{e}'_1$ не ортогонален вектору $\overline{e}'_3$.

Поскольку

$$F^2(\overline{e}'_1)=1$$

то длина вектора

$$\overline{e}''_1=\overline{e}'_1=\begin{pmatrix}1\\0\\0\end{pmatrix}$$

равна 1. Поскольку

$$F^2(\overline{e}'_2)=(-1)^2+2^2+b^2+\frac{(-1)^3 2}{(-1)^2+2^2+b^2}=5+b^2-\frac{2}{5+b^2}$$

[3.8] Условие

$$d=b \tag{3.4.50}$$

означает равенство векторов $\overline{e}'_3$, $\overline{e}'_3$. Однако из равенств (3.4.48), (3.4.50) следует

$$b^2+5=\frac{4(b^2+5)}{(5+b^2)^2+2} \tag{3.4.51}$$

Поскольку

$$\frac{4}{(5+b^2)^2+2}<1 \tag{3.4.52}$$

то равенство (3.4.50) невозможно.

то длина вектора

$$\overline{e}''_2=\sqrt{\frac{5+b^2}{(5+b^2)^2-2}}\overline{e}'_2=\begin{pmatrix}-\sqrt{\dfrac{5+b^2}{(5+b^2)^2-2}}\\2\sqrt{\dfrac{5+b^2}{(5+b^2)^2-2}}\\b\sqrt{\dfrac{5+b^2}{(5+b^2)^2-2}}\end{pmatrix}$$

равна 1. Поскольку

$$F^2(\overline{e}'_3)=(-1)^2+2^2+d^2+\frac{(-1)^3 2}{(-1)^2+2^2+d^2}=5+d^2-\frac{2}{5+d^2}$$

то длина вектора

$$\overline{e}''_3=\sqrt{\frac{5+d^2}{(5+d^2)^2-2}}\overline{e}'_3=\begin{pmatrix}-\sqrt{\dfrac{5+d^2}{(5+d^2)^2-2}}\\2\sqrt{\dfrac{5+d^2}{(5+d^2)^2-2}}\\d\sqrt{\dfrac{5+d^2}{(5+d^2)^2-2}}\end{pmatrix}$$

равна 1. Следовательно, базис $\overline{\overline{e}}''$ является ортонормированным базисом.

Пусть

$$f_b=\sqrt{\frac{5+b^2}{(5+b^2)^2-2}}\tag{3.4.54}$$

$$f_d=\sqrt{\frac{5+d^2}{(5+d^2)^2-2}}\tag{3.4.55}$$

Базисы $\overline{\overline{e}}''$ и $\overline{\overline{e}}$ связаны пассивным преобразованием

$$\begin{pmatrix}\overline{e}''_1&\overline{e}''_2&\overline{e}''_3\end{pmatrix}=\begin{pmatrix}\overline{e}_1&\overline{e}_2&\overline{e}_3\end{pmatrix}\begin{pmatrix}1&-f_b&-f_d\\0&2f_b&2f_d\\0&bf_b&df_d\end{pmatrix}\tag{3.4.56}$$

Теорема 3.4.5. *Рассмотрим пространство Минковского с финслеровой метрикой, определённой равенством (3.4.1). Пусть вектор $\overline{v}$ имеет координаты*

$$\overline{v}=\begin{pmatrix}-\alpha\\2\alpha\\\gamma\end{pmatrix}\tag{3.4.57}$$

*Координаты* $g(\overline{v})$ *имеют вид*

(3.4.58)
$$\left\{\begin{aligned}
g_{11}(\overline{v})&=1-\frac{88\alpha^6+46\alpha^4\gamma^2+6\alpha^2\gamma^4}{(5\alpha^2+\gamma^2)^3}\\
g_{12}(\overline{v})&=\frac{-23\alpha^6+4\alpha^4\gamma^2+3\alpha^2\gamma^4}{2(5\alpha^2+\gamma^2)^3}\\
g_{13}(\overline{v})&=\frac{-22\alpha^5\gamma-6\alpha^3\gamma^3}{(5\alpha^2+\gamma^2)^3}\\
g_{22}(\overline{v})&=1+\frac{-2\alpha^6+6\alpha^4\gamma^2}{(5\alpha^2+\gamma^2)^3}\\
g_{23}(\overline{v})&=\frac{-11\alpha^5\gamma+\alpha^3\gamma^3}{(5\alpha^2+\gamma^2)^3}\\
g_{33}(\overline{v})&=1+\frac{-6\alpha^4\gamma^2+10\alpha^6}{(5\alpha^2+\gamma^2)^3}
\end{aligned}\right.$$

Доказательство. Согласно теореме 3.4.1 и равенству (3.4.57),

(3.4.59)
$$\left\{\begin{aligned}
g_{11}(\overline{v})&=1+\frac{6(-\alpha)(2\alpha)^3\gamma^2+3(-\alpha)2\alpha\gamma^4}{(5\alpha^2+\gamma^2)^3}\\
&\quad+\frac{3(-\alpha)(2\alpha)^5-(-\alpha)^3(2\alpha)^3-(-\alpha)^3 2\alpha\gamma^2}{(5\alpha^2+\gamma^2)^3}\\
&=1+\frac{-48\alpha^4\gamma^2-6\alpha^2\gamma^4-96\alpha^6+8\alpha^6+2\alpha^4\gamma^2}{(5\alpha^2+\gamma^2)^3}\\
g_{12}(\overline{v})&=\frac{3(-\alpha)^4(2\alpha)^2+2(-\alpha)^4\gamma^2}{(5\alpha^2+\gamma^2)^3}\\
&\quad+\frac{1}{2}\frac{\alpha^6-3(-\alpha)^2(2\alpha)^4+3(-\alpha)^2\gamma^4}{(5\alpha^2+\gamma^2)^3}\\
g_{13}(\overline{v})&=\frac{(-\alpha)^4 2\alpha\gamma-3(-\alpha)^2(2\alpha)^3\gamma-3(-\alpha)^2 2\alpha\gamma^3}{(5\alpha^2+\gamma^2)^3}\\
&=\frac{2\alpha^5\gamma-24\alpha^5\gamma-6\alpha^3\gamma^3}{(5\alpha^2+\gamma^2)^3}\\
g_{22}(\overline{v})&=1+\frac{(-\alpha)^3(2\alpha)^3-3(-\alpha)^5 2\alpha-3(-\alpha)^3 2\alpha\gamma^2}{(5\alpha^2+\gamma^2)^3}\\
&=1+\frac{-8\alpha^6+6\alpha^6+6\alpha^4\gamma^2}{(5\alpha^2+\gamma^2)^3}\\
g_{23}(\overline{v})&=\frac{3(-\alpha)^3(2\alpha)^2\gamma-(-\alpha)^5\gamma-(-\alpha)^3\gamma^3}{(5\alpha^2+\gamma^2)^3}\\
&=\frac{-12\alpha^5\gamma+\alpha^5\gamma+\alpha^3\gamma^3}{(5\alpha^2+\gamma^2)^3}\\
g_{33}(\overline{v})&=1+\frac{3(-\alpha)^3 2\alpha\gamma^2-(-\alpha)^5 2\alpha-(-\alpha)^3(2\alpha)^3}{(5\alpha^2+\gamma^2)^3}
\end{aligned}\right.$$

Равенства (3.4.58) следуют из равенств (3.4.59). □

Теорема 3.4.6. *Компоненты метрического тензора относительно базиса* $\overline{\overline{e}}''$ *имеют вид*

(3.4.60)
$$g''(\overline{e}''_1)=\begin{pmatrix}1&0&0\\0&f_b^2(3+b^2)&f_bf_d(3+bd)\\0&f_bf_d(3+bd)&f_d^2(3+d^2)\end{pmatrix}$$

$$
g''(\overline{e}''_2) =
$$

$$
= \begin{pmatrix}
1-\dfrac{88+46b^2+6b^4}{(5+b^2)^3} & -\dfrac{(b^2+3)(b^2+4)}{(5+b^2)^2}f_b & \boxed{\left(\dfrac{-b^4-b^2+10}{(5+b^2)^2} - \dfrac{88+24b^2}{(5+b^2)^2((5+b^2)^2+2)}\right)f_d} \\
-\dfrac{(b^2+3)(b^2+4)}{(5+b^2)^2}f_b & 1 & 0 \\
\boxed{\left(\dfrac{-b^4-b^2+10}{(5+b^2)^2} - \dfrac{88+24b^2}{(5+b^2)^2((5+b^2)^2+2)}\right)f_d} & 0 & \boxed{\left(5+\left(1+\dfrac{-6b^2+10}{(5+b^2)^3}\right)d^2 - \dfrac{50+110b^2+12b^4}{(5+b^2)^3} + \dfrac{64b^2}{(5+b^2)^2\,((5+b^2)^2+2)}\right)f_d^2}
\end{pmatrix} \tag{3.4.61}
$$

$$
g''(\overline{e}''_3) =
$$

$$
= \begin{pmatrix}
1-\dfrac{88+46d^2+6d^4}{(5+d^2)^3} & \boxed{f_b\left(-1+\dfrac{9d^2+35}{(5+d^2)^2} - \dfrac{8(11+3d^2)(5+b^2)}{(5+d^2)^3((5+b^2)^2+2)}\right)} & -\dfrac{(d^2+3)(d^2+4)}{(5+d^2)^2}f_d \\
\boxed{f_b\left(-1+\dfrac{9d^2+35}{(5+d^2)^2} - \dfrac{8(11+3d^2)(5+b^2)}{(5+d^2)^3((5+b^2)^2+2)}\right)} & \boxed{f_b^2\left(5+\left(1+\dfrac{-6d^2+10}{(5+d^2)^3}\right)b^2 - \dfrac{50+110d^2+12d^4}{(5+d^2)^3} + \dfrac{64d^2(5+b^2)}{(5+d^2)^3((5+b^2)^2+2)}\right)} & \boxed{f_bf_d\left(-\dfrac{4}{(5+d^2)^2} + \frac{4(5+b^2)}{(5+b^2)^2+2} * \left(1+\dfrac{2}{(5+d^2)^2}\right)\right)} \\
-\dfrac{(d^2+3)(d^2+4)}{(5+d^2)^2}f_d & \boxed{f_bf_d\left(-\dfrac{4}{(5+d^2)^2} + \frac{4(5+b^2)}{(5+b^2)^2+2} * \left(1+\dfrac{2}{(5+d^2)^2}\right)\right)} & 1
\end{pmatrix} \tag{3.4.62}
$$

Доказательство. Из равенств (2.5.3), (3.4.18), (3.4.56), следует

$$
\begin{cases}
g''_{11}(\overline{e}''_1) = g_{11}(\overline{e}''_1) = 1\\
g''_{12}(\overline{e}''_1) = g_{11}(\overline{e}''_1)A^1_2 + g_{12}(\overline{e}''_1)A^2_2 = 0\\
g''_{13}(\overline{e}''_1) = g_{11}(\overline{e}''_1)A^1_3 + g_{12}(\overline{e}''_1)A^2_3 = 0\\
g''_{22}(\overline{e}''_1) = g_{11}(\overline{e}''_1)A^1_2A^1_2 + 2g_{12}(\overline{e}''_1)A^1_2A^2_2\\
\qquad +g_{22}(\overline{e}''_1)A^2_2A^2_2 + g_{33}(\overline{e}''_1)A^3_2A^3_2\\
\qquad = f_b^2 - 2f_b^2 + 4f_b^2 + b^2f_b^2 = f_b^2(3+b^2)\\
g''_{23}(\overline{e}''_1) = g_{11}(\overline{e}''_1)A^1_2A^1_3 + g_{12}(\overline{e}''_1)A^1_2A^2_3 + g_{12}(\overline{e}''_1)A^1_3A^2_2\\
\qquad +g_{22}(\overline{e}''_1)A^2_2A^2_3 + g_{33}(\overline{e}''_1)A^3_2A^3_3\\
\qquad = (-f_b)(-f_d) + \frac{1}{2}(-f_b)2f_d + \frac{1}{2}2f_b(-f_d) + 4f_bf_d + bf_bdf_d\\
\qquad = f_bf_d(3+bd)\\
g''_{33}(\overline{e}''_1) = g_{11}(\overline{e}''_1)A^1_3A^1_3 + 2g_{12}(\overline{e}''_1)A^1_3A^2_3\\
\qquad +g_{22}(\overline{e}''_1)A^2_3A^2_3 + g_{33}(\overline{e}''_1)A^3_3A^3_3\\
\qquad = f_d^2 - 2f_d^2 + 4f_d^2 + d^2f_d^2 = f_d^2(3+d^2)
\end{cases}
\tag{3.4.63}
$$

Равенство (3.4.60) следует из (3.4.63). Из равенств (2.5.3), (3.4.38), (3.4.56), следует

$$
(3.4.64)\quad
\left\{
\begin{aligned}
g''_{11}(\overline{e}''_2) &= g_{11}(\overline{e}''_2) = 1-\frac{88+46b^2+6b^4}{(5+b^2)^3}\\
g''_{12}(\overline{e}''_2) &= g_{11}(\overline{e}''_2)A^1_2+g_{12}(\overline{e}''_2)A^2_2+g_{13}(\overline{e}''_2)A^3_2\\
&=\left(1-\frac{88+46b^2+6b^4}{(5+b^2)^3}\right)(-f_b)+\frac{-23+4b^2+3b^4}{2(5+b^2)^3}2f_b+\frac{-22b-6b^3}{(5+b^2)^3}bf_b\\
&=\left(-1+\frac{65+28b^2+3b^4}{(5+b^2)^3}\right)f_b=\left(-1+\frac{3b^2+13}{(5+b^2)^2}\right)f_b\\
&=\left(\frac{-b^4-7b^2-12}{(5+b^2)^2}\right)f_b=-\frac{(b^2+3)(b^2+4)}{(5+b^2)^2}f_b\\
g''_{13}(\overline{e}''_2) &= g_{11}(\overline{e}''_2)A^1_3+g_{12}(\overline{e}''_2)A^2_3+g_{13}(\overline{e}''_2)A^3_3\\
&=\left(1-\frac{88+46b^2+6b^4}{(5+b^2)^3}\right)(-f_d)+\frac{-23+4b^2+3b^4}{2(5+b^2)^3}2f_d+\frac{-22b-6b^3}{(5+b^2)^3}df_d\\
&=\left(-1+\frac{65+50b^2+9b^4}{(5+b^2)^3}-\frac{22+6b^2}{(5+b^2)^3}bd\right)f_d\\
&=\left(-1+\frac{65+50b^2+9b^4}{(5+b^2)^3}-\frac{22+6b^2}{(5+b^2)^3}\left(\frac{4(5+b^2)}{(5+b^2)^2+2}-5\right)\right)f_d\\
&=\left(-1+\frac{175+80b^2+9b^4}{(5+b^2)^3}-\frac{22+6b^2}{(5+b^2)^3}\frac{4(5+b^2)}{(5+b^2)^2+2}\right)f_d\\
&=\left(-1+\frac{9b^2+35}{(5+b^2)^2}-\frac{22+6b^2}{(5+b^2)^2}\frac{4}{(5+b^2)^2+2}\right)f_d\\
&=\left(\frac{-b^4-b^2+10}{(5+b^2)^2}-\frac{88+24b^2}{(5+b^2)^2((5+b^2)^2+2)}\right)f_d\\
g''_{22}(\overline{e}''_2) &= g_{11}(\overline{e}''_2)A^1_2A^1_2+2g_{12}(\overline{e}''_2)A^1_2A^2_2+2g_{13}(\overline{e}''_2)A^1_2A^3_2\\
&\quad+g_{22}(\overline{e}''_2)A^2_2A^2_2+2g_{23}(\overline{e}''_2)A^2_2A^3_2+g_{33}(\overline{e}''_2)A^3_2A^3_2\\
&=\left(1-\frac{88+46b^2+6b^4}{(5+b^2)^3}\right)(-f_b)(-f_b)\\
&\quad+2\frac{-23+4b^2+3b^4}{2(5+b^2)^3}(-f_b)2f_b\\
&\quad+2\frac{-22b-6b^3}{(5+b^2)^3}(-f_b)bf_b+\left(1+\frac{-2+6b^2}{(5+b^2)^3}\right)2f_b2f_b\\
&\quad+2\frac{-11b+b^3}{(5+b^2)^3}2f_bbf_b+\left(1+\frac{-6b^2+10}{(5+b^2)^3}\right)bf_bbf_b\\
&=(5+b^2-2\frac{25+10b^2+b^4}{(5+b^2)^3})f_b^2=1
\end{aligned}
\right.
$$

$$
\begin{aligned}
g''_{23}(\overline{e}''_2) =& g_{11}(\overline{e}''_2)A_2^1A_3^1+g_{12}(\overline{e}''_2)A_2^1A_3^2+g_{12}(\overline{e}''_2)A_3^1A_2^2\\
&+g_{13}(\overline{e}''_2)A_2^1A_3^3+g_{13}(\overline{e}''_2)A_3^1A_2^3\\
&+g_{22}(\overline{e}''_2)A_2^2A_3^2+g_{23}(\overline{e}''_2)A_2^2A_3^3+g_{23}(\overline{e}''_2)A_3^2A_2^3+g_{33}(\overline{e}''_2)A_2^3A_3^3\\
=&\left(1-\frac{88+46b^2+6b^4}{(5+b^2)^3}\right)(-f_b)(-f_d)\\
&+\frac{-23+4b^2+3b^4}{2(5+b^2)^3}(-f_b)2f_d+\frac{-23+4b^2+3b^4}{2(5+b^2)^3}(-f_d)2f_b\\
&+\frac{-22b-6b^3}{(5+b^2)^3}(-f_b)df_d+\frac{-22b-6b^3}{(5+b^2)^3}(-f_d)bf_b\\
&+\left(1+\frac{-2+6b^2}{(5+b^2)^3}\right)2f_b2f_d+\frac{-11b+b^3}{(5+b^2)^3}2f_bdf_d\\
&+\frac{-11b+b^3}{(5+b^2)^3}2f_dbf_b+\left(1+\frac{-6b^2+10}{(5+b^2)^3}\right)bf_ddf_d\\
=&\left(5+bd-\frac{50+30b^2+4b^4}{(5+b^2)^3}\right)f_bf_d+\frac{10+2b^2}{(5+b^2)^3}bdf_bf_d\\
=&\left(\frac{4(5+b^2)}{(5+b^2)^2+2}-\frac{50+30b^2+4b^4}{(5+b^2)^3}\right.\\
&\left.+\frac{2}{(5+b^2)^2}\left(\frac{4(5+b^2)}{(5+b^2)^2+2}-5\right)\right)f_bf_d\\
=&\left(\frac{4(5+b^2)}{(5+b^2)^2+2}-\frac{50+30b^2+4b^4}{(5+b^2)^3}-\frac{10}{(5+b^2)^2}\right.\\
&\left.+\frac{8}{(5+b^2)\left((5+b^2)^2+2\right)}\right)f_bf_d\\
=&\left(\frac{4(5+b^2)}{(5+b^2)^2+2}-\frac{100+40b^2+4b^4}{(5+b^2)^3}\right.\\
&\left.+\frac{8}{(5+b^2)\left((5+b^2)^2+2\right)}\right)f_bf_d\\
=&\left(\frac{4(5+b^2)}{(5+b^2)^2+2}-\frac{4}{5+b^2}+\frac{8}{(5+b^2)\left((5+b^2)^2+2\right)}\right)f_bf_d=0
\end{aligned}
\tag{3.4.65}
$$

$$
\begin{aligned}
g''_{33}(\overline{e}''_2) &= g_{11}(\overline{e}''_2)A_3^1A_3^1 + 2g_{12}(\overline{e}''_2)A_3^1A_3^2 + 2g_{13}(\overline{e}''_2)A_3^1A_3^3 \\
&+g_{22}(\overline{e}''_2)A_3^2A_3^2 + 2g_{23}(\overline{e}''_2)A_3^2A_3^3 + g_{33}(\overline{e}''_2)A_3^3A_3^3 \\
&= \left(1 - \frac{88+46b^2+6b^4}{(5+b^2)^3}\right)(-f_d)(-f_d) \\
&+2\frac{-23+4b^2+3b^4}{2(5+b^2)^3}(-f_d)2f_d \\
&+2\frac{-22b-6b^3}{(5+b^2)^3}(-f_d)df_d + \left(1+\frac{-2+6b^2}{(5+b^2)^3}\right)2f_d2f_d \\
&+2\frac{-11b+b^3}{(5+b^2)^3}2f_ddf_d + \left(1+\frac{-6b^2+10}{(5+b^2)^3}\right)df_ddf_d \\
&= \left(5+\left(1+\frac{-6b^2+10}{(5+b^2)^3}\right)d^2 - \frac{50+30b^2+12b^4}{(5+b^2)^3}\right. \\
&+\left.\frac{16b^3d}{(5+b^2)^3}\right)f_d^2 \\
&= \left(5+\left(1+\frac{-6b^2+10}{(5+b^2)^3}\right)d^2 - \frac{50+30b^2+12b^4}{(5+b^2)^3}\right. \\
&+\left.\frac{16b^2}{(5+b^2)^3}\left(\frac{4(5+b^2)}{(5+b^2)^2+2}-5\right)\right)f_d^2 \\
&= \left(5+\left(1+\frac{-6b^2+10}{(5+b^2)^3}\right)d^2 - \frac{50+30b^2+12b^4}{(5+b^2)^3} - \frac{80b^2}{(5+b^2)^3}\right. \\
&+\left.\frac{64b^2}{(5+b^2)^2\left((5+b^2)^2+2\right)}\right)f_d^2 \\
&= \left(5+\left(1+\frac{-6b^2+10}{(5+b^2)^3}\right)d^2 - \frac{50+110b^2+12b^4}{(5+b^2)^3} + \frac{64b^2}{(5+b^2)^2\left((5+b^2)^2+2\right)}\right)f_d^2
\end{aligned}
\tag{3.4.66}
$$

Равенство (3.4.61) следует из (3.4.64), (3.4.65), (3.4.66). Из равенств (2.5.3), (3.4.53), (3.4.56), следует

$$
(3.4.67)\quad \left\{
\begin{aligned}
g''_{11}(\overline{e}''_3)&=g_{11}(\overline{e}''_3)=1-\frac{88+46d^2+6d^4}{(5+d^2)^3}\\
g''_{12}(\overline{e}''_3)&=g_{11}(\overline{e}''_3)A^1_2+g_{12}(\overline{e}''_3)A^2_2+g_{13}(\overline{e}''_3)A^3_2\\
&=\left(1-\frac{88+46d^2+6d^4}{(5+d^2)^3}\right)(-f_b)+\frac{-23+4d^2+3d^4}{2(5+d^2)^3}2f_b+\frac{-22d-6d^3}{(5+d^2)^3}bf_b\\
&=\left(-1+\frac{65+50d^2+9d^4}{(5+d^2)^3}-\frac{22+6d^2}{(5+d^2)^3}bd\right)f_b\\
&=\left(-1+\frac{65+50d^2+9d^4}{(5+d^2)^3}-\frac{22+6d^2}{(5+d^2)^3}\left(\frac{4(5+b^2)}{(5+b^2)^2+2}-5\right)\right)f_b\\
&=\left(-1+\frac{175+80d^2+9d^4}{(5+d^2)^3}-\frac{8(11+3d^2)(5+b^2)}{(5+d^2)^3((5+b^2)^2+2)}\right)f_b\\
&=\left(-1+\frac{9d^2+35}{(5+d^2)^2}-\frac{8(11+3d^2)(5+b^2)}{(5+d^2)^3((5+b^2)^2+2)}\right)f_b\\
g''_{13}(\overline{e}''_3)&=g_{11}(\overline{e}''_3)A^1_3+g_{12}(\overline{e}''_3)A^2_3+g_{13}(\overline{e}''_3)A^3_3\\
&=\left(1-\frac{88+46d^2+6d^4}{(5+d^2)^3}\right)(-f_d)+\frac{-23+4d^2+3d^4}{2(5+d^2)^3}2f_d\\
&+\frac{-22d-6d^3}{(5+d^2)^3}df_d\\
&=\left(-1+\frac{65+28d^2+3d^4}{(5+d^2)^3}\right)f_d=\left(-1+\frac{3d^2+13}{(5+d^2)^2}\right)f_d\\
&=\left(\frac{-d^4-7d^2-12}{(5+d^2)^2}\right)f_d=-\frac{(d^2+3)(d^2+4)}{(5+d^2)^2}f_d\\
g''_{22}(\overline{e}''_3)&=g_{11}(\overline{e}''_3)A^1_2A^1_2+2g_{12}(\overline{e}''_3)A^1_2A^2_2+2g_{13}(\overline{e}''_3)A^1_2A^3_2\\
&+g_{22}(\overline{e}''_3)A^2_2A^2_2+2g_{23}(\overline{e}''_3)A^2_2A^3_2+g_{33}(\overline{e}''_3)A^3_2A^3_2\\
&=\left(1-\frac{88+46d^2+6d^4}{(5+d^2)^3}\right)(-f_b)(-f_b)\\
&+2\frac{-23+4d^2+3d^4}{2(5+d^2)^3}(-f_b)2f_b\\
&+2\frac{-22d-6d^3}{(5+d^2)^3}(-f_b)bf_b+\left(1+\frac{-2+6d^2}{(5+d^2)^3}\right)2f_b2f_b\\
&+2\frac{-11d+d^3}{(5+d^2)^3}2f_bbf_b+\left(1+\frac{-6d^2+10}{(5+d^2)^3}\right)bf_bbf_b\\
&=\left(5+\left(1+\frac{-6d^2+10}{(5+d^2)^3}\right)b^2-\frac{50+30d^2+12d^4}{(5+d^2)^3}+\frac{16d^2}{(5+d^2)^3}bd\right)f_b^2\\
&=\left(5+\left(1+\frac{-6d^2+10}{(5+d^2)^3}\right)b^2-\frac{50+30d^2+12d^4}{(5+d^2)^3}+\frac{16d^2}{(5+d^2)^3}\left(\frac{4(5+b^2)}{(5+b^2)^2+2}-5\right)\right)f_b^2\\
&=\left(5+\left(1+\frac{-6d^2+10}{(5+d^2)^3}\right)b^2-\frac{50+30d^2+12d^4}{(5+d^2)^3}-\frac{80d^2}{(5+d^2)^3}\right.\\
&\left.+\frac{64d^2(5+b^2)}{(5+d^2)^3((5+b^2)^2+2)}\right)f_b^2\\
&=\left(5+\left(1+\frac{-6d^2+10}{(5+d^2)^3}\right)b^2-\frac{50+110d^2+12d^4}{(5+d^2)^3}+\frac{64d^2(5+b^2)}{(5+d^2)^3((5+b^2)^2+2)}\right)f_b^2
\end{aligned}
\right.
$$

$$
(3.4.68)\quad
\left\{
\begin{aligned}
g''_{23}(\overline{e}''_3)&=g_{11}(\overline{e}''_3)A_2^1A_3^1+g_{12}(\overline{e}''_3)A_2^1A_3^2+g_{12}(\overline{e}''_3)A_3^1A_2^2\\
&+g_{13}(\overline{e}''_3)A_2^1A_3^3+g_{13}(\overline{e}''_3)A_3^1A_2^3\\
&+g_{22}(\overline{e}''_3)A_2^2A_3^2+g_{23}(\overline{e}''_3)A_2^2A_3^3+g_{23}(\overline{e}''_3)A_3^2A_2^3+g_{33}(\overline{e}''_3)A_3^3A_3^3\\
&=\left(1-\frac{88+46d^2+6d^4}{(5+d^2)^3}\right)(-f_b)(-f_d)\\
&+\frac{-23+4d^2+3d^4}{2(5+d^2)^3}(-f_b)2f_d+\frac{-23+4d^2+3d^4}{2(5+d^2)^3}(-f_d)2f_b\\
&+\frac{-22d-6d^3}{(5+d^2)^3}(-f_b)df_d+\frac{-22d-6d^3}{(5+d^2)^3}(-f_d)bf_b\\
&+\left(1+\frac{-2+6d^2}{(5+d^2)^3}\right)2f_b2f_d+\frac{-11d+d^3}{(5+d^2)^3}2f_bdf_d\\
&+\frac{-11d+d^3}{(5+d^2)^3}2f_dbf_b+\left(1+\frac{-6d^2+10}{(5+d^2)^3}\right)df_ddf_d\\
&=\left(\left(5+bd-\frac{50+30d^2+4d^4}{(5+d^2)^3}\right)+\frac{10+2d^2}{(5+d^2)^3}bd\right)f_bf_d\\
&=\left(-\frac{50+30d^2+4d^4}{(5+d^2)^3}+\frac{4(5+b^2)}{(5+b^2)^2+2}+\frac{2(5+d^2)}{(5+d^2)^3}\left(\frac{4(5+b^2)}{(5+b^2)^2+2}-5\right)\right)f_bf_d\\
&=\left(-\frac{50+30d^2+4d^4}{(5+d^2)^3}+\frac{4\,(5+b^2)}{(5+b^2)^2+2}+\frac{2}{(5+d^2)^2}\left(\frac{4\,(5+b^2)}{(5+b^2)^2+2}-5\right)\right)f_bf_d\\
&=\left(-\frac{50+30d^2+4d^4}{(5+d^2)^3}+\frac{4(5+b^2)}{(5+b^2)^2+2}-\frac{10}{(5+d^2)^2}+\frac{8(5+b^2)}{(5+d^2)^2((5+b^2)^2+2)}\right)f_bf_d\\
&=\left(-\frac{50+30d^2+4d^4}{(5+d^2)^3}-\frac{10(5+d^2)}{(5+d^2)^3}+\frac{4(5+b^2)}{(5+b^2)^2+2}\left(1+\frac{2}{(5+d^2)^2}\right)\right)f_bf_d\\
&=\left(-\frac{100+40d^2+4d^4}{(5+d^2)^3}+\frac{4(5+b^2)}{(5+b^2)^2+2}\left(1+\frac{2}{(5+d^2)^2}\right)\right)f_bf_d\\
&=\left(-\frac{4(5+d^2)^2}{(5+d^2)^3}+\frac{4(5+b^2)}{(5+b^2)^2+2}\left(1+\frac{2}{(5+d^2)^2}\right)\right)f_bf_d\\
&=\left(-\frac{4}{(5+d^2)^2}+\frac{4(5+b^2)}{(5+b^2)^2+2}\left(1+\frac{2}{(5+d^2)^2}\right)\right)f_bf_d\\
g''_{33}(\overline{e}''_3)&=g_{11}(\overline{e}''_3)A_3^1A_3^1+2g_{12}(\overline{e}''_3)A_3^1A_3^2+2g_{13}(\overline{e}''_3)A_3^1A_3^3\\
&+g_{22}(\overline{e}''_3)A_3^2A_3^2+2g_{23}(\overline{e}''_3)A_3^2A_3^3+g_{33}(\overline{e}''_3)A_3^3A_3^3\\
&=\left(1-\frac{88+46d^2+6d^4}{(5+d^2)^3}\right)(-f_d)(-f_d)\\
&+2\frac{-23+4d^2+3d^4}{2(5+d^2)^3}(-f_d)2f_d\\
&+2\frac{-22d-6d^3}{(5+d^2)^3}(-f_d)df_d+\left(1+\frac{-2+6d^2}{(5+d^2)^3}\right)2f_d2f_d\\
&+2\frac{-11d+d^3}{(5+d^2)^3}2f_ddf_d+\left(1+\frac{-6d^2+10}{(5+d^2)^3}\right)df_ddf_d\\
&=\left(5+d^2-\frac{50+20d^2+2d^4}{(5+d^2)^3}\right)f_d^2=\left(5+d^2-\frac{2(5+d^2)^2}{(5+d^2)^3}\right)f_d^2\\
&=\left(5+d^2-\frac{2}{5+d^2}\right)f_d^2=\frac{1}{f_d^2}f_d^2=1
\end{aligned}
\right.
$$

Равенство (3.4.62) следует из (3.4.67), (3.4.68). □

Следствие 3.4.7. *Относительно базиса $\overline{\overline{e}}''$, определённого равенством (3.4.56), матрица (2.3.9), соответствующая финслеровой метрике $F^2$, определённой равенством (3.4.1), имеет вид*

$$\begin{pmatrix} 1 & 0 & 0 \\ -\frac{(b^2+3)(b^2+4)}{(5+b^2)^2}f_b & 1 & 0 \\ -\frac{(d^2+3)(d^2+4)}{(5+d^2)^2}f_d & \boxed{\begin{matrix} f_bf_d\left(-\frac{4}{(5+d^2)^2}\right. \\ +\frac{4(5+b^2)}{(5+b^2)^2+2} \\ *\left.\left(1+\frac{2}{(5+d^2)^2}\right)\right) \end{matrix}} & 1 \end{pmatrix} \tag{3.4.69}$$

□

### 3.4.3. Пример движения в пространстве Минковского размерности 3.

Теорема 3.4.8. *Рассмотрим пространство Минковского с финслеровой метрикой, определённой равенством (3.4.1).*

*Согласно построению в подразделе 3.4.1, базис $\overline{\overline{e}}_1$, определённый равенством (3.4.28),*

$$\begin{pmatrix}\overline{e}_{1\cdot 1} & \overline{e}_{1\cdot 2} & \overline{e}_{1\cdot 3}\end{pmatrix} = \begin{pmatrix}\overline{e}_1 & \overline{e}_2 & \overline{e}_3\end{pmatrix}\begin{pmatrix} 1 & 0 & -\sqrt{\frac{5}{23}} \\ 0 & 0 & 2\sqrt{\frac{5}{23}} \\ 0 & 1 & 0\end{pmatrix} \tag{3.4.70}$$

$$\overline{e}_{1\cdot k} = e_{1\cdot}{}^{i}_{k}\,\overline{e}_i \tag{3.4.71}$$

*является ортонормированным базисом.*

*Согласно построению в подразделе 3.4.2, базис $\overline{\overline{e}}_2$, определённый равенством (3.4.56),*

$$\begin{pmatrix}\overline{e}_{2\cdot 1} & \overline{e}_{2\cdot 2} & \overline{e}_{2\cdot 3}\end{pmatrix} = \begin{pmatrix}\overline{e}_1 & \overline{e}_2 & \overline{e}_3\end{pmatrix}\begin{pmatrix} 1 & -f_b & -f_d \\ 0 & 2f_b & 2f_d \\ 0 & bf_b & df_d\end{pmatrix} \tag{3.4.72}$$

$$\overline{e}_{2\cdot k} = e_{2\cdot}{}^{i}_{k}\,\overline{e}_i \tag{3.4.73}$$

*является ортонормированным базисом.*

*Рассмотрим движение $\overline{A}$, отображающее базис $\overline{\overline{e}}_1$ в базис $\overline{\overline{e}}_2$*

$$e_{2\cdot}{}^{i}_{k} = A^i_j e_{1\cdot}{}^{j}_{k} \tag{3.4.74}$$

*Пусть движение $\overline{A}$ отображает базис $\overline{\overline{e}}_2$ в базис $\overline{\overline{e}}_3$*

$$e_{3\cdot}{}^{i}_{k} = A^i_j e_{2\cdot}{}^{j}_{k} \tag{3.4.75}$$

*Базис $\overline{\overline{e}}_3$ не является ортонормированным базисом.*

Доказательство. Из равенства (3.4.74) следует

$$A^i_j = e_{2\cdot}{}^{i}_{k} e_1^{-1}{}_{\cdot j}^{k} \tag{3.4.76}$$

Из равенства (3.4.70) следует

$$e_1^{-1}=\begin{pmatrix}1&\frac{1}{2}&0\\0&0&1\\0&\frac{1}{2}\sqrt{\frac{23}{5}}&0\end{pmatrix}\tag{3.4.77}$$

Из равенств (3.4.72), (3.4.76), (3.4.77) следует

$$A=\begin{pmatrix}1&-f_b&-f_d\\0&2f_b&2f_d\\0&bf_b&df_d\end{pmatrix}\begin{pmatrix}1&\frac{1}{2}&0\\0&0&1\\0&\frac{1}{2}\sqrt{\frac{23}{5}}&0\end{pmatrix}=\begin{pmatrix}1&\frac{1}{2}-\frac{1}{2}f_d\sqrt{\frac{23}{5}}&-f_b\\0&f_d\sqrt{\frac{23}{5}}&2f_b\\0&\frac{1}{2}df_d\sqrt{\frac{23}{5}}&bf_b\end{pmatrix}\tag{3.4.78}$$

Из равенств (3.4.72), (3.4.78), (3.4.75) следует

$$\begin{aligned}\overline{\overline{e}}_3&=\begin{pmatrix}1&\frac{1}{2}-\frac{1}{2}f_d\sqrt{\frac{23}{5}}&-f_b\\0&f_d\sqrt{\frac{23}{5}}&2f_b\\0&\frac{1}{2}df_d\sqrt{\frac{23}{5}}&bf_b\end{pmatrix}\begin{pmatrix}1&-f_b&-f_d\\0&2f_b&2f_d\\0&bf_b&df_d\end{pmatrix}\\&=\left(\begin{pmatrix}1\\0\\0\end{pmatrix}\quad\begin{pmatrix}-f_bf_d\sqrt{\frac{23}{5}}-bf_b^2\\2f_bf_d\sqrt{\frac{23}{5}}+2bf_b^2\\df_bf_d\sqrt{\frac{23}{5}}+b^2f_b^2\end{pmatrix}\quad\begin{pmatrix}-f_d^2\sqrt{\frac{23}{5}}-df_bf_d\\2f_d^2\sqrt{\frac{23}{5}}+2df_bf_d\\df_d^2\sqrt{\frac{23}{5}}+bdf_bf_d\end{pmatrix}\right)\end{aligned}\tag{3.4.79}$$

Пусть

$$\begin{aligned}\alpha&=f_bf_d\sqrt{\frac{23}{5}}+bf_b^2&\gamma&=df_bf_d\sqrt{\frac{23}{5}}+b^2f_b^2\\\beta&=f_d^2\sqrt{\frac{23}{5}}+df_bf_d&\rho&=df_d^2\sqrt{\frac{23}{5}}+bdf_bf_d\end{aligned}\tag{3.4.80}$$

Тогда

$$\overline{e}_{3\cdot 2}=\begin{pmatrix}-\alpha\\2\alpha\\\gamma\end{pmatrix}\qquad\overline{e}_{3.3}=\begin{pmatrix}-\beta\\2\beta\\\rho\end{pmatrix}\tag{3.4.81}$$

Поскольку последующие вычисления очень сложны, то для того, чтобы получить ответ для конкретных значений $a$, мы написали программу на $C\#$. Расчёты выполнены на основании теоремы 3.4.5 и равенств (3.4.48), (3.4.54), (3.4.55), (3.4.80). Утверждение теоремы следует из расчётов, приведенных в примерах 3.4.9, 3.4.10. □

Пример 3.4.9. Пусть $b = 1$. Тогда

$$d = -4.3684210526 \quad f_b = 0.4200840252 \quad f_d = 0.2041239040$$

$$\overline{e}_{3\cdot 2} = \begin{pmatrix} -0.3603821145 \\ 0.7207642289 \\ -0.6269323948 \end{pmatrix} \quad \overline{e}_{3\cdot 3} = \begin{pmatrix} 0.2852237394 \\ -0.5704474788 \\ -0.7649717899 \end{pmatrix}$$

$$g(\overline{e}_{3\cdot 2}) = \begin{pmatrix} 0.4543053517 & 0.0160337367 & 0.1351074203 \\ 0.0160337367 & 1.031248994 & 0.0268266657 \\ 0.1351074203 & 0.0268266657 & 0.9842227918 \end{pmatrix}$$

$$g_{ij}(\overline{e}_{3\cdot 2}) e_{3\cdot 2}^{\ i} e_{3\cdot 3}^{\ j} = 0.0157973688$$

□

Пример 3.4.10. Пусть $b = 100$. Тогда

$$d = -0.049996002 \quad f_b = 0.009997501 \quad f_d = 0.4661156563$$

$$\overline{e}_{3\cdot 2} = \begin{pmatrix} -0.019989571554 \\ 0.0399791431 \\ 0.9990005814 \end{pmatrix} \quad \overline{e}_{3\cdot 3} = \begin{pmatrix} -0.4657459676 \\ 0.9314919353 \\ -0.0465951802 \end{pmatrix}$$

$$g(\overline{e}_{3.2}) = \begin{pmatrix} 0.9976047372 & 0.0005972998 & -0.0000785154 \\ 0.0005972998 & 1.000000956 & 0.0000079285 \\ -0.0000785154 & 0.0000079285 & 0.9999990446 \end{pmatrix}$$

$$g_{ij}(\overline{e}_{3\cdot 2}) e_{3\cdot 2}^{\ i} e_{3\cdot 3}^{\ j} = -0.000013181$$

□

## 3.5. Пример пространства Минковского размерности 4

Мы рассмотрим пример финслеровой метрики пространства Минковского размерности 4 (рассматриваемая метрика является частным случаем метрики [12]-(3); мы полагаем $c = 1$)

$$\begin{aligned} F^2(\overline{v}) &= -(v^0)^2 + ((v^1)^2 + (v^2)^2 + (v^3)^2)\left(1 + \frac{(v^1)^3 v^2}{((v^1)^2 + (v^2)^2 + (v^3)^2)^2}\right) \\ &= -(v^0)^2 + (v^1)^2 + (v^2)^2 + (v^3)^2 + \frac{(v^1)^3 v^2}{(v^1)^2 + (v^2)^2 + (v^3)^2} \end{aligned} \tag{3.5.1}$$

Соответствующее пространство Минковского $\overline{V}$ можно представить в виде прямой суммы

$$\overline{V} = \overline{V}_3 \oplus \overline{R}$$

где $\overline{V}_3$ - пространство Минковского, рассмотренное в разделе 3.4.

Пусть

$$\begin{aligned} F_r^2(v^0) &= -(v^0)^2 \\ F_3^2(v^1, v^2, v^3) &= (v^1)^2 + (v^2)^2 + (v^3)^2 + \frac{(v^1)^3 v^2}{(v^1)^2 + (v^2)^2 + (v^3)^2} \end{aligned} \tag{3.5.2}$$

Соответственно, мы можем представить функцию $F^2$ в виде суммы

(3.5.3)
$$F^2(v^0,v^1,v^2,v^3)=F_3^2(v^1,v^2,v^3)+F_r^2(v^0)$$

Функция $F_3^2$ является финслеровой метрикой пространства Минковского, рассмотренного в разделе 3.4.

Теорема 3.5.1. *Функция $F_3^2$ непрерывна в точке $(0,0,0)$ и может быть продолжена по непрерывности*[3.9]

(3.5.4)
$$F_3^2(0,0,0)=0$$

Доказательство. Так как
$$(v^1)^2\le(v^1)^2+(v^2)^2+(v^3)^2$$
то $((v^1,v^2,v^3)\ne(0,0,0))$

(3.5.5)
$$\frac{1}{(v^1)^2+(v^2)^2+(v^3)^2}\le\frac{1}{(v^1)^2}$$

Из неравенства (3.5.5) и равенства (3.5.2) следует $((v^1,v^2,v^3)\ne(0,0,0))$

(3.5.6)
$$|F_3^2(v^1,v^2,v^3)|\le(v^1)^2+(v^2)^2+(v^3)^2+\left|\frac{(v^1)^3v^2}{(v^1)^2}\right|=(v^1)^2+(v^2)^2+(v^3)^2+|v^1v^2|$$

Для заданного $\epsilon>0$ пусть $\delta=\dfrac{\sqrt{\epsilon}}{2}$. Тогда из неравенства (3.5.6) и требования
$$|v^1|<\delta\quad|v^2|<\delta\quad|v^3|<\delta$$
следует $((v^1,v^2,v^3)\ne(0,0,0))$

(3.5.7)
$$|F_3^2(v^1,v^2,v^3)|\le\frac{\epsilon}{4}+\frac{\epsilon}{4}+\frac{\epsilon}{4}+\frac{\epsilon}{4}=\epsilon$$

Следовательно, функция $F_3^2$ непрерывна в точке $(0,0,0)$ и
$$\lim_{(v^1,v^2,v^3)\to(0,0,0)}F_3^2(v^1,v^2,v^3)=0$$

□

Теорема 3.5.2.

(3.5.8)
$$F^2(a,0,0,0)=F_r^2(a)=-a^2$$

Доказательство. Равенство (3.5.8) следует из равенства (3.5.3) и теоремы 3.5.1. □

3.9 Значение функции $F_3^2$ в точке $(0,0,0)$ плохо определено, так как выражение включает отношение двух бесконечно малых.

Теорема 3.5.3. *Финслерова метрика $F^2$, определённая равенством (3.5.1) порождает метрический тензор пространства Минковского*

$$g_{00}(\overline{v}) = -1 \tag{3.5.9}$$

$$\begin{aligned} g_{11}(\overline{v}) = 1 &- \frac{(v^1)^3(v^2)^3+(v^1)^3v^2(v^3)^2}{((v^1)^2+(v^2)^2+(v^3)^2)^3} \\ &+ \frac{6v^1(v^2)^3(v^3)^2+3v^1(v^2)^5+3v^1v^2(v^3)^4}{((v^1)^2+(v^2)^2+(v^3)^2)^3} \end{aligned} \tag{3.5.10}$$

$$g_{12}(\overline{v}) = \frac{1}{2}\frac{(v^1)^6+6(v^1)^4(v^2)^2-3(v^1)^2(v^2)^4+4(v^1)^4(v^3)^2+3(v^1)^2(v^3)^4}{((v^1)^2+(v^2)^2+(v^3)^2)^3} \tag{3.5.11}$$

$$g_{13}(\overline{v}) = \frac{(v^1)^4v^2v^3-3(v^1)^2(v^2)^3v^3-3(v^1)^2v^2(v^3)^3}{((v^1)^2+(v^2)^2+(v^3)^2)^3} \tag{3.5.12}$$

$$g_{22}(\overline{v}) = 1+\frac{(v^1)^3(v^2)^3-3(v^1)^5v^2-3(v^1)^3v^2(v^3)^2}{((v^1)^2+(v^2)^2+(v^3)^2)^3} \tag{3.5.13}$$

$$g_{23}(\overline{v}) = \frac{3(v^1)^3(v^2)^2v^3-(v^1)^5v^3-(v^1)^3(v^3)^3}{((v^1)^2+(v^2)^2+(v^3)^2)^3} \tag{3.5.14}$$

$$g_{33}(\overline{v}) = 1+\frac{3(v^1)^3v^2(v^3)^2-(v^1)^5v^2-(v^1)^3(v^2)^3}{((v^1)^2+(v^2)^2+(v^3)^2)^3} \tag{3.5.15}$$

Доказательство. Равенства (3.5.10), (3.5.11), (3.5.12), (3.5.13), (3.5.14), (3.5.15) являются следствием теоремы 3.4.1 и равенства (3.5.3). Из равенства (3.5.3) следует

$$\frac{\partial F^2(\overline{v})}{\partial v^i} = 2F_3(v^1,v^2,v^3)\frac{\partial F_3(v^1,v^2,v^3)}{\partial v^i} \quad i=1,2,3 \tag{3.5.16}$$

Из равенств (2.2.2), (3.5.1), (3.5.16) следует

$$g_{0j}(\overline{v}) = 0 \quad i=1,2,3 \tag{3.5.17}$$

□

Теорема 3.5.4. *Если векторы*

$$\overline{e}_1 = \begin{pmatrix} 0 \\ e_1^1 \\ e_1^2 \\ e_1^3 \end{pmatrix} \quad \overline{e}_2 = \begin{pmatrix} 0 \\ e_2^1 \\ e_2^2 \\ e_2^3 \end{pmatrix} \quad \overline{e}_3 = \begin{pmatrix} 0 \\ e_3^1 \\ e_3^2 \\ e_3^3 \end{pmatrix} \tag{3.5.18}$$

*порождают ортонормированный базис пространства Минковского $\overline{V}_3$, то векторы*

$$\overline{e}_1 = \begin{pmatrix} 0 \\ e_1^1 \\ e_1^2 \\ e_1^3 \end{pmatrix} \quad \overline{e}_2 = \begin{pmatrix} 0 \\ e_2^1 \\ e_2^2 \\ e_2^3 \end{pmatrix} \quad \overline{e}_3 = \begin{pmatrix} 0 \\ e_3^1 \\ e_3^2 \\ e_3^3 \end{pmatrix} \quad \overline{e}_0 = \begin{pmatrix} 1 \\ 0 \\ 0 \\ 0 \end{pmatrix} \tag{3.5.19}$$

*порождают ортонормированный базис пространства Минковского $\overline{V}$.*

Доказательство. Согласно равенству (3.5.3), если векторы $\overline{a}, \overline{b} \in \overline{V}_3$ ортогональны в пространстве Минковского $\overline{V}_3$, то эти векторы также ортогональны в пространстве Минковского $\overline{V}$. Следовательно, векторы $\overline{e}_1$, $\overline{e}_2$, $\overline{e}_3$ принадлежат ортонормированному базису пространства Минковского $\overline{V}$. Из теоремы 3.5.3 следует, что

$$g_{kl}(\overline{e}_i)e_i^k e_0^l = g_{k0}(\overline{e}_i)e_i^k = 0 \quad i = 1, 2, 3$$

Следовательно, вектор $\overline{e}_0$ ортогонален вектору $\overline{e}_i$, $i = 1,\ 2,\ 3$. Согласно равенству (3.5.8)

$$F^2(\overline{e}_0) = -1$$

□

Глава 4

# Список литературы

Глава 5

# Предметный указатель